\hoffset=0.2in
\magnification=1100
\input amstex
\documentstyle{amsppt}
\nologo

%\input amstex
% \nologo
%\voffset -.1in
%\hcorrection {.2in}

%baselineskip 18pt

%\documentstyle{amsppt}
 
%\pageno=-1
%\voffset -0.2in
%NoRunningHeads
%\Monograph
\topmatter

\title
Homogeneous vector bundles on symplectic Grassmannians
\endtitle
\author Martina Bode
\endauthor
\address {Department of mathematics, Northwestern University, IL 60208, USA}
\endaddress
\email bode\@math.northwestern.edu
\endemail
\endtopmatter

\topmatter

%\title
%Contents
%\endtitle
%\endtopmatter

\vskip 1in

\document

\define\oh{{\operatorname H}}
\define\={{\  = \ }}
\define\tor{{\operatorname {Tor}}}
\define\sym{\operatorname {Sym}}
\define\st{\operatorname {S}^\perp}

%Macro for young diagrams 

\def\sqr#1#2{{\vcenter{\vbox{\hrule height.#2pt
\hbox{\vrule width.#2pt height#1 
\kern #1  \vrule width.#2pt}
\hrule height.#2pt}}}}
\def\square{\sqr{\dimen1}4}
\def\b#1{
\setbox1=\hbox{$\square$}
\countdef\I=10
\I=#1
\loop
\ifnum\I>1
\setbox0=\box1
\setbox1=\hbox{$\square$\kern-.3pt\box0}
\advance\I by -1
\repeat}
\def\ch#1,{
\if#10 \let\next=\relax
\else
\advance\count0 by1 {\ifnum\count0=1
\b{#1}
\box1\vskip\dimen0
\else
\b{#1} 
\box1\vskip\dimen0
\fi}
\let\next=\ch
\fi
\next} 

\def\yng#1{{\count0=0  \dimen0= -4.2\jot  \dimen1=8pt
\vcenter{\hbox{\vbox{\vskip -6pt\ch#1,0, }}}\, }}

\def\young#1{{\count0=0  \dimen0= -2.2\jot  \dimen1=7pt
\vcenter{\hbox{\vbox{\vskip -20pt\ch#1,0, }}}\, }}

\def\lskipamount{12pt}
\def\addline{\vskip\lskipamount plus3pt minus2pt}

\pageheight {7.2in}
\hcorrection {.2in}

\baselineskip 18pt

%\documentstyle{amsppt}
%\Monograph
%\NoRunningHeads
\topmatter
\title Introduction \endtitle
\endtopmatter

Let $V$ be a vector space of dimension $N$ over an algebraically closed field $
\Bbbk$ of characteristic 0 and denote by $Gr(k,V)$  the grassmannian of
$k$-planes in
$V$. Let
$$ J= \left( \matrix  0 & I \\ -I & 0 \endmatrix \right) \ 
$$  and denote by $X=SpGr(k,V)\subset Gr(k,V)$   the symplectic grassmannian of 
iso-tropic
 $k$~-planes with respect to the symplectic form $< \ , \ >_J $.

In this paper, we will discuss the K-theory and the category of homogeneous vector
bundles on the symplectic grassmannian $SpGr(2,N)$ of isotropic 2-planes.
Atiyah and Hirzebruch [AH], Araki [A] and Hodgkin [Ho] have shown that the
topological 
$K_0$ of homogeneous spaces is generated by homogeneous bundles.
This implies that the algebraic $K_0$ is also generated by algebraic homogeneous
bundles, however this fact relies on the topological methods of Atiyah and
Hirzebruch [AH], Araki [A] and Hodgkin [Ho]; at present there is no purely algebraic
proof.

Kapranov gives a description of the bounded derived category $D^b_{coh}$ of coherent
sheaves of
$Gr(k,V)$ in [KI], and in [L-S-W], Marc Levine, V.Srinivas, Jerzy Weyman compute the
$K$-groups of Twisted Grassmannians.
 In both articles this is done by constructing a minimal resolution of the structure
sheaf
$\Cal O_{\Delta_{Gr}}$ of the diagonal $\Delta_{Gr}$ inside of the product
$Gr(k,V)\times Gr(k,V)$. 
The explicit resolution of the diagonal gives both a purely algebraic
computation of the K-theory of the grassmannian $Gr(k,V)$ in terms of the
 representation theory
and homogeneous vector bundles as well as a description of the bounded derived
category $D^b_{coh} $ of coherent sheaves on $Gr(k,V) $. In addition, this resolution
gives a generalization of the Beilinson spectral sequence [B] to the case of coherent
sheaves on grassmannians.  We will be following these methods. Our goal is to find a
minimal resolution of the  structure sheaf of the diagonal $\Delta \subset X \times
X$.

\specialhead Resolutions for  ordinary Grassmannians
\endspecialhead

Let $S_{Gr}$ be the tautological bundle on $Gr(k,V)$ and let $S_{Gr}^\perp $ be the
dual of the quotient bundle, that is
$ ((V \otimes \Cal O_{Gr}) / S_{Gr} )^*  $ .
 Denote by $p_1$ and
$p_2$  the projections on the first respectively second factor from
$Gr(k,V) \times Gr(k,V)$ to $Gr(k,V)$, let $\Cal F$ and $\Cal G$ be vector bundles on
$Gr(k,V)$ and denote by $ \Cal F \boxtimes
\Cal G$ the external
 tensor product $p_1^* \Cal F \otimes p_2^* \Cal G$.

We have the well known Koszul complex $C_.$:
$$
 0 \to \bigwedge^{k(N-k)}(S \boxtimes S^\perp) \to \dots \to \bigwedge^2 (S
\boxtimes S^\perp) \to S \boxtimes S^\perp @>d>> \Cal O_{Gr\times Gr}
$$  
This is a resolution of $\Cal O_{\Delta_{Gr}}$,
in addition the Cauchy formula  gives
$$
 C_i = \bigwedge^i(S \boxtimes S^\perp)\cong 
\underset { \underset |\alpha|=i \to { N-k\geq
\alpha_1 \geq \alpha_2 \geq ... \geq \alpha_k \geq 0} } \to
\bigoplus
\Sigma^\alpha S  \boxtimes \Sigma^{\alpha^*}S^\perp ,
$$
where $\alpha^*$  is the conjugate Young diagram to $\alpha $, and 
 $\Sigma^\alpha S $ and $\Sigma^{\alpha^*} S^\perp $  are the   Schur functors.
Let $\Cal F $ be a locally free sheaf
 such that 
$
\oh^i (Gr,\Cal F
\otimes
\Sigma^{\alpha^*} S^{\perp})=0  $ for all $ i \geq 1 
$
 { and all}
$ 
N-k\geq
\alpha_1 \geq \alpha_2 \geq ... \geq \alpha_k \geq 0 \ ,
$
 then the Koszul complex gives rise
 to a finite resolution of $\Cal F$:
$$
\gather 0 \to \oh^0 (Gr, \Cal F \otimes \Sigma^{(k,...,k)} S^{\perp}) \otimes
\Sigma^{(N-k,...,N-k)} S  \to ...\\  \to\bigoplus_{|\alpha|=i} \oh^0 (Gr, \Cal F
\otimes \Sigma^{\alpha^*} S^{\perp}) \otimes
\Sigma^{\alpha} S \to ... \to \oh^0 (Gr, \Cal F ) \otimes\Cal O_X \to \Cal F
\to 0
\endgather 
$$

This resolution can be applied to give a computation of the algebraic K-theory as well
as to describe the structure of the derived category of the grassmannian $(Gr(k,V))$.
For  the K-theory this results in a
description of 
$K_*(Gr(k,V)) $ as a free $K_*(\Bbbk )$-module, with generators 
$$
 \Cal X=\{\, [ \Sigma^\alpha S_{Gr}]\  | \ \alpha=(\alpha_1,\alpha_2, ...
,\alpha_k)\  \& \  N-k\ge \alpha_1 \ge \alpha_2 \ge ... \ge \alpha_k \ge 0
\}
$$
and
$$ \Cal Y=\{\, [ \Sigma^\beta S_{Gr}]\  | \ \beta=(\beta_1,\beta_2, ...
,\beta_{N-k})\  \& \  k\ge \beta_1 \ge \beta_2 \ge ... \ge \beta_{N-k}
\ge 0 \}.
$$

\specialhead The symplectic Grassmannian $X$ of isotropic 2-planes in V
\endspecialhead

In this work we attempt to extend these calculations to the symplectic Grassmannian
$X=SpGr(2,V) \subset Gr(2,V)$. We will use the above Koszul complex
$C_.$ of
$\Cal O_{\Delta_{Gr}} $ and restrict it to $X \times X$,  but now $C_.\otimes \Cal
O_{X \times X}$ is no longer exact. The obstruction to the exactness results
from the torsion in degree one. We use the Tate construction to kill the extra torsion
and obtain a resolution of $\Cal O_\Delta $. Specifically, in degree two,
$$ C_2 =\overset 2 \to \bigwedge (S\boxtimes S^\perp) \cong \sym_2 S \  \boxtimes
\wedge^2 S^\perp \  \oplus \ \wedge^2 S \boxtimes \sym_2 S^\perp \ , 
$$
 we extend $ \sym_2 S^\perp$ to the non-trivial extension $\Psi $ of $ \sym_2
S^\perp$ and $\Cal O_X$. This kills the extra torsion in degree one.
Then we define non-trivial extensions $\Psi_\beta$ of
$\Sigma^\beta S $ and $\Psi_{\beta-2} $ to
 extend this construction to the whole complex. 

Let us denote the thus obtained complex by $D_.$.

\proclaim{Theorem}
 $D_.$ is a resolution of
$\Cal O_\Delta .$
\endproclaim 

(Chapter 2, Theorem 2.10)

 Now this complex does not terminate but instead becomes periodic in large degrees.
This is analogous to Kapranov's construction in the case of quadrics.
In the case of a quadric $Q$, when terminating this resolution, the last term can be
understood using the Clifford algebra of $Q$.
It turns out that this terminated
resolution has the right number of generators. 

 However, a similar construction applied to the   the resolution $D_.$ does not lead
to a minimal resolution. There are  summands in $D_.$ , which are not independent
in $K_0$.  Nevertheless, it turns out that $D_.$ contains a finite sub-complex, let us
call it
$B_.$, which seems to be a good candidate for a minimal resolution.  
The main result of this work is now to show that this sub-complex is  exact.

\proclaim{Main-Theorem} The sub-complex $B_.$ is exact.
\endproclaim 

(Chapter 3, Theorem 3.2)

 We will prove this in chapter 6 by showing that the quotient complex $D_. / B_. $ is
exact up to   degree $2N-7$, and therefore the sub-complex is exact up to degree
$2N-8$. The proof of exactness of the quotient complex   is an application of the geometric
techniques of calculating syzygies, see [K] and [PW], which we discuss in
chapter 4.

Since the codimension of the diagonal $\Delta \subset X \times X$ is equal to
the dimension of $X$, thus $codim_{X \times X} \Delta = 2N-5$, it follows that
  the Kernel $K$ of $B_{2N-6} \to B_{2N-7}$ has to be a vector bundle. This results
in  a finite resolution of $ \Cal O_\Delta$, with i-th entry
$$ B_i=\bigoplus_{|\alpha|=i} \Sigma^\alpha S \ \boxtimes \  \Psi_{\alpha^*},\  N-2
> \alpha_1 \geq \alpha_2 \ge 0 
$$ for $ i\leq m=2N-6$ and $(2N-5)$-th entry $K$.
Although we don't prove it here, the vector bundles
$\Psi_{\alpha^*},\  N-2
> \alpha_1 \geq \alpha_2 \ge 0 $ have no relations in $K_0$. The same applies to the
bundles
$\Sigma^{\alpha} S,\  N-2
> \alpha_1 \geq \alpha_2 \ge 0 $.

The hope is that we can filter the kernel $K$ as a direct sum of tensor products of
the form $\Cal F \boxtimes \Cal G $ and that these plus the terms in $B_.$ will lead
to a purely algebraic computation of the K-theory and a description of the derived
category $D^b_{coh} $ of symplectic grassmannians.

As for now we will complete the computation of the K-theory for the symplectic
grassmannian of isotropic 2-planes in 4-space, Chapter 5,  $ X=SpGr(2,4)$, and leave
the other cases for future research.

\newpage

\topmatter
\title\chapter{1} K-theory of ordinary Grassmannians \endtitle
\endtopmatter

In this chapter we will discuss the K-theory of the Grassmannians by illustrating the work of
Kapranov [KI], On the derived category of coherent sheaves on Grassmann Manifolds, and the work of
Marc Levine, V. Srinivas and Jerzy Weyman [LSW], K-theory of some twisted Grassmannians.

\head Notations
\endhead

Let $V$ be a vector space of dimension $N$ over an algebraically closed field $ \Bbbk$ of
characteristic 0, let $Gr=Gr(k,V)$ be the grassmannian of 
$k$-planes in $V$ and denote by
$\Delta_{Gr}$  the diagonal inside of $Gr(k,V) \times Gr(k,V)$. Let $S$ be the tautological bundle
of the grassmannian $Gr(k,V)$ and denote by
$S^\perp$ the dual of the quotient of the tautological bundle, that is
$$ S^\perp \cong (( V \otimes \Cal O_X )/ S )^*\ .
$$ Thus $S^\perp$ is the vector bundle that fits into the short exact sequences:
$$
\gather 0 \to S \to V \otimes \Cal O_X \to ( S^\perp)^* \to 0 \quad \text {and} \\ 0 \to S^\perp \to
V^* \otimes \Cal O_X \to S^* \to 0 .
\endgather
$$

Let $\alpha=(\alpha_1, \dots , \alpha_k) $ and
$ \beta =( \beta_1, \dots ,\beta_{N-k} )  $ be  ordered partitions, that is $\alpha_1 \geq \dots
\geq \alpha_k $ and $  \beta_1\geq \dots \geq \beta_{N-k}  $. Denote by $\alpha^* $ the
conjugate partition defined by interchanging rows and columns in the Young diagram
$\alpha$, that is if
$\alpha
$ is represented by a Young diagram with  $ \alpha_i $ boxes in the i-th row, then the Young diagram
 of the conjugate partition has $ \alpha_i $ boxes in the i-th column.

Let $\Sigma^\alpha S $ and $\Sigma^\beta S^\perp $ be the  corresponding Schur functors, see
appendix for the definition of   Young symmetrizers and the Schur functors.

Denote by $p_1$ and
$p_2$  the projections on the first respectively second factor from
$Gr(k,V) \times Gr(k,V)$ to $Gr(k,V)$, let $\Cal F$ and $\Cal G$ be vector bundles on $Gr(k,V)$ and
denote by $ \Cal F \boxtimes
\Cal G$ the external
 tensor product $p_1^* \Cal F \otimes p_2^* \Cal G$.

\addline

The main objective for the K-theory of the Grassmannian $Gr(k,V)$  is to define a minimal resolution
of the  structure sheaf $\Cal O_{\Delta_{Gr}}$ of the diagonal 
$\Delta_{Gr}\subset Gr(k,V) \times Gr(k,V)$. The terms of this resolution decompose as a sum of
external tensor products of the form
$\Sigma^\alpha S \boxtimes \Sigma^{\alpha^*} S^\perp $. 
Using a standard trick, 
 we will construct resolutions  for certain sheaves
$\Cal F$, which in turn can be used to find a generating system of the  K-theory of the
Grassmannian
$Gr(k,V)$. This leads to  a description of the K-theory of the grassmannian as a direct sum over
the generating system of the K-theory of the field $\Bbbk$, see [LSW]. See [KI] for a
description of the bounded derived category $D^b_{coh}$ of coherent sheaves of
$Gr(k,V)$.

\specialhead { Resolution of $\Cal O_{\Delta_{Gr}} $}
\endspecialhead

We will define  a homomorphism from $  S \boxtimes S^\perp$ to $\Cal O_{Gr
\times Gr} $ with cokernel $\Cal O_{\Delta_{Gr}} $ and then define the associated Koszul complex.
We will use some well known facts about Koszul complexes and demonstrate that this defines a
resolution of
$\Cal O_{\Delta_{Gr}}
$.

Define a
homomorphism of coherent sheaves
$d: S
\boxtimes S^\perp
\to
\Cal O_{Gr
\times Gr} $ as follows:
Choose a  basis $v_1, v_2 ... v_{N}$ of   
$V$ and the corresponding dual  basis $x_1,x_2 ... x_{N}$ of    $V^*$. Then on the fibers over the
point $W_1 \times W_2 \in Gr \times Gr$ define the map 
$$
\gather
 d: W_1
\boxtimes(V/W_2)^* \to \Bbbk \\
d(v\boxtimes f) =\sum_{i=1}^{N} x_i(v) \boxtimes f(v_i) = f(v).
\endgather
$$ This definition does not depend on the base choice. Note over a fiber this map is just the map
given by evaluation.

\proclaim{Remark} An alternative definition is given by Kapranov [KI],1.1, page 185:

 Consider the cohomology of the dual bundle of $S
\boxtimes S^\perp$:
$$
\oh^* (Gr \times Gr, ( S \boxtimes S^\perp )^*) \cong \oh^0(Gr, S^*) \otimes
\oh^0(Gr,(S^\perp)^*) \cong V^* \otimes V \cong \text{End}(V)
$$ Choose the section $\ s: Gr \times Gr \to (S \boxtimes S^\perp)^*$ corresponding to the identity
element in End$(V)$, such that
  the set of zeroes of $s$ is equal to the diagonal $\Delta_{Gr}$. Therefore $s^*$  gives rise to
the homomorphism $d: S \boxtimes S^\perp \to \Cal O_{Gr \times Gr}$ with cokernel
$\Cal O_{\Delta _{Gr}}$. 

Note that these two different definitions in fact define the same map, because the identity element
in $V^* \otimes V $ is given by $ \sum_{i=1}^{N} x_i\boxtimes v_i$ and
$f(v)=\sum_{i=1}^{N} x_i(v) \boxtimes f(v_i)$ for $v \in V$ and $f \in V^*$.

\endproclaim

Denote by $C_.$ the Koszul complex of the map $d$, that is the complex 
$$
\gather 0 \to \bigwedge^{k(N-k)}(S \boxtimes S^\perp) \to \dots \to \bigwedge^2 (S
\boxtimes S^\perp) \to S \boxtimes S^\perp @>d>> \Cal O_{Gr\times Gr} \\
\text{with} \\  d((v_1 \boxtimes f_1 ) \wedge \dots \wedge (v_k \boxtimes f_k)) \\ =
\underset 0 \leq i \leq k \to \Sigma (-1)^i d (v_i \boxtimes f_i) (v_1 \boxtimes f_1) \wedge
\dots \wedge
\widehat {( v_i \boxtimes f_i)}
\wedge  \dots \wedge ( v_k \boxtimes f_k) .
\endgather
$$ This is a well defined complex. For ordinary grassmannians this Koszul complex is exact as well.
The differences in the symplectic case arise  from the fact that this complex will no longer be 
exact.

\proclaim {Theorem 1.1} The Koszul complex 
$$
C_.: \ 0 \to  C_{k ( N-k)} 
\to ... \to   C_i
\to ... \to \Cal O_{X \times X} 
$$ is a resolution of $\Cal O_{\Delta_{Gr}}$.
\endproclaim

\demo{Proof} It will be enough to show that the image of $d$ is an ideal sheaf and that
$S \boxtimes S^\perp$ has the right rank. Note that the image of $d$ is nothing else but the ideal
sheaf of 
$ \Delta_{Gr}$, and  the rank of $S \boxtimes S^\perp$ is equal to
$ k (N-k)$ which is the same as the codimension of 
$\Delta_{Gr} \subset Gr \times Gr $.

First let us discuss in general terms when a Koszul complex is exact.

\enddemo

\head {Koszul complexes}
\endhead

Let $X$ be a noetherian scheme over $\Bbbk$, and $\Cal F$ a vector bundle on $X$. Given a
homomorphism
$\Psi : \Cal F \to \Cal O_X$, we can form the Koszul complex $K.(\Cal F, \Psi) $ of $\Cal F$ and $
\Psi$.

We define the Koszul complex to be the complex:
$$
\gather
 K.( \Cal F, \Psi ) :  0 \to \overset top \to \Lambda \Cal F \to ... \to 
\overset k \to \Lambda \Cal F @> d >> \overset k-1 \to \Lambda \Cal F \to ...
\Cal F @> \Psi >> \Cal O_X \quad
\\
\text {where}\quad d(f_1 \wedge f_2 \wedge ... \wedge f_k) =
\underset 0 \leq i \leq k \to \Sigma (-1)^i \Psi (f_i) f_1 \wedge ...\wedge \widehat { f_i }
\wedge  ... \wedge f_k .
\endgather
$$

This is a well defined complex.

\proclaim {Lemma 1.2} Suppose that the image of $\Psi$ is equal to the ideal sheaf $ \Cal I _Y $ of
a locally complete intersection $Y\subset X$. Further more if 
$rank(\Cal F)=codim_X Y$, then the Koszul complex $K(\Cal F ,\Psi )$ is exact.  
\endproclaim

\demo {Proof}  This is just a local question, therefore choose an open cover $\Cal U = \{ U_i
\} $ of
$X$, such that for all open sets $U \in \Cal U$:
$Y|_U $ is a complete intersection and $\Cal F|_U$ is free. Then locally,
$$ U = Spec (A), \Cal F |_U = \overset \sim \to M, M \cong A x_1 \oplus ...\oplus A x_r  
$$ for a noetherian ring $A$  and denote by $\overset \sim \to M$  the sheaf associated to $M$ on 
$U$. The map $\Psi|_U:\Cal F |_U \to \Cal O_U $ is equal to 
$(\psi: M \to A)^\sim$.

Let $\psi(x_i)=f_i$, then
$$ K(U, \Psi|_U) \cong K(A, f_1, ... ,f_r)^\sim ,
$$ which is given by: 
$$
\gather
 K.( A, f_1, ... ,f_r) :  0 \to \overset top \to \Lambda M \to ...
\to 
\overset k \to \Lambda M @> d >> \overset k-1 \to \Lambda M \to ... M @> d >> A \quad
\\
\text {where}\quad d(x_{i_1} \wedge x_{i_2} \wedge ... \wedge x_{i_r}) =
\underset 0 \leq j \leq r \to \Sigma (-1)^i f_{i_j} x_{i_1} \wedge ...\wedge \widehat
{x_{i_j}}
\wedge  ... \wedge x_{i_r}.
\endgather
$$

This is exact if and only if $f_1, .. , f_r$ form a regular sequence, see [H II,7.10A] . Let $I_Y$
be the ideal  generated by
$f_1, ... ,f_r$ , then since $Y|_U$ is a complete intersection,
 the $f_i$'s form a regular sequence [H] II, 8.21 A(c). Thus the Koszul complex of $f_1, ...,f_r$ is
exact. q.e.d.
\enddemo

\demo{Proof of Theorem 1.1} Since the image of $d$ is equal to the ideal sheaf of 
$ \Delta_{Gr}$, and  the rank of $S \boxtimes S^\perp$ is equal to
 the codimension of 
$\Delta_{Gr} \subset Gr \times Gr $, the exactness follows from Lemma 1.2.

q.e.d.
\enddemo

\subhead {The Cauchy Formula}
\endsubhead

 The i-th term of the Koszul complex
$$
C_.: \ 0 \to  C_{k ( N-k)} 
\to ... \to   C_i
\to ... \to \Cal O_{X \times X} 
$$  can be expressed as a direct sum of external tensor
products of vector bundles over the grassmannian $Gr(k,V)$. We'll do this by applying the Cauchy
formula to this complex, (char $ \Bbbk=0$):
$$ C_i = \bigwedge^i(S \boxtimes S^\perp)\cong \bigoplus_{|\alpha| = i}
\Sigma^\alpha S  \boxtimes \Sigma^{\alpha^*}S^\perp \  
$$
{  (see [M] I,4,Ex.6 and [KI], Lemma 0.5)},
where
$$
\alpha\in I= \{ \alpha=(\alpha_1,\alpha_2, ...,\alpha_k) \  | \qquad N-k\geq
\alpha_1 \geq \alpha_2 \geq ... \geq \alpha_k \geq 0 \ \}  
$$ runs over all Young diagrams with i-cells, and 
$\alpha^*$ is the conjugate Young diagram to $\alpha $.

Set 
$$
 C_\alpha = \Sigma^\alpha S  \boxtimes
\Sigma^{\alpha^*}S^\perp.
$$
This implies:

\proclaim{Corollary 1.3}
The Koszul complex $C_.$ is isomorphic to the complex
$$ 0 \to  C_{(k,...,k)} 
\to ... \to  \bigoplus_{|\alpha| = i} C_\alpha
\to ... \to \Cal O_{X \times X} \ .
$$

\endproclaim

We will apply this resolution to the K-theory of the grassmannian $Gr(k,V)$.

\specialhead{K-theory of $Gr(k,V)$ }
\endspecialhead

 Let $ \Cal P$ be the category of vector bundles on Spec $\Bbbk $. Following  Quillen's proof for
the K-theory of projective spaces  [Q], define 
$$
 U_{\alpha} : \Cal P \to \Cal C_{Gr} , \ U_{\alpha} (W) = q^* W \otimes
\Sigma^\alpha S 
$$  for all  Young diagrams $\alpha \in I$.  Now $U_\alpha$ is inducing a homomorphism $u_\alpha$ on
the K-theory, that is
$$
 u_\alpha :  K_*( \Bbbk ) \to K_*(Gr(k,V)) 
$$   Set $u= \underset {\alpha \in I } \to \bigoplus \  u_\alpha $, then

\proclaim{Theorem 1.4}
 
$$ u: \underset \alpha \in I \to \bigoplus K_*( \Bbbk ) \to K_*(Gr(k,V)) \quad
\text {is an isomorphism.}
$$ 
\endproclaim

\demo{Proof} We will prove this in several steps. 

Define
$$  I= \{ \alpha=(\alpha_1,\alpha_2, ...,\alpha_k) \  | \qquad N-k\geq
\alpha_1 \geq \alpha_2 \geq ... \geq \alpha_k \geq 0 \ \}.
$$

Let
 $ \Cal C_{Gr}$ be the category of locally free coherent sheaves on $Gr(k,V)$ and define  $ \Cal C$ 
as the full subcategory of  the category $ \Cal C_{Gr}$ as  follows:
$$
 \Cal C = \left\{ {  \Cal F \in \Cal C_{Gr} \biggm| \ \ \aligned &\oh^i (Gr,\Cal F
\otimes
\Sigma^{\alpha^*} S^{\perp})=0 \ \  \text {for all}\  i \geq 1 \\ &\text { and all}\ 
 \alpha \in I \endaligned } \right\} .
$$ 

First we prove that u is surjective by
 using a standard trick  to construct resolutions for all  sheaves
$\Cal F \in \Cal C$ from the above resolution of 
$\Cal O_{\Delta_{Gr}} $ . 

Then we will use cohomology computations in order to show that u is injective.
\enddemo

\head u surjective:
\endhead

We will apply the standard method to construct resolutions for all sheaves  $\Cal F \in \Cal C$. We
will tensor the resolution of $\Cal O_{\Delta_{Gr}}$ with the pullback of $\Cal F$ to the product of
$Gr
\times Gr$ on one factor and then push this new complex forward on the other factor.

\subhead {The pushforward of a complex of coherent sheaves}
\endsubhead

\proclaim {Lemma 1.5} Let $q: X \to Y$ be a projective morphism of noetherian schemes. Let 
$ 0 \to \Cal F_m @> d_m >> \to ... @> d_1 >> \Cal F_0 
$ be an exact sequence of coherent $\Cal O_X$~-modules.  Suppose $R^i q_* \Cal F_j =0 $ for all
$i>0$ and all $j$, then
$$ 0 \to q_* \Cal F_m \to ... \to q_* \Cal F_0
$$ is an exact complex of  coherent $\Cal O_Y$~-modules.
\endproclaim

\demo {Proof}
First we  split the complex $\Cal F.$ into short exact sequences,
$$ 0 \to ker \ d_{m-k} \to \Cal F_{m-k} \to Im \ d_{m-k} \to 0 . $$
In order to show that $q_* \Cal F.$ is exact, we also split this complex into short sequences,
$ 0 \to ker \ q_* d_{m-k} \to q_* \Cal F_{m-k} \to q_* Im \ d_{m-k} \to 0
$, and show that  these are exact for all $k$. The sequence 
$ 0 \to ker \ q_* d_{m-k} \to q_* \Cal F_{m-k} \to q_* Im \ d_{m-k} \to 0
$ is exact if $R^1 q_* ker \ d_{m-k} =0$, or since $ ker \ d_{m-k} \cong Im \ d_{m-(k-1)}$, this is
exact if $R^1 q_* Im \ d_{m-(k-1)} =0$

We  show by induction over $k$ that 
$ R^i q_* Im \ d_{m-(k-1)} =0$ for all $ 0 \leq k \leq m$ and all $i>0$.
\item {$k=0:$} The first case of $k=0$ is obvious, because $ Im \ d_{m+1} =0 $ and therefore 

$R^i q_* Im \ d_{m+1} =0 $.
\item{$k \rightsquigarrow k+1:$} Now suppose $R^i q_* Im \ d_{m-(k-1)} =0$ for all $i>0$, then we
need to show that

$R^i q_* Im \ d_{m-k} =0$. 

Consider the short exact sequence:
$$
 0 \to Im \ d_{m-(k-1)} \to \Cal F_{m-k} \to Im \ d_{m-k} \to 0 \ ,
$$ and consider the resulting long  exact sequence: 
 $$
 ... \to R^i q_* \Cal F_{m-k} \to R^i q_* Im \ d_{m-k} \to R^{i+1} q_* Im
\ d_{m-(k-1)} \to ... \ .
$$

From the assumptions we know that $R^i q_* \Cal F_{m-k}=0$ and from the induction hypothesis it
follows  that $R^{i+1} q_* Im \ d_{m-(k-1)}=0$, therefore
$R^i q_* Im \ d_{m-k}=0$. q.e.d.
\enddemo

\specialhead {Resolution for $\Cal F \in \Cal C$}
\endspecialhead

Let $\Cal F \in \Cal C$ be a locally free sheaf, that is
$
 \Cal F \in \Cal C_{Gr} $
 such that 
$$
\oh^i (Gr,\Cal F
\otimes
\Sigma^{\alpha^*} S^{\perp})=0 \ \  \text {for all}\  i \geq 1 \text { and all}\ 
 \alpha \in I \ .
$$
Recall that $p_1$ and
$p_2$ are the projections
$$
 p_i : Gr(k,V) \times  Gr(k,V) \to  Gr(k,V)
$$ to the i-th factor. Denote by $q$  the structure map
$q:  Gr(k,V) \to Spec(\Bbbk).$ Then
$$ R^i p_{1*}(p_2^* \Cal F \otimes \Cal O_{\Delta_{Gr}}) \cong 
\left\{ \aligned  & \Cal F \ \ 
\text {for} \ i=0,  \\ &  0  \ \ \text{ otherwise }\endaligned \right.
$$ and 
$$\aligned  R^i p_{1*}(p_2^* \Cal F \otimes C_\alpha)  & \cong R^i p_{1*}(p_2^* \Cal F \otimes p_2^*
\Sigma^{\alpha^*} S^\perp  \otimes p_1^*\Sigma^{\alpha}S)
 \\ & \cong  R^i p_{1*} ( p_2^*( \Cal F \otimes  \Sigma^{\alpha^*} S^\perp )
\otimes p_1^*\Sigma^{\alpha}S) 
 \\ & \cong  R^i p_{1*} ( p_2^*( \Cal F \otimes  \Sigma^{\alpha^*} S^\perp ))
\otimes  \Sigma^{\alpha}S \ \ \text { (by projection formula)}
\\ & \cong q^* R^i q_* ( \Cal F \otimes  \Sigma^{\alpha^*} S^\perp )
\otimes  \Sigma^{\alpha}S \ \ \ \ \  (q \ \text {flat} )
\\  & \cong  H^i  ( \Cal F \otimes  \Sigma^{\alpha^*} S^\perp )
\otimes  \Sigma^{\alpha}S, 
\endaligned
$$ which in turn implies:

\proclaim {Proposition 1.6}

For $\Cal F \in \Cal C$ , $R^i p_{1*} ( p_2^* \Cal F \otimes C_. ) $ gives rise to a finite
resolution of $\Cal F$:
$$
\gather 0 \to \oh^0 (Gr, \Cal F \otimes \Sigma^{(k,...,k)} S^{\perp}) \otimes
\Sigma^{(N-k,...,N-k)} S  \to ...\\  \to\bigoplus_{|\alpha|=i} \oh^0 (Gr, \Cal F
\otimes \Sigma^{\alpha^*} S^{\perp}) \otimes
\Sigma^{\alpha} S \to ... \to \oh^0 (Gr, \Cal F ) \otimes\Cal O_X \to \Cal F
\to 0
\endgather 
$$

\endproclaim

\demo {Proof} For $\Cal F \in \Cal C , H^i  ( \Cal F \otimes  \Sigma^{\alpha^*} S^\perp )=0 $, thus
all the higher direct images vanish. The statement follows from Lemma 1.4. q.e.d.

\enddemo

We will show next that it suffices to define $K_*(Gr(k,V))$ on the subcategory
$\Cal C $ rather than defining it on the whole category $\Cal C_{Gr} $.

\subhead {K$_* (\Cal C) \cong $ K$ _*(Gr)$}
\endsubhead

Let $g:Y \to Z$ be a quasi-projective map of schemes, and let $\Cal X$ be a set of locally free
sheaves on $Y$. Denote by $\Cal C_Y$ the category of locally free sheaves $\Cal F$ on $Y$ and let
$\Cal C_{\Cal X}$  be 
 the subcategory of $\Cal C_Y$ consisting of all sheaves $\Cal F$ with
$$ R^i g_* (\Cal F \otimes \Cal G ) = 0 \ \text {for all} \ i > 0 
\ \text {and for all} \ \Cal G \in \Cal X \ .
$$

\proclaim{Lemma 1.7} The inclusion $\Cal C_{\Cal X} \subset \Cal C_Y$ induces an isomorphism on the
K-theory, that is
$$ K_* (\Cal C_{\Cal X}) \cong K_*(\Cal C_Y) \ .
$$

\endproclaim

\demo{Proof} This follows from  [LSW] Lemma 4.3, which implies that 
$q$ induces a homotopy equivalence on the classifying spaces
$BQ \Cal C_{\Cal X} \to BQ \Cal C_Y$. q.e.d.

\enddemo

\subhead{u surjective}
\endsubhead

 Since $K_*(Gr) \cong K_* (\Cal C)$ by Lemma 1.7, it is enough to show that
$u$ is surjective on $K_* (\Cal C)$. Let $\Cal F \in \Cal C$.
 Consider the finite resolution of $\Cal F $. Denote by
$\lbrack {\Cal F } \rbrack$ the class of $\Cal F $ in  $ K_*(\Cal C)$. Then
$$
 \lbrack {\Cal F } \rbrack = \sum_ {0 \neq   \alpha \in I } (-1)^{|\alpha|}
\lbrack \oh^0(\Cal F \otimes \Sigma^{\alpha ^*} S^\perp ) \otimes
\Sigma^\alpha S \rbrack,
$$ thus $u$ is surjective. 

q.e.d.

\head{u injective}
\endhead

To show that $ u $ is injective we  define a homomorphism 
$$
 v: K_*(Gr(k,V)) \to \underset \alpha \in I \to \bigoplus K_*( \Bbbk ),
$$ and then show   that $v
\circ u$ is injective, hence that $u$ is an isomorphism.

Let
$$
 \Cal C^* = \left\{ {  \Cal F \in \Cal C_{Gr} \biggm| \ \ \aligned &\oh^i (Gr,\Cal F
\otimes
\Sigma^{\alpha} S^*)=0 \ \  \text {for all}\  i \geq 1 \\ &\text { and all}\ 
 \alpha \in I \endaligned } \right\} .
$$

By Lemma 1.7, it is enough to define the K-theory on the subcategory $\Cal C^*$.

Set
$$ V_\alpha : \Cal C^* \to \Cal P , \ \ V_\alpha(\Cal F ) = q_*(\Cal F \otimes
\Sigma^\alpha S^* ) ,
$$
 for all  $\alpha \in I$. Note that  $V_\alpha$ on the  category
$\Cal C_{Gr} $ is not exact, because $R^i q_* (\Cal F \otimes \Sigma^\alpha S^* ) $ is in general
not zero for $ i>0$, but $V_\alpha$ is exact on the subcategory $\Cal C^*$. The $V_\alpha$'s  induce
homomorphisms
$v_\alpha$  on the K-theory.
$$  v_\alpha :  K_*(Gr(k,V)) \to  K_*( \Bbbk ).
$$   Let $ v=(...,v_\alpha,...)$. Suppose that the higher cohomology of 
$\Sigma ^\alpha S \otimes ( \Sigma ^\beta S )^* $ vanishes and 
consider  $v \circ u $ :
$$ V_\beta \circ U_\alpha (W) =V_\beta ( q^* W \otimes \Sigma^\alpha S ) = \lbrack
\oh^0 ( W \otimes \Sigma ^\alpha S \otimes ( \Sigma ^\beta S )^* ) \rbrack \ .
$$

We will show that the matrix of $v \circ u $ with respect to the lexicographical ordering of the
weights $\alpha$, is upper triangular with  ones down the diagonal, thereby showing that $u$ is an
isomorphism.

First we will recall some cohomology computations. Since we will need the same techniques for the
symplectic Grassmannians we will do this  more detailed then needed at this point.

\head {Cohomology Computations}
\endhead

Let $G$ be a Lie group, $B$ a Borel subgroup containing the torus $T$ of $G$. We denote by $R
\subset X$ the set of all roots of $G$ with respect to $T$ and call negative the roots of $B$. Let
$\rho$ be the half sum over all positive roots $\alpha$, i.e. $\rho =  {1 \over 2} \underset {\alpha
> 0} \to \Sigma
\alpha =  \Sigma \lambda_i $, where the last sum is over the fundamental weights $\lambda_i$  .Let
$\lambda$ be a character of $T$, then the one dimensional representation $ L_{\lambda}$ of $T$,
extends trivially to $B$, and we can form the holomorphic line bundle 
$$
\Cal L_{\lambda} = G \times_B L_{\lambda}= \left\{ {(g,v) \in G \times L_{\lambda} | (g,v) \sim
(gb,b^{-1} v),  b \in B} \right \}.
$$ We have Bott's theorem for the vanishing of the cohomology of this line bundle:
   
\proclaim {Borel-Weil-Bott  Theorem 1.8}
\roster
\item  If $\lambda $ is a dominant weight, then $\Cal L_{\lambda}$ has no higher cohomology and the
space of sections \ $\oh^0(G/B, \Cal L_{\lambda}) $ \ is the dual of the irreducible representation
with highest weight $\lambda$.
\item  If $\lambda$ is not dominant add $\rho$ to $\lambda$, then find \
$\sigma$ in the Weyl group such that $\sigma ( \lambda + \rho) $ is dominant, then subtract again $
\rho$. 

If $\sigma ( \lambda + \rho) - \rho $ is dominant, then there is only cohomology in degree $
l(\sigma)=lenght(\sigma)$ and $\oh^{l(\sigma)} ( G/B,
 \Cal L_{\lambda}) \cong  \oh^0 ( G/B, \Cal L_{\sigma ( \lambda +
\rho) - \rho})$ which by (1) is given by the dual of the irreducible representation with highest
weight $\sigma ( \lambda + \rho) - \rho $.

If $\sigma ( \lambda + \rho) - \rho $  is not dominant, then $\oh^{i} ( G/B,\Cal L_{\lambda}) =0 $ \
for all $ i \ge 0$. 
\endroster

\endproclaim
\demo{Proof} See \cite{D}.
\enddemo

Here we will use $G=Gl_N(V)$ for our cohomology computations for $Gr(k,N)$. Let
$F=G/B$ be the total flag manifold, i.e. 
$$ F=\left\{ 0 \subset W_1 \subset W_2 \subset ... \subset W_N = V \ \text { with dimension} \ W_i =
i  \right \} 
$$ Then $Gr(k,N) $ is isomorphic to $G/P$, where $P \supset B$ is the parabolic subgroup preserving
the space $ V_k $  in the standard representation. $ P$ can be described as the parabolic subgroup,
corresponding to omitting one node of the Dynkin diagram:( see also [F-H]  \$ 23.3 Homogeneous
spaces)

$\circ$---$\circ$---$\circ$---$\circ$---$\circ$---$\bullet$---$\circ$---$\circ$---$\circ$---$\circ$

From this description it is easy to get the Levi decomposition of $P$ into Levi factors
$$ P \cong (L_1 \times L_2) \ltimes U_P
$$ where $U_P$ is the unipotent radical of $P$,\, $L_1 \cong GL(V_k)$, and $L_2
\cong GL(V/V_k)$ are the Levi factors with $V_k$ the vector space of dimension
$k$ of the standard representation. We denote by $\pi$ the projection $\pi:G/B
\to G/P$. Now the vector bundles $\Sigma^\alpha S$ and $\Sigma^\beta S^\perp$ on
$Gr(k,V)$ come from line bundles on $F =G/B$. Set 
$$\aligned &\lambda(\alpha) = (-\alpha_k,-\alpha_{k-1},
\dots ,-\alpha_1,0, \dots,0) \\ &\lambda(\beta) = (0, \dots,0,\beta_1,
\dots,\beta_{N-k-1},\beta_{N-k}).
\endaligned
$$

\proclaim{Lemma 1.9} 
\roster
\item \qquad $\pi_*\Cal L_{\lambda(\alpha)} \cong \Sigma^\alpha S$
\item \qquad $\pi_*\Cal L_{\lambda(\beta)} \cong \Sigma^\beta S^\perp$
\item \qquad $\pi_*\Cal L_{\lambda(\alpha) + \lambda(\beta)} \cong 
\Sigma^\alpha S \otimes \Sigma^\beta S^\perp$
\endroster
\endproclaim
\demo{Proof} See \cite{KI}  2.5.
\enddemo We order the partitions $\alpha =(\alpha_1, \dots,\alpha_m)$ lexicographically, that is
$\alpha \leq \beta$ if $\alpha_1 < \beta_1$, or $\alpha_1 = \beta_1$ but
$\alpha_2 < \beta_2$, etc.

\proclaim{Lemma 1.10} For $\alpha =(\alpha_1, \dots,\alpha_k)$ and $\beta = (\beta_1, \dots,
\beta_k)$ with $N-k \geq \alpha_1 \geq \dots \geq \alpha_k \geq 0$ and $N-k \geq \beta_1 \geq \dots
\geq \beta_k \geq 0$:
\roster
\item \qquad $\oh^i(Gr,(\Sigma^\beta S)^*\otimes \Sigma^\alpha S) = 0 \ \
$ for all $i>0$
\item \qquad $\oh^0(Gr,(\Sigma^\beta S)^*\otimes \Sigma^\alpha S) = 0 \ \
$ for $\alpha > \beta$
\item \qquad $\oh^0(Gr,(\Sigma^\alpha S)^*\otimes \Sigma^\alpha S) = \Bbbk $
\endroster
\endproclaim

\demo{Proof} See \cite{KI} 2.2.

Since we will be using these techniques again, we will outline Kapranov's proof.

Note that $(\Sigma^\beta S )^* \cong \Sigma^{(-\beta_k, \dots , -\beta_1)} S$. Consider the tensor
product:
$$
\Sigma^\alpha S \otimes (\Sigma^\beta S )^* \cong
\Sigma^\alpha S \otimes \Sigma^{(-\beta_k, \dots , -\beta_1)} S 
\cong \bigoplus \Sigma^\gamma S
$$ Since all the $\alpha_i$'s and $\beta_i$'s are between $0$ and $N-k$, it follows that  every
summand $\gamma$ occurring in this product satisfies $ -(N-k) \leq \gamma_i \leq N-k$, see appendix,
Remark (3).

Consider the cohomology of 
$$
\Sigma^\gamma S \cong \pi_* (\Cal O_F (-\gamma_k, \dots , -\gamma_1, 0 \dots 0))
$$ Set $\overset \wedge \to \gamma=(-\gamma_k, \dots , -\gamma_1, 0 \dots 0)$.
\item{(a)} If $\gamma_1 \leq 0$, thus $-\gamma_k \geq \dots \geq -\gamma_1 \geq 0 = \dots = 0 $, then
$\overset \wedge \to \gamma $ is a dominant weight and by the Borel-Weil-Bott Theorem has only $H^0
\cong \Sigma^{(\gamma_1, \dots , \gamma_k)} V^*$.
\item{(b)} If $\gamma_1 > 0 $, then $\overset \wedge \to \gamma $ is not dominant. Consider
$$
\overset \wedge \to \gamma + \rho = (N-\gamma_k, \dots , N-k+1-\gamma_1, N-k, \dots ,1) \ .
$$ Since $\gamma_1 > 0 $, it follows that $N-k+1-\gamma_1$ is between $1$ and $N-k$, thus $\overset
\wedge \to \gamma + \rho$ has a repetition and therefore
$\Cal O_F (\overset \wedge \to \gamma) $ has no cohomology.

\item{(1)} From (a) and (b) it follows that for all summands $\gamma$, $\Sigma^\gamma S$ can have at
most $H^0$.
\item{(2)} If $\alpha > \beta$, then $\gamma_1 > 0$, thus $\Sigma^\gamma S$ has no cohomology.
\item{(3)} If $\alpha = \beta$, consider the tensor product
$$
\Sigma^\alpha S \otimes (\Sigma^\alpha S )^* \cong \bigoplus \Sigma^\gamma S \ ,
$$ then $\gamma$ is either $(0, \dots , 0)$ or at least $\gamma_1 > 0 $. Thus 
$$ H^0(Gr,\Sigma^\alpha S \otimes (\Sigma^\alpha S )^* ) \cong H^0(Gr, \Cal O_{Gr})\cong \Bbbk \ .
$$ q.e.d.
\enddemo

\proclaim{Lemma 1.11} For $\alpha =(\alpha_1, \dots,\alpha_{N-k})$ and
$\beta = (\beta_1, \dots, \beta_{N-k})$ with $k \geq \alpha_1 \geq \dots
\geq
\alpha_{N-k}
\geq 0$ and $k \geq \beta_1 \geq \dots \geq \beta_{N-k} \geq 0$ we have
\roster
\item \qquad $\oh^i(Gr,(\Sigma^\beta S^\perp)^*\otimes \Sigma^\alpha S^\perp) = 0 \ \
$ for all $i>0$
\item \qquad $\oh^0(Gr,(\Sigma^\beta S^\perp)^*\otimes \Sigma^\alpha S^\perp) = 0 \ \
$ for $\alpha > \beta$
\item \qquad $\oh^0(Gr,(\Sigma^\alpha S^\perp)^*\otimes \Sigma^\alpha S^\perp) =
\Bbbk $
\endroster
\endproclaim
\demo{Proof} See \cite{LSW}  Lemma 4.2.
\enddemo

\subhead {u injective}
\endsubhead

Now we finish the proof of Theorem 1.3.
 Consider the matrix of $ v \circ u $ with respect to the lexicographical ordering of the
weights $\alpha$. Then Lemma 1.9 implies:
$$
\gather
\oh^i ( W \otimes \Sigma ^\alpha S \otimes ( \Sigma ^\beta S )^* )= 0 
\ \text {for} \ i > 0 \ \text{and} 
\\
  V_\beta \circ U_\alpha (W)  =
\lbrack
\oh^0 ( W \otimes \Sigma ^\alpha S \otimes ( \Sigma ^\beta S )^* ) \rbrack =
\left\{ { 
\aligned 
 & 0 \  \ \text { for } \alpha > \beta \\
 & \Bbbk \ \ \text { for }  \alpha = \beta 
\endaligned } \right. \ .
\endgather 
$$  Thus the matrix of $v \circ u $ is upper triangular with ones down the diagonal, therefore $ u $
is injective, hence an isomorphism.

\newpage

\topmatter
\title\chapter{2} The Tate construction \endtitle
\endtopmatter

Let $V$ be an $N=2n$ dimensional vector space over an algebraically closed
field
$\Bbbk$ of characteristic 0. In this section $X$ shall denote the symplectic
grassmannian
$SpGr(2,V)$ of two dimensional  symplectic planes and $Gr(k,V)$ shall denote
the grassmannian  of 
$k$-planes in $V$. In the previous chapter we discussed the definition of the
Koszul complex
$C.$ of the structure sheaf 
$\Cal O_{\Delta_{Gr}}$ of the diagonal $\Delta_{Gr}$ inside the product
$ Gr(k,V) \times Gr(k,V)$. Using the inclusion $ X \subset Gr(2,V)$ and 
the complex
$C.$  restricted to $X \times X$, we will construct  a  resolution of
the structure sheaf 
$ \Cal O_\Delta$ of the diagonal $\Delta$ inside  $X \times X$.

Consider the Koszul complex $C.$ restricted to $X \times X$. In degree $i$ the
terms are:
$$ C_i \otimes \Cal O_{X \times X} = \bigwedge^i(S \boxtimes S^\perp)\cong
\bigoplus_{|\alpha| = i}
\Sigma^\alpha S  \boxtimes \Sigma^{\alpha^*}S^\perp \
$$

This  Koszul complex is no longer exact. In this section
 we will use  Tate's techniques of
 adjoining an element in degree two to kill the extra torsion in degree one.

Accordingly we will define a non-trivial  extension $\Psi$ of $Sym_2 \st$
and
$\Cal O_X$ and then extend this construction to the whole complex, that is we
will extend 
$\overset i \to \Lambda (S \boxtimes \st) $ to vector bundles $\Phi_i$. This
extended complex defines  a resolution of $\Cal O_\Delta$.

Analogous to the case of the ordinary grassmannians, the
bundles $\Phi_i$ break up into  direct sums over  external tensor
products of homogeneous bundles, that is
$\Phi_i \cong \bigoplus \Sigma^{\alpha}S \boxtimes \Psi_{\alpha^*}$ for some homogeneous
bundles
 $\Psi_{\alpha^*}$ on $X$.

\specialhead {Notations and Preliminaries}
\endspecialhead

Let $G=Sp_N(V)$ be the symplectic group of $V$ consisting of all 
$x \in Gl_N (V) $ satisfying
$$ {}^\tau x J x = J \ ,
$$

$$
\text {where} \ J=\left(
\matrix 0 &I \\ -I &0
\endmatrix
\right) \ .
$$ Let $< \ , \ >_J $ be the symplectic form defined by
$<v,w>_J= {}^\tau v J w$ . We call a subspace $W$ of $V$ isotropic if it is
isotropic with respect to the form $< \ , \ >_J $. 

Let $X=SpGr(2,V)$ be the symplectic grassmannian of isotropic two planes in
$V$ inside the ordinary grassmannian $Gr(2,V)$ of two planes in $V$. Let us
fix a symplectic basis
$v_1, v_2, \dots , v_{2n} $ of $V$, i.e.
$$ <v_i,v_j>_J =0 \ \text {if} \ j \neq i \pm n \  \text{and} \ <v_i,v_i+n>_J
=1
\
\text{for} \ 1 \leq i \leq n \ .
$$  Let $x_1, x_2, \dots , x_{2n} $
 be the corresponding dual basis of $V^*$. Set
$$ z=x_1 \wedge x_{n+1} + x_2 \wedge x_{n+2} + \dots + x_n \wedge x_{2n} \ .
$$

\proclaim {Lemma 2.1}
The symplectic grassmannian 
$X$ is given by the hyperplane section $z$ inside the grassmannian $Gr(2,V)$.
\endproclaim

\demo{Proof} Consider the Pl\"ucker embedding that  embeds $Gr(2,V)$ in the
projective space 
$\Bbb P( \overset 2 \to \bigwedge V)$. The condition for an element 
$ v \wedge w$ in $Gr(2,V)$ to be isotropic is given by
$ <v,w>_J=0$. This is equivalent to  $z( v \wedge w ) =0$.  Thus $z$ is the
defining equation of $X$.

Let us show that $z$ is independent of the base choice:

Let $\Cal A= \left( \matrix A & B \\ C & D  \endmatrix \right) = (s_{i j})_{1
\leq i,j \leq 2n }
$ be a symplectic matrix and $y'=\Cal A \cdot y$, then
$$
\aligned 
& \ y_1 \wedge y_{n+1} + y_2 \wedge y_{n+2} + ... + y_n \wedge y_{2n} \\ &=
\underset 1
\leq j,k \leq 2n \to \sum \ \underset 1
\leq i \leq n \to \sum s_{i \ j} y'_j \wedge s_{n+i \ k } y'_k \\ &=
\underset 1 \leq j<k
\leq 2n \to \sum\  \underset 1
\leq i \leq n \to \sum ( s_{i \ j} s_{n+i \ k } - s_{i \ k} s_{n+i \ j } ) y'_j
\wedge  y'_k 
\endaligned 
$$ 
\item{}
If $1 \leq j<k \leq n$, then the coefficient of $y'_j \wedge y'_k $is
given by 
$( ^\tau A C -^\tau C A )_{j \ k } = 0$,
\item{}
 if $ n+1 \leq j < k \leq 2n$, then
the coefficient equals 
$( ^\tau B D -^\tau D B )_{j-n \ k-n } = 0$ and finally
 if 
$1 \leq j \leq n , n+1 \leq k \leq 2n $, then the coefficient is 
$( ^\tau A D -^\tau C B )_{j \ k-n } = \delta_{j \ k-n}$, hence 
$$ y_1 \wedge y_{n+1} + y_2 \wedge y_{n+2} + ... + y_n \wedge y_{2n}=y'_1
\wedge y'_{n+1} + y'_2 \wedge y'_{n+2} + ... + y'_n \wedge y'_{2n},
$$
 thus z is well defined. 

q.e.d.

\enddemo

Since $X$ is a hyperplane in $Gr(2,V)$,  we obtain the short exact sequence:
$$ 0 \to \Cal O_{Gr} (-1)  @>  z >> \Cal O_{Gr} \to \Cal O_X \to 0 ,
$$ which we will use  for cohomology computations.

For our purposes it will be convenient to give a description of $X$ as
homogeneous space. Let $V_i = Span \{v_1, \dots , v_i\} $ and denote by 
$0 \subset V_1 \subset V_2 \subset \dots \subset V_n \subset V $ the standard
symplectic flag. Let $B$ be the Borel preserving the  standard symplectic flag
and let $P$ be the parabolic subgroup preserving the point $V_2 \in X$. Then
$$  X \cong G/P \ . 
$$

Let $\Cal V $ be a homogeneous vector bundle on $X$. Denote by
$\Cal V_{eP}\cong \Cal V_{V_2} $ the zero fiber of $\Cal V$ over the point
$V_2 \in X$. Recall that homogeneous vector bundles are completely determined
by its 
 representation of $P $ on the zero fiber $\Cal V_{eP}\cong \Cal V_{V_2} $, that is
$
\Cal V \cong G \times_P \Cal V_{eP}
$.  Conversely each representation $W$ of $P $ defines a homogeneous vector bundle
$$ G \times_P W = G \times W /  (gp, w) \sim (g,pw) , p \in P  \ .
$$

Let $S$ be the tautological bundle on $X$  and $\st$ the dual of the quotient
bundle on
$X$, i.e. the bundle that fits into the  short exact sequence:
$$ 0 \to S^\perp \to V^* \otimes \Cal O_X \to S^* \to 0
$$
Analogously denote by $S_{Gr}$ and $\st_{Gr}$ the corresponding bundles on the
grassmannian
$Gr(2,V)$.

Following the notations of the previous chapter, let
$\alpha=(\alpha_1, \alpha_2) $ and 
$\beta=(\beta_1, \dots , \beta _{N-2} ) $ be ordered partitions, and let 
$\Sigma^\alpha S,\  \Sigma^\beta S^\perp,\   \Sigma^\alpha S_{Gr}  $ and 
$\Sigma^\beta \st_{Gr}  $ be the associated Schur functors. Denote by
$\alpha^* $ the conjugate or dual of 
$\alpha $, with $\alpha_j= $ number of entries in the j-th column of the Young
diagram given by $\alpha $.

Note that all these bundles are homogeneous:
$$
\gather S \cong G \times_P V_2 \quad S^\perp \cong G \times_P (V/ V_2)^* \\
\Sigma^\alpha S \cong G \times_P \Sigma^\alpha V_2 \quad 
\Sigma^\beta S^\perp \cong G \times_P \Sigma^\beta (V/ V_2)^* \ . 
\endgather
$$

$$
\gather
\text{Let} \ \  C_\alpha = \Sigma^\alpha S \boxtimes \Sigma^{\alpha^*} \st 
\ \text{for} \ N-2 \geq \alpha_1 \geq \alpha_2 \geq 0,
\\ C_i = \overset i \to
\Lambda (S \boxtimes \st) \cong 
\underset |\alpha|=i \to \bigoplus C_\alpha \ \ 
 \text {and} \ \  C= \underset i \geq 0 \to \bigoplus C_i \ .
\endgather
$$

\specialhead The Tate construction $\Phi$ in degree 2
\endspecialhead

Let us first discuss the Tate construction in degree 2. We will
 adjoin an element in degree 2 in order to  kill the extra torsion in degree
1, that is we will define an extension $\Psi \subset V^* \otimes \st$ of
$Sym_2
\st$ and 
$\Cal O_X$. Then 
$$
\gather
\Phi= \overset 2 \to \Lambda S \boxtimes \Psi \ \oplus \ Sym_2 S \boxtimes
\overset 2 
\to \Lambda S^\perp 
\ \text{ is an extension of } \ \
\\
 C_2=\overset 2 \to \Lambda (S \boxtimes \st) \cong 
\  
\overset 2 \to \Lambda S \boxtimes Sym_2 \st \ \oplus
\ Sym_2 S \boxtimes \overset 2 
\to \Lambda S^\perp
\endgather 
$$
 in the Koszul complex $C_.$. At the end of this chapter  we will 
complete the discussion of $\Psi$ and show that
$\Psi$ is a non-trivial extension.

Let
$$
\eta=(x_1 \otimes x_{n+1} + x_2 \otimes x_{n+2} + \dots + x_n \otimes x_{2n}
) \in V^* \otimes (V/V_2)^*, 
$$ and define $W$ to be the subspace of $V^* \otimes (V/V_2)^*$,
$$ W = Sym_2 (V/V_2)^* \oplus \Bbbk \eta \ .
$$ Recall that the natural action of $P$ on $V^*$ is given by:
$$ p \cdot f(v) = f ( p^{-1} v ) \ .
$$  The corresponding action on $V^* \otimes V^*$ leaves $Sym_2 (V/V_2)^*  $
invariant, which is defining the vector bundle $Sym_2 S^\perp \cong G \times_P
Sym_2 (V/V_2)^*  $. Suppose that this action is leaving $W$ invariant as well,
then the
$P-$ representation $W$ defines the homogeneous vector bundles:

\proclaim {Definition 2.2}
$$
\Psi = G \times_P W \ \subset V^* \otimes \st  \text{and} \ 
\Phi = \overset 2 \to \Lambda S \boxtimes \Psi \oplus \ Sym_2 S \boxtimes
\overset 2 
\to \Lambda S^\perp \ .
$$

\endproclaim

These are well defined as long as $W$ is a  representation of $P $, which we
will show next.

Let $\frak {g}$ be the Lie algebra of the Lie group $G$ and let
$\frak p$ be the parabolic algebra corresponding to the parabolic group $P
\subset G$. The action of $P$ on $V^*$ translates to the natural action of
$\frak p$ on $V^*$, that is  the action of $\frak p$ on $V^*$ is given by 
 $p\cdot f(v) = f(-^\tau p v )$.

\proclaim {Lemma 2.3}
\roster
\item The Levi  factors of the parabolic $P$ are
$ L_1 \cong Gl_2 \ \text{and} \ L_2 \cong Sp_{2n-4} \ .
$

\item For
$ p \in \frak p, p \cdot \eta \in Sym_2 (V/V_2)^*$ and  

\item
$P \to Gl(W)$ is a   subrepresentation of $P $ of $V^* \otimes
(V/V_2)^*$.

\endroster

\endproclaim

\demo{Proof} It is enough to do these computations on the algebra level.

\item{(1)} Let us first study the Levi factors of $P$ and $\frak p$. Consider
the Levi decomposition of $P$   into Levi factors:
$$ P \cong (L_1 \times L_2) \ltimes U_P ,
$$ where $L_1$ and $L_2$ are the Levi factors and $U_P$ is the unipotent
radical. This composition corresponds to the Levi decomposition of $\frak p$
as 
$$
\frak p \cong \frak l_1 \oplus \frak l_2 \oplus \frak u_{\frak p} .
$$ Observe that the condition for an element 
$ x =
\left(
\matrix m &n \\ p &q
\endmatrix
\right)
 \in \frak {gl}_{2n}$ to be symplectic is that
$$
\left(
\matrix 0 &I \\ -I &0
\endmatrix
\right) \
\left(
\matrix m &n \\ p &q
\endmatrix
\right)  =
\left(
\matrix -^\tau m &-^\tau p \\ -^\tau n &-^\tau q
\endmatrix
\right) \
\left(
\matrix 0 &I \\ -I &0
\endmatrix
\right),
$$ i.e. that
$ {}^\tau n = n, {}^\tau p = p,$ and $ {}^\tau m = -q$.

Thus for $ p \in \frak p$, this Levi decomposition corresponds to:
$$
\gather
 p= \left( 
\matrix
 A & F && B & G \\
 0 & C && {}^\tau G & D \\ \\
 0 & 0 && -^\tau A & 0 \\
 0 & E && -^\tau F & -^\tau C
\endmatrix
\right) 
\\ =
\left( 
\matrix  A & 0 & 0 & 0 \\  0 & 0 & 0 & 0 \\   0 & 0 & -^\tau A & 0 \\  0 & 0
& 0 & 0
\endmatrix
\right) +
\left( 
\matrix  0 & 0 & 0 & 0 \\  0 & C & 0 & D \\   0 & 0 & 0 & 0 \\  0 & E & 0 &
-^\tau C
\endmatrix
\right)    +
\left( 
\matrix  0 & F & B & G \\  0 & 0 & {}^\tau G & 0 \\   0 & 0 & 0 & 0 \\  0 & 0
& -^\tau F & 0
\endmatrix
\right).
\endgather
$$ where
$A$ and $B$ are   $2 \times 2$ matrices,  $F$ and $G$ are $ \ 2 \times (n-2)$
matrices and 
$C,D$ and
$E$ are 
$(n-2) \times (n-2)$ matrices. Since $ p$ is a symplectic matrix, $B,D$ and
$E$ are symmetric, i.e. $B={}^\tau B$,
$D={}^\tau D$ and $E={}^\tau E$. Therefore the Levi factors of the parabolic
algebra
$\frak p$ are 
$$
\frak l_1 \cong \frak {gl}_2 \ \text {and} \ \frak l_2 \cong \frak
{sp}_{2n-4} 
$$ and  the  Levi factors of the parabolic group $P$ are 
$$ L_1 \cong Gl_2 \ \text{ and} \ L_2 \cong Sp_{2n-4}.
$$

\item{(2)} The action of $\frak p $ on $V^* \otimes V^*$ is given by:
$$ p \cdot ( f \otimes g ) = (p \cdot f) \otimes g + f \otimes (p \cdot g),  
$$ and the action of $\frak p$ on  $f= \Sigma a_i x_i=\left(
\matrix a_1 \\ \vdots \\ a_{2n}
\endmatrix
\right)$: 
$$ (p \cdot f) = - {}\tau  p 
\left(
\matrix a_1  \\ \vdots \\ a_{2n}
\endmatrix
\right) .
$$ 

We will use the Levi decomposition of $\frak p$ into the two levi factors
$\frak l_1$ and
$\frak l_2$ and the unipotent ideal
$\frak u_{{\frak p}}$.
Let
$p \in \frak p$ be as above with 
$p=p_1 + p_2 + u$. Then 
$$p \cdot \eta = 
 p_1 \cdot \eta + p_2 \cdot \eta + u \cdot \eta
\ .
$$ All the following matrix
multiplications are the obvious ones, for example:
$$
\gather {}^\tau A
\left( 
\matrix x_1 \\ x_2 \endmatrix \right) \otimes
 \left( 
\matrix x_{n+1} \\ x_{n+2} \endmatrix \right) \\ = (a_{1,1} x_1 + a_{2,1} x_2
)
\otimes x_{n+1} + (a_{1,2} x_1 + a_{2,2} x_2 ) \otimes x_{n+2}
\endgather
$$ Using the symmetry of $B,D$ and $E$, we get:
 
$$
\align p_1 \cdot \eta  &= 
\left( 
\matrix  A & 0 & 0 & 0 \\  0 & 0 & 0 & 0 \\   0 & 0 & -^\tau A & 0 \\  0 & 0
& 0 & 0
\endmatrix
\right) \cdot (x_1 \otimes x_{n+1} + x_2 \otimes x_{n+2} + \dots + x_n \otimes
x_{2n} ) \\ &=- {}^\tau A  \left( 
\matrix x_1 \\ x_2 \endmatrix \right) 
 \otimes \left( \matrix  x_{n+1} \\ x_{n+2} \endmatrix \right)  + 
\left( 
\matrix x_1 \\ x_2 \endmatrix \right) 
 \otimes A \left( 
\matrix x_{n+1} \\ x_{n+2} \endmatrix \right)  \\ & = 0
\endalign
$$

$$
\align p_2 \cdot \eta &= 
\left( 
\matrix  0 & 0 & 0 & 0 \\  0 & C & 0 & D \\   0 & 0 & 0 & 0 \\  0 & E & 0 &
-^\tau C
\endmatrix
\right) \cdot (x_1 \otimes x_{n+1} + x_2 \otimes x_{n+2} + \dots + x_n \otimes
x_{2n} ) \\ &= -{}^\tau C  \left( 
\matrix x_3 \\ \vdots \\ x_n \endmatrix \right) 
 \otimes \left( \matrix  x_{n+3} \\ \vdots \\ x_{2n} \endmatrix \right) + 
\left( 
\matrix x_3 \\ \vdots \\ x_n \endmatrix \right) 
 \otimes C \left( 
\matrix x_{n+3} \\ \vdots \\ x_{2n} \endmatrix \right) \\ &- {}^\tau D 
\left( 
\matrix x_3 \\ \vdots \\ x_n \endmatrix \right) 
 \otimes \left( \matrix  x_{n+3} \\ \vdots \\ x_{2n} \endmatrix \right)  - 
\left( 
\matrix x_3 \\ \vdots \\ x_n \endmatrix \right) 
 \otimes {}^\tau E \left( 
\matrix x_{n+3} \\ \vdots \\ x_{2n} \endmatrix \right) \\ &= 0 - \underset 2
\leq i<j
\leq n \to \Sigma (d_{i,j}+e_{i,j})  ( x_{n+j} \otimes x_{n+i} + x_{n+i}
\otimes x_{n+j} ) \\ & \in Sym_2 (V/V_2)^*
\endalign
$$
$$
\align & u \cdot \eta \\ &=
\left( 
\matrix & 0 & F & B & G \\ & 0 & 0 & {}^\tau G & 0 \\  & 0 & 0 & 0 & 0 
\\ & 0 & 0 & -^\tau F & 0
\endmatrix
\right) \cdot x_1 \otimes (x_{n+1} + x_2 \otimes x_{n+2} + \dots + x_n \otimes
x_{2n} ) \\ &= -{}^\tau F 
\left( 
\matrix x_1 \\ x_2 \endmatrix \right) 
 \otimes \left( \matrix  x_{n+1} \\ x_{n+2} \endmatrix \right)  + 
\left( 
\matrix x_1 \\ x_2 \endmatrix \right) 
 \otimes F \left( 
\matrix x_{n+1} \\ x_{n+2} \endmatrix \right) 
\\ &-  {}^\tau G 
\left( 
\matrix x_1 \\ x_2 \endmatrix \right) 
 \otimes \left( \matrix  x_{n+1} \\ x_{n+2} \endmatrix \right) - G
\left( 
\matrix x_3 \\ \vdots \\ x_n \endmatrix \right) 
 \otimes \left( \matrix  x_{n+3} \\ \vdots \\ x_{2n} \endmatrix \right) 
\endalign
$$
$$
\align
 &-  {}^\tau B 
\left( 
\matrix x_1 \\ x_2 \endmatrix \right) 
 \otimes \left( \matrix  x_{n+1} \\ x_{n+2} \endmatrix \right) \\ &= 0 \\ &-
\underset 0
\leq i \leq n-2, 1 \leq j \leq 2 \to \Sigma g_{i,j} (x_{n+(i+2)} \otimes
x_{n+j} + x_{n+j} \otimes x_{n+(i+2)}) \\ &- 
\underset 1 \leq i<j \leq 2 \to \Sigma b_{i,j}  ( x_{n+j} \otimes x_{n+i} +
x_{n+i}
\otimes x_{n+j} ) \\ & \in Sym_2 (V/V_2)^*
\endalign
$$ This implies that $p \cdot \eta $ is in $Sym_2(V/V_2)^* $, and thus that 
$P \to Gl ( W) $ is an indecomposable  representation of $P$.

\item{(3)} Since $\frak p$ leaves the subspace $(V/V_2)^*$ and thus the space
$Sym_2 (V/V_2)^* $ invariant, we only need to check the action of
$\frak p$ on $ \Bbbk \eta$. Since
$ p \cdot \eta $ is in $Sym_2(V/V_2)^* \subset W $ for $ p \in \frak p$, it
follows that
$W$ is a representation of $P$.

q.e.d.
\enddemo 
\proclaim {Remark}

\roster
\item
$\Psi$ is a subbundle of the homogeneous bundle $V^* \otimes \st$, this
follows from Lemma 2.3(3).

\item Since $X$ is a homogeneous space isomorphic to $G/P$, where $P \supset
B$ is the parabolic subgroup preserving the point $ V_2 \in X $, 
 $ P$ can also be described as the parabolic subgroup, corresponding to
omitting one node of the Dynkin diagram:
 
( see also [FH]  \$ 23.3 Homogeneous spaces)

$\circ$---$\bullet$---$\circ$---$\circ$---$\circ$---$\circ$---$\circ$---$\circ
\Longleftarrow
\circ$

Thus the Levi factors of  $P$ correspond to the two Dynkin diagrams:

$\circ \ $ and $\ \circ$---$\circ$---$\circ$---$\circ$---$\circ$---$\circ
\Longleftarrow
\circ$ , therefore 
$ L_1 \cong GL_2 \ \text{and} \  L_2 \cong Sp_{2n-4} \ .
$

\endroster

\endproclaim

\specialhead The Tate construction in higher degrees
\endspecialhead

We will extend the construction of $\Psi$ and $\Phi$ to the whole complex.
We will do this in several steps, let us outline the
procedure.

Let $\beta$ be an ordered partition 
 with $ 2 \geq \beta_1 \geq \dots \geq \beta_{N-2} \geq 0 
$. 
 
In the first part  we define the homogeneous bundles 
$\Psi_\beta = G \times_P W_\beta$ for some   representations  $W_\beta$ of $P$.

In the second part we discuss the homogeneous
$\Cal O_{X
\times X}$-bundles
$$\Phi_k= ( G \times G) \times_{P \times P} M_k . 
$$ These will be extensions of
$\overset k \to \Lambda (S \boxtimes \st)$ in the Koszul complex and are
needed to extend the Koszul complex to a resolution of $\Cal
O_\Delta$.

In the third part we will relate  these two  definitions.  We will show that 
$$
\Phi_k \cong \underset |\beta|=k \to \bigoplus 
\Sigma^{\beta^*} S \boxtimes \Psi_\beta \ .
$$

Finally in the last part, we extend the map $d$ to the Tate construction
$$
D_. = \underset k \leq 0 \to \bigoplus \Phi_k
$$ and show that this 
 complex defines a resolution of the structure sheaf of the diagonal.

\head {The extensions $\Psi_\beta$}
\endhead

We will now define the extensions $\Psi_\beta$ of $\Sigma^\beta \st$ 
 {for} $2 \geq \beta_1 \geq
\dots \geq \beta_{N-2} \geq 0$ .

Denote by $\beta-2$  the ordered partition
$$
\gather
\beta -2 = \left \{ 
\matrix & ( \beta_2, \beta_3, .. ,\beta_{N-2}, 0 ) && \text {if}  \ \beta_1=2
\\ & 0 && \text {otherwise}
\endmatrix \right.  \ .
\endgather
$$
Set $
 \Sigma_\beta =  \Sigma^{\beta} (V/V_2)^* $
  and let
$Sym_2=Sym_2 (V/V_2)^*$.

\proclaim{Definition 2.4}

Define $W_\beta$ as
$$
\gather W_\beta=\underset  i \geq 0 
\to 
\bigoplus \Sigma_{\beta -2i}  \otimes \eta^i 
\subset \bigoplus \Sigma_{\beta -2i}  \otimes Sym_i( V^* \otimes (V/V_2)^* )
\ .
\endgather
$$
\endproclaim

Note that this definition is extending the definition of $W$ in the previous
section, since
$$
W_2= \Sigma_{(2,0, \dots 0)} \oplus \Sigma_{(0, \dots 0)} \otimes \eta
\cong Sym_2 (V/V_2)^* \oplus \Bbbk \eta = W.
$$

We will define an action on $W_\beta$. 
We cannot use the natural action of
$P$ on $W_\beta$ as subspace of $\bigoplus \Sigma_{\beta -2i}  \otimes Sym_i(
V^* \otimes (V/V_2)^* )$, since this action does not leave $W_\beta$
invariant, for example:

Let $ \beta=(2,1,0, \dots , 0),$ then
$ W_\beta = \Sigma^{(2,1,0, \dots , 0)} (V/V_2)^* \oplus (V/V_2)^* \otimes
\eta$. Choose $ p \in \frak p$ and $ f \in (V/V_2)^*$, then
$$
\gather
 p \cdot ( f \otimes \eta ) = ( p \cdot f ) \otimes \eta + f \otimes
( p \cdot \eta).
\endgather
$$
Since $p \cdot \eta \in Sym_2 (V/V_2)^*  $, see Lemma 2.3, it follows that 
$$
f \otimes ( p \cdot \eta) \in (V/V_2)^* \otimes Sym_2 (V/V_2)^*  
$$
which is isomorphic to $\Sigma^{(2,1,0, \dots , 0)} (V/V_2)^* \oplus 
Sym_3 (V/V_2)^*
$, see appendix on the tensor product of Young diagrams.
Here the tensor product is associated to the tensor product
$$
\gather
\young{1} \boxtimes \young{2} \cong \yng{2,1} \oplus \young{3} \ .
\\
\endgather
$$
Note that
$$
f \otimes ( p \cdot \eta ) \in
\Sigma^{(2,1,0, \dots , 0)} (V/V_2)^* \oplus  Sym_3 (V/V_2)^*
$$
 is  not
contained in 
$W_{(2,1,0, \dots , 0)}$, but the projection of $f \otimes ( p \cdot
\eta) $ in
$$
\Sigma^{(2,1,0, \dots , 0)} (V/V_2)^* \oplus  Sym_3 (V/V_2)^* \to
\Sigma^{(2,1,0, \dots , 0)} (V/V_2)^* 
$$ is.
In order to define an action of $P$ on $W_\beta$, we need to define 
a projection  
$$
\pi_\beta: \ 
\Sigma_{\beta-2} \otimes Sym_2 \to
 \Sigma_\beta 
$$
in general.

\subhead {The $*$-multiplication}
\endsubhead

Let us first discuss the tensor product of 
$\Sigma_{\beta-2}$ and $Sym_2$.

\proclaim{Lemma 2.5}
Let $\beta$ be an ordered partition $N-2 \geq \beta_1 \geq \dots \geq
\beta_{N-2} \geq 0$ such that 
$\Sigma_\beta$ and $\Sigma_{\beta-2}$ are non-zero. 

Then $\Sigma_\beta$ is a summand of $\Sigma_{\beta-2} \otimes Sym_2$ with
multiplicity one.

\endproclaim

\demo{Proof}

Consider the corresponding tensor product of the Young diagrams:
$$
\beta-2 \  \otimes \young{2}  
$$
This is given by adding
the two boxes of $\young{2}$ to $\beta-2$
according to the Littlewood-Richardson  rule, (see appendix), that is at most
one box to each column. This adds either   no,  one or   two boxes to
the first two columns of the Young diagram
$\beta-2$.
If both boxes are added to the first two columns, then there is only one way
of doing so and the resulting weight is equal to
$\beta$. Thus $\beta$ appears with multiplicity one.

q.e.d.

\enddemo

Since $\Sigma_\beta$ is a direct summand of 
$\Sigma_{\beta-2} \otimes Sym_2$ with
multiplicity one, we can define the projections
$$
\pi_\beta: \ 
\Sigma_{\beta-2} \otimes Sym_2 \to
 \Sigma_\beta \ .
$$
These are well defined and determined up to a constant.
Note that these  projections project $\Sigma_{\beta-2} \otimes Sym_2$ onto
the factors 
$\Sigma_\beta$ for all partitions $\beta$ that do appear in
the Koszul complex.

\proclaim{Definition 2.6}
Let $\ N-2 \geq \alpha_1 \geq \alpha_2 \geq 0$.

\item {(1)}{(a)}
If $ \alpha^* -2$ is not zero, then let $\pi_\alpha$ be the
projection
$$
\pi_\alpha: \Sigma_{\alpha^*-2}  \otimes Sym_2
\to
\Sigma_{\alpha^*} \ ,
$$
\item {}{(b)}
 if $ \alpha^* -2$ is  zero, then set $\pi_\alpha=0$.
\item{(2)}
 {Set} 
$
\pi=
\Sigma \pi_\alpha : 
\bigoplus \Sigma_{\alpha^*-2}  \otimes Sym_2
\to \bigoplus
\Sigma_{\alpha^*} \ ,
$ where the sum goes over all 
$ \ N-2 \geq \alpha_1 \geq \alpha_2 \geq 0 $ and

\item{(3)}
for $F \in \bigoplus \Sigma_{\alpha^*-2}$, define the $*$-multiplication:
$$
 F *  ( f \cdot g)= 
\pi (F \otimes (f \cdot g ) ) \ .
$$

\endproclaim

\demo{$\pi$ is well defined }

$\pi$ is well defined  since 
$ \Sigma_\alpha$ appears as the summand  of the tensor product
$\Sigma_{\alpha^*-2}  \otimes Sym_2$   with multiplicity one, Lemma 2.5,
which determines the projections $\pi_ {\alpha^*}$ up to multiplication by a
constant.

\enddemo

We will see later that the choice of the scalars does not effect our
computations.
We will use the $*$-multiplication to define the Lie algebra action on
$W_\beta$.

\subhead {The  action of $P$ on $W_\beta$}
\endsubhead

\proclaim {Definition 2.7}

For $p \in \frak p, f_\beta \cdot \eta^i \in W_\beta$ let
$$
\gather p \cdot \eta^i = i ( p \cdot \eta ) \cdot \eta^{i-1} , \tag 1 \\ p
\cdot (f_\beta \cdot
\eta^i )  = (p \cdot f_\beta) \cdot \eta^i +  i (f_\beta * ( p
\cdot \eta ) )  \cdot \eta^{i-1} \ , i \geq 1 \tag 2 \\
p \cdot f_\beta = p \cdot f_\beta \ , f_\beta \in W_\beta \ .
\endgather
$$

\endproclaim 

\demo{The  action of $P$ is well defined}
 
This is a well defined action leaving $W_\beta$ invariant. However we need
to show that this is a Lie algebra action.
$$
\gather
\text { Claim:} \\
 \pi_\beta \ \text{ commutes with the }  \text{ action of}\ P \ \text{and} \tag a
\\
(f_\beta * ( p
\cdot \eta )) * ( g
\cdot \eta ) = (f_\beta * ( g
\cdot \eta )) * ( p
\cdot \eta ) \tag b
\endgather
$$
\item{(a)}
$\Sigma_{\beta-2} \otimes Sym_2 \cong \Sigma_\beta \ \oplus \ $ other
representations of $(V/V_2)^*$. This is an isomorphism of 
representations of $Gl((V/V_2)^*)$. Since the  action of $P$ leaves $V_2$ invariant, this is also
an isomorphism of  representations of $P$.
\item{(b)}
Consider the tensor product of 
$\Sigma_{\gamma-4}  $ by $ \Sigma_{(2,2,0, \dots , 0 )}$, for
$2 \geq \gamma_1 \geq \dots \geq \gamma_{N-2} \geq 0$.
This is given by adding the four boxes of the Young diagram
$$
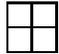
\gather
\yng{2,2} \ \ \text{to the diagram of } \ \gamma-4   , 
\ \text{ according to the Littlewood-Richardson rule.}
\endgather
$$
Similar to the proof of the previous Lemma, there is precisely one way of
adding two boxes to each of the two columns of $\gamma-4$,
thus $\Sigma_\gamma$ appears with multiplicity one in the tensor product
of 
$\Sigma_{\gamma-4}  $ by $ \Sigma_{(2,2,0, \dots , 0 )}$.

Denote by $\overset \wedge \to \pi$ the projection
$$
\overset \wedge \to \pi: 
\Sigma_{\gamma-4} \otimes \Sigma_{(2,2,0, \dots , 0 )}
\to \Sigma_\gamma \ . 
$$
This is determined uniquely up to a non-zero constant.
Note that $\Sigma_\gamma$ also appears with multiplicity one in 
$\Sigma_{\gamma-2} \otimes Sym_2$. Therefore the two projections 
$$
\gather
\pi_\gamma \circ ( \pi_{\gamma-2} \otimes id|_{Sym_2} ) \ , \
\overset \wedge \to \pi \circ ( id|_{\Sigma_{\gamma-4}} \otimes \pi_{(2,2,0
\dots 0)} ) :
\Sigma_{\gamma-4} \otimes Sym_2 \otimes Sym_2 \to \Sigma_\gamma
\endgather
$$
only differ by a constant $c$. Let $f \in \Sigma_{\gamma-4}$ and $ x,y \in
Sym_2$, then it is enough to prove that
$$
\gather
x * y = y * x \\
\text {since} \ 
( f * x ) * y =
\pi ( \pi ( f \otimes x ) \otimes y )
= c \cdot \overset \wedge \to \pi ( f \otimes   ( x * y ) \\
\text {and} \
( f * y ) * x
= \pi ( \pi ( f \otimes y ) \otimes x ) 
= c \cdot \overset \wedge \to \pi ( f \otimes  ( y *
x ) \ .
\endgather
$$

Consider the projection 
$$
\pi_{(2,0 \dots 0)} : Sym_2 \otimes Sym_2 \to\Sigma_{(2,2,0,\dots,0)} \ .
$$
The Young symmetrizer of (2,2,0,\dots,0) is equal to
$$
\gather
c_{(2,2,0, \dots ,0)} 
= [(1+e_{(1 \ 2)} ) (1+e_{(3 \ 4)} )][(1-e_{(1 \ 3)} )(1-e_{(2 \ 4)} )] \\
=[(1+e_{(1 \ 2)} ) (1+e_{(3 \ 4)} )][ 1-e_{(1 \ 3)}   -e_{(2 \ 4)} +
e_{(1 \ 3)(2 \ 4)}] ,
\endgather
$$
thus we get an explicit description of the projection
$
\pi_{(2,0 \dots 0)}$ up to a multiplication by a non-zero constant $k$ as
follows:
$$
\gather
\pi_{(2,0 \dots 0)} ( (f_1 \cdot f_2) \otimes (g_1 \cdot g_2 )) \\
= k
[(f_1 \cdot f_2) \otimes (g_1 \cdot g_2 ) 
- (g_1 \cdot f_2) \otimes (f_1 \cdot g_2 ) \\ 
- (f_1 \cdot g_2) \otimes (g_1 \cdot f_2 ) 
+ (g_1 \cdot g_2) \otimes (f_1 \cdot f_2 )] \\
=\pi_{(2,0 \dots 0)} ( (g_1 \cdot g_2) \otimes (f_1 \cdot f_2 )) .
\endgather
$$
This proves the claim. q.e.d.

For $ g,p \in \frak P $, 
denote by $[p,g]= pg-gp$ the  Lie algebra bracket. Consider
$$
\gather
p \cdot ( g \cdot ( f_\beta \otimes \eta^i )) - 
g \cdot ( p \cdot ( f_\beta \otimes \eta^i )) \\
=
p \cdot (
( g \cdot f_\beta) \otimes \eta^i 
+   (f_\beta * ( g
\cdot \eta ) )  \otimes i \eta^{i-1} ) \\ -
g \cdot (
( p \cdot f_\beta) \otimes \eta^i +   (f_\beta * ( p
\cdot \eta ) )  \otimes i \eta^{i-1} ) \\
=
(p 
 g \cdot f_\beta) \otimes \eta^i + (g \cdot f_\beta) * (p \cdot \eta)
\otimes i \eta^{i-1}
\\
+ p \cdot (f_\beta * ( g
\cdot \eta ) )  \otimes i \eta^{i-1}  
+ (f_\beta * ( g
\cdot \eta ) * ( p
\cdot \eta ) ) \otimes i (i-1) \eta^{i-2} \\
- 
(g 
 p \cdot f_\beta) \otimes \eta^i - (p \cdot f_\beta) * (g \cdot \eta)
\otimes i \eta^{i-1}
\\
- g \cdot (f_\beta * ( p
\cdot \eta ) )  \otimes i \eta^{i-1}  
- (f_\beta * ( p
\cdot \eta ) * ( g
\cdot \eta ) ) \otimes i (i-1) \eta^{i-2} \\
\text{ using (a) and (b) above}
\\
=([p,g] \cdot f_\beta ) \otimes \eta^i  
+ " f_\beta \otimes [p,g] \cdot \eta^i \  "
= [p,g] \cdot ( f_\beta \otimes \eta^i )
\endgather
$$
Thus the $*$-multiplication does define a Lie algebra action.

q.e.d.
\enddemo

\proclaim{Remark} Note that the choice of the scalars for the projections
$\pi_\beta$ define isomorphisms of representations of $P$. We can realize this in
two ways.

Let $W_\beta$   and $\overset \wedge \to W_\beta$  be representations of $P$ that
differ by  the choices of the projections $\pi_\beta$.
Suppose $\overset \wedge \to \pi_\beta = \lambda_\beta \cdot \pi_\beta $.
\roster
\item
The representation
$  P \to Gl (W_\beta) $ and $ P \to Gl (\overset \wedge \to W_\beta) $ are defined
by
$$
\rho, \overset \wedge \to \rho :P \to  \underset i \geq 0 \to \bigoplus 
Hom ( \Sigma_{\beta - 2i}, \Sigma_{\beta - 2i} ) \ \oplus \
Hom ( \Sigma_{\beta - 2i}, \Sigma_{\beta - 2(i-1)} )  \ .
$$
The difference between $\rho$ and $\overset \wedge \to \rho $ is given by the
relation:
$$
\overset \wedge \to \rho = 
\underset i \geq 0 \to \bigoplus 
(id|_{Hom ( \Sigma_{\beta - 2i}, \Sigma_{\beta - 2i} )} \oplus \lambda_{\beta -
2i} \cdot id|_{Hom ( \Sigma_{\beta - 2i}, \Sigma_{\beta - 2(i-1)} ) } ) \circ \rho
\ .
$$
\item
We can define an isomorphism between $W_\beta$ and $\overset \wedge \to W_\beta $
directly. Set
$$
\mu_\gamma= \underset j \geq 0 \to \Pi \lambda_{\gamma-2j} \ .
$$
Define $i:W_\beta  \to \overset \wedge \to W_\beta $ by
$$
i ( f_{\beta-2i}\ \eta^i ) = \mu_{\beta-2i} \cdot f_{\beta-2i}\ \eta^i \ .
$$
Then
$$
\gather
 p \cdot i ( f_{\beta-2i}\ \eta^i ) \\
= \mu_{\beta-2i} p \cdot f_{\beta-2i}\ \eta^i \\
=\mu_{\beta-2i} (  p \cdot f_{\beta-2i}) \ \eta^i +
\mu_{\beta-2i} (f_{\beta-2i} \overset \wedge \to * p \cdot \eta ) \ \eta^{i-1} \\
=\mu_{\beta-2i} (  p \cdot f_{\beta-2i}) \ \eta^i +
\mu_{\beta-2i} \lambda_{\beta - 2i + 2} (f_{\beta-2i} * p \cdot \eta ) \
\eta^{i-1} \\
=\mu_{\beta-2i} (  p \cdot f_{\beta-2i}) \ \eta^i +
\mu_{\beta-2i + 2}  (f_{\beta-2i} * p \cdot \eta ) \
\eta^{i-1} \\
=i ( p \cdot   ( f_{\beta-2i}\ \eta^i ) ) \ .
\endgather
$$
Thus $i$ defines an isomorphism of  representations of $P$ between
$W_\beta$ and $\overset \wedge \to W_\beta$.

\endroster
\endproclaim

Let us define the homogeneous bundles associated to the   representations $W_\beta$ of $P$
, that is

\proclaim{Definition 2.8}
$$
\Psi_\beta = G \times_P W_\beta
$$
\endproclaim

Next we define the homogeneous bundles $\Phi_k$ and then show how this
definition relates to the definition of the bundles $\Psi_\beta$.

\head
The homogeneous bundles $\Phi_k$
\endhead

 We will define 
$P \times P$ representations $M_k$ and and use these to define the
homogeneous bundles
$\Phi_k$.

\subhead{Notations}
\endsubhead
 Recall that
$
\eta= x_1 \otimes x_{n+1} + \dots x_n \otimes x_{2n} \in V^* \otimes
(V/V_2)^* .
$ Let $\eta^i \in Sym_i (V^* \otimes (V/V_2)^*) $ and set
$$
\gather 
\Lambda_\alpha = \Sigma^\alpha V_2 \boxtimes \Sigma^{\alpha^*}
(V/V_2)^* \ ,
\\
\Lambda_i =
\overset i \to
\Lambda (V_2 \boxtimes (V/V_2)^*) 
\cong  
\underset |\alpha|=i \to \bigoplus \Lambda_\alpha \ \text{and} \\ 
\Lambda= \underset i \geq 0 \to
\bigoplus \Lambda_i \ ,
\\
\theta^i = ( \overset 2 \to \Lambda V_2 )^{\otimes i} \boxtimes \eta^i 
\subset ( \overset 2 \to \Lambda V_2 )^{\otimes i} \boxtimes Sym_i (V^*
\otimes
(V/V_2)^*) \ \ 
\text{and} \\ M=\underset i \geq 0 \to \bigoplus \Lambda [-2i] \otimes
\theta^i
\ .
\endgather
$$

\proclaim{Remark}
\roster
\item
 There is a natural  action of $P \times P$ on $M$ coming from the inclusion of
$M$ in the  representation of $P \times P$:
$$ N= \underset i \geq 0 \to \bigoplus \Lambda [-2i] \otimes  (( \overset 2
\to
\Lambda V_2 )^{\otimes i} \boxtimes Sym_i (V^* \otimes (V/V_2)^*) ) \ .
$$
\item
$N$ has a natural multiplication induced by the wedge product on 
$\Lambda$ and the symmetric product on $Sym (V^* \otimes (V/V_2)^*) $,
that is
 for
$$
\gather
\lambda_1,\lambda_2 \in \Lambda, \mu_1 \in Sym_i (V^* \otimes (V/V_2)^*),
  \mu_2 \in Sym_j (V^* \otimes (V/V_2)^*),
\\
(\lambda_1 \otimes ((v_1 \wedge v_2)^{\otimes i}  \boxtimes \mu_1 )) \cdot 
(\lambda_2 \otimes ((v_1 \wedge v_2)^{\otimes j}  \boxtimes\mu_2 )
) \\ = (\lambda_1 \cdot \lambda_2 ) \otimes ((v_1 \wedge v_2)^{\otimes i+j}) \boxtimes (\mu_1 \cdot
\mu_2)) \ .
\endgather
$$
This restricts to give a product on $M$.

The multiplicative structure of $\Lambda$ is essential for the Tate
construction. The Tate construction only applies to skew-commutative, graded
$R$-algebras over a commutative noetherian ring $R$ satisfying certain
conditions.
\endroster
\endproclaim

\proclaim{Lemma 2.9} The  $P\times P$ action leaves 
$M$  and 
$M_k$,  the $k$-th graded piece of $M$, that is
$$ M_k=\underset i \geq 0 \to \bigoplus \Lambda_{k-2i} \otimes \theta^i 
$$  invariant.

\endproclaim

\demo{Proof}
 It is sufficient to do these computations on the algebra level. Let 

$\lambda_k \in \Lambda_k, v_1 \wedge v_2 
\in \overset 2 \to
\Lambda V_2, g,h \in \frak p $ and $i \geq 1 $, then
$$ w=\lambda_k \otimes ((v_1 \wedge v_2 )^{\otimes i} \boxtimes \eta^i) \in 
\Lambda_k \otimes \theta^i
$$ and
$$
\gather  (g,h) \cdot w \\  =((g,h) \cdot \lambda_k) \otimes ((v_1 \wedge v_2
)^{\otimes i}
\boxtimes \eta^i )  +  \lambda_k \cdot (g,h) \cdot ((v_1 \wedge v_2
)^{\otimes i}
\boxtimes \eta^i) \\  =((g,h) \cdot \lambda_k) \otimes ((v_1 \wedge v_2
)^{\otimes i}
\boxtimes \eta^i)
  + \lambda_k \cdot ((g \cdot (v_1 \wedge v_2 )^{\otimes i}) \boxtimes
\eta^i) \\ + \lambda_k \cdot  ((v_1 \wedge v_2 )^{\otimes i} \boxtimes h \cdot
\eta^i)
\\  =((g,h) \cdot \lambda_k) \otimes ((v_1 \wedge v_2 )^{\otimes i} \boxtimes
\eta^i ) +
\lambda_k \cdot ((g \cdot (v_1 \wedge v_2 )^{\otimes i}) \boxtimes \eta^i) \\
+ i
\lambda_k \cdot  ((v_1 \wedge v_2 ) \boxtimes h \cdot \eta ) \otimes 
  ((v_1 \wedge v_2 )^{\otimes i-1} \boxtimes  \eta^{i-1})
\endgather
$$ Note that $(v_1 \wedge v_2 ) \boxtimes h \cdot \eta \in 
\overset 2 \to \Lambda V_2 \boxtimes Sym_2 (V/V_2)^* \subset \Lambda_2$, see
Lemma 2.3. Therefore 
$$ i \lambda_k \cdot  ((v_1 \wedge v_2 ) \boxtimes h \cdot \eta ) \in
\Lambda_k
\cdot \Lambda_2 \subset \Lambda_{k+2}
$$ and thus for $i\geq1$ 
$$ (g,h) \cdot w \in (\Lambda_{k} \otimes \theta^i) \ \oplus \  (\Lambda_{k+2}
\otimes
\theta^{i-1})
$$ and for $w\in \Lambda_k, (g,h) \cdot w \in \Lambda_k $.

This proves that $M$ and $M_k, k\geq 0 $ are   invariant under the action of $P \times P$ .

q.e.d. 
\enddemo

Since the $M_k$'s are  representations of $P \times P $, we can define the
corresponding homogeneous vector bundles.

\proclaim{Definition 2.10} Set 
$$
\Phi_k=(G \times G)_{P \times P} M_k
$$ and set
$$ D_. = \underset k \geq 0 \to \bigoplus  \Phi_k \ .
$$
\endproclaim

The bundles $\Phi_k$ are the extensions needed to extend the Koszul complex
$C_.$ to an exact complex. Before extending the map $d$ to $D_.$, we
show that these bundles are isomorphic to a direct sum of homogeneous
bundles, that is 
$$ \Psi_k \cong 
\underset  {\underset
 N-2 \geq \alpha_1 \geq \alpha_2 \geq 0 \to {|\alpha|=k} }
\to
\bigoplus
\Sigma^\alpha \boxtimes \Psi_{\alpha^*}
\ .
$$

\proclaim{Remark 2.11}
\roster
\item
Note that
$$
\Phi_2 \cong \Phi \ ,
$$
since the  representations
$$
\gather
M_2 = \Lambda_2 \oplus \theta = 
\overset 2 \to \Lambda  ( V_2 \boxtimes (V/V_2)^* ) \oplus 
\overset 2 \to \Lambda V_2 \boxtimes \eta \\
\cong Sym_2 V_2 \boxtimes \overset 2 \to \Lambda (V/V_2)^*
\ \oplus \ \overset 2 \to \Lambda V_2 \boxtimes Sym_2 (V/V_2)^* \oplus
\overset 2 \to \Lambda V_2 \boxtimes \eta \\
\cong Sym_2 V_2 \boxtimes \overset 2 \to \Lambda (V/V_2)^*
\ \oplus \  \overset 2 \to \Lambda V_2 \boxtimes W 
\endgather
$$ of $P \times P$
  are isomorphic as subrepresentations
of $V^*
\otimes V^*$, thus $\Phi_2 \cong \Phi$.

\item
 The inclusion of $\Lambda_k \subset M_k$ defines a short
exact sequence of
representations of $P \times P$:
$$ 0 \to \Lambda_k \to M_k \to \overset 2 \to \Lambda V_2 \otimes M_{k-2}  
\to 0 \ ,
$$ which in turn defines the short exact sequence of homogeneous vector
bundles.
$$ 0 \to C_k \to \Phi_k \to \Cal O_X(-1) \otimes \Phi_{k-2} \ .
$$

This follows from

(a)
$
 M_k / \Lambda_k \cong 
\Lambda_{k-2} \otimes \theta \oplus \Lambda_{k-4} \otimes \theta^2 \oplus
\dots 
\
$

and

(b)
 the action of $P \times P$ on $M_k / \Lambda_k$ is compatible to the action 
of
$P \times P$ on
$\overset 2 \to \Lambda V_2 \otimes M_{k-2} $.

\endroster
\endproclaim

Similarly to the definition of $\Phi_k$, we will define homogeneous bundles
$\Phi_\alpha$ and then discuss the relation to the homogeneous bundles
$\Psi_{\alpha^*}$.

\head The extensions $\Phi_\alpha$ and $\Psi_{\alpha^*}$
\endhead

\subhead {Notations}
\endsubhead
Set $ \ M_\alpha =\underset i\geq 0 \to \bigoplus 
\Lambda_{\alpha-2i} \otimes
\theta^i $
{for} $ N-2 \geq \alpha_1 \geq \alpha_2 \geq 0 $ and recall that 
$$
\gather
\beta -2 = \left \{ 
\matrix & ( \beta_2, \beta_3, .. ,\beta_{N-2}, 0 ) && \text {if}  \ \beta_1=2
\\ & 0 && \text {otherwise}
\endmatrix \right. \ .
\endgather
$$

{ Let} $
 \Sigma_\beta =  \Sigma^{\beta} (V/V_2)^* $
  and set
$Sym_2=Sym_2 (V/V_2)^*$.

\head Definition of $\Phi_\alpha$
\endhead

\proclaim{Lemma 2.12}

$M_\alpha $ is invariant under the action of $P \times P$.
\endproclaim

\demo{Proof}

 It is sufficient to do these computations on the algebra level. Let 

$\lambda \in \Lambda_\gamma, v_1 \wedge v_2 
\in \overset 2 \to
\Lambda V_2, g,h \in \frak p $ and $i \geq 1 $, then
$$
\gather
 w=\lambda \otimes ((v_1 \wedge v_2 )^{\otimes i} \boxtimes \eta^i)
\in 
\Lambda_\gamma \otimes \theta^i
\ \ 
\text{ and}
\\ (g,h) \cdot w \\ =
  ((g,h) \cdot \lambda) \otimes ((v_1 \wedge v_2 )^{\otimes i} \boxtimes
\eta^i ) +
\lambda \cdot ((g \cdot (v_1 \wedge v_2 )^{\otimes i}) \boxtimes \eta^i) \\ +
i
\lambda \cdot  ((v_1 \wedge v_2 ) \boxtimes h \cdot \eta ) \otimes 
  ((v_1 \wedge v_2 )^{\otimes i-1} \boxtimes  \eta^{i-1})
\endgather
$$   Note that 
$$ i \lambda \cdot  ((v_1 \wedge v_2 ) \boxtimes h \cdot \eta ) \in
\Lambda_\gamma
\cdot \Lambda_{(1,1)} 
$$  and 
$\Lambda_{(1,1)}=(\overset 2 \to \Lambda V_2 \boxtimes Sym_2 )
$. Before we finish the proof, we will discuss the multiplication of 
$\Lambda_\gamma$ by
$ \Lambda_{(1,1)} $.
\enddemo

Let 
$|\alpha|=k,N-2 \geq \alpha_1 \geq  \alpha_2 \geq 0, \gamma^*=\alpha^*-2,$
that is
$\gamma=(\alpha_1-1, \alpha_2-1)$. Assume that $\gamma \geq (0,0)$.

Consider the multiplication $m_\gamma$:
$$
\alignat 8  & \Lambda_\gamma & \otimes &\quad \ \Lambda_{(1,1)} \ \to
&\Lambda_k
\qquad  &\qquad \qquad \qquad 
\qquad \qquad
\\  & \cap &&\quad \ \cap\  &\Vert \qquad \ 
\\
 & \Lambda_\gamma & \otimes &\quad \
\Lambda_2 \quad \ \ \to &\Lambda_k  \qquad
\\
 & (\lambda) \ &\otimes &\  (\lambda_1
\wedge \lambda_2)\ \mapsto &\lambda    \wedge \lambda_1 \wedge
\lambda_2
\endalignat
$$

\proclaim{Lemma 2.13}
Let $|\alpha|=k,N-2 \geq \alpha_1 \geq  \alpha_2 \geq 1, \gamma^*=\alpha^*-2 $, then
the image of $m_\gamma:\Lambda_\gamma  \otimes  \Lambda_{(1,1)} \ \to
\Lambda_k$
\roster 
\item is contained in $\Lambda_\alpha $ and
\item is not zero.
\endroster
\endproclaim

\demo{Proof}
\roster

\item Consider the tensor product 
$$
\gather
\Lambda_\gamma  \otimes  \Lambda_{(1,1)} \\
\cong (\overset 2 \to \Lambda V_2 \boxtimes Sym_2 )
\otimes (\Sigma^\gamma V_2 \boxtimes \Sigma_{\gamma^*} ) \\
\cong (\Sigma^{(\gamma_1+1, \gamma_2+1)} V_2)  \boxtimes 
(\Sigma_{(\gamma_1+1,
\gamma_2+1)^*} 
\oplus 
\Sigma_{(\gamma_1+1, \gamma_2, 1)^*}  
\oplus
\Sigma_{(\gamma_1, \gamma_2+1,1)^*}
\oplus
\Sigma_{(\gamma_1, \gamma_2,1,1)^*} ) \\
 =  (\Sigma^{(\alpha_1,
\alpha_2)} V_2)  \boxtimes  (\Sigma_{(\alpha_1, \alpha_2)^*}  
\oplus 
\Sigma_{(\alpha_1, \alpha_2-1,1)^*}  \oplus
\Sigma_{(\alpha_1-1, \alpha_2,1)^*}
\oplus
\Sigma_{(\alpha_1-1, \alpha_2-1,1,1)^*} )
\endgather
$$
See also appendix, Lemma 2(b), for this tensor product.

The image is contained in $\Lambda_\gamma  \otimes  \Lambda_{(1,1)} $ modulo
the relations of  the k-th exterior power $\Lambda^k$. Since 
$\Lambda_k \cong  
\underset |\beta|=i \to \bigoplus 
\Sigma^\beta V_2 \boxtimes \Sigma_{\beta^*} $, it follows that the image is 
 inside of
$\Sigma^\alpha V_2 \boxtimes \Sigma_{\alpha^*} $.

\item Let $v_1 \dots v_N$ be the chosen and fixed symplectic basis with its
dual basis
$x_1 \dots x_N$. Recall that $V_j=Span\{ v_1 \dots v_j \} $. 

Let $f=x_3$, then $f \in (V/V_2)^* $.

Consider
$$
\gather
\qquad ( v \wedge w ) \boxtimes f^2 \in \overset 2 \to \Lambda V_2 \boxtimes
Sym_2 \\
\Vert \qquad \qquad \qquad \qquad \cap \\ 2 ( v \boxtimes f ) \wedge ( w
\boxtimes f ) \in \overset 2 \to \Lambda ( V_2 \boxtimes  (V/V_2)^*)
\endgather
$$

Recall that $\gamma=(\alpha_1-1, \alpha_2-1)$, therefore $\gamma^*$ is of the
form

$(\gamma_1^*, \dots , \gamma_{N-3}^*, 0) $. Set
$$
\overset \wedge \to \gamma = (\gamma_1^*, \dots , \gamma_{N-3}^*) \ .
$$ Then  $V_2 \subset V_3$ and $\gamma \geq (0,0)$ imply
$$ 0 \neq
\Sigma^\gamma S \boxtimes \Sigma^{\overset \wedge \to \gamma^*}(V/V_3)^*
\subset \Sigma^\gamma S \boxtimes \Sigma^{\gamma^*}(V/V_2)^*
$$

Finally consider the multiplication of 
$\Sigma^\gamma S \boxtimes \Sigma^{\overset \wedge \to \gamma^*}(V/V_3)^* $
by 
$(v \wedge w) \boxtimes f^2 $:
$$
\gather (v \wedge w) \boxtimes f^2 \cdot w_{\overset \wedge \to \gamma} \\ =
2 ( v \boxtimes f )
\wedge ( w \boxtimes f ) \wedge w_{\overset \wedge \to \gamma}
\ .
\endgather
$$  The image is not zero, because $w_{\overset \wedge \to \gamma} $ is a sum
over elements of the form
$(v_1  \boxtimes f_1) \wedge \dots \wedge  (v_{i-2}  \boxtimes f_{i-2}) \neq 0
$ and all of these elements $v_j \boxtimes f_j $ and  $ v \boxtimes f$ and $ 
w
\boxtimes f $ are linearly independent.

Thus multiplication of  
$\Sigma^\gamma S \boxtimes \Sigma^{\overset \wedge \to \gamma^*}(V/V_3)^* $
by 
$(v \wedge w) \boxtimes f^2 $ is non-zero and therefore the multiplication  
$\Lambda_{(1,1)}  \otimes  \Lambda_\gamma  \to \Lambda_\alpha $ is non-zero.

q.e.d.

\endroster

\enddemo

\demo{Proof of Lemma 2.12}

Since $\Lambda_\gamma
\cdot \Lambda_{(1,1)} \subset \Lambda_{(\gamma_1 +1, \gamma_2 + 1)}$,  it
follows that for
 $i\geq1, w \in \Lambda_\gamma \otimes \theta^i$,$\gamma=(\alpha_1-i,\alpha_2
- i)$, 
$$ (g,h) \cdot w \in (\Lambda_\gamma \otimes \theta^i)
\ \oplus \  (\Lambda_{{(\gamma_1 +1, \gamma_2 + 1)}}
\otimes
\theta^{i-1})
$$ and for $w\in \Lambda_\alpha, (g,h) \cdot w \in \Lambda_\alpha $. This
proves that $M_\alpha $ is
 invariant under the action of $P \times P $.

q.e.d
\enddemo

\proclaim{Definition 2.14} 
$$
\text {Define} \ \ 
\Phi_\alpha=(G \times G)_{P \times P} M_\alpha \ .
$$
\endproclaim

\subhead{Relation between $M_\alpha$ and $W_{\beta^*}$}
\endsubhead

\proclaim{Proposition 2.15} Let 
$N-2 \geq \alpha_1 \geq  \alpha_2 \geq 0, $ 
$\gamma=(\alpha_1-1, \alpha_2-1) \geq (0,0)$ and let
$\alpha^*-2=\gamma^*.$ The multiplication $m_\gamma:\Lambda_\gamma   \otimes 
\Lambda_{(1,1)}\ \to \
\Lambda_\alpha$  and the $*$-multiplication are compatible, that is:

Denote by $r$ the isomorphism 
$$ r: 
\underset  {\underset
 N-2 \geq \alpha_1 \geq \alpha_2 \geq 0 \to {|\alpha|=k} }
\to
\bigoplus
\Sigma^\alpha V_2
\boxtimes
\Sigma^{\alpha^*} (V/V_2)^* 
\cong
\Lambda ,
$$   then there is an isomorphism 
$$
\gather
i= \oplus i_\alpha : \bigoplus \Lambda_\alpha \to \Lambda_\alpha
\ \text { such that} \\
 m_\gamma (((v \wedge w) \boxtimes f \cdot g ) \otimes ( v_\gamma \boxtimes
f_\gamma)) = i ( r (( v \wedge w \otimes v_\gamma ) \boxtimes (  f_\gamma * ( f
\cdot g ) )) ) \ .
\endgather
$$
\endproclaim

\demo{Proof} The multiplication $m_\gamma$ is given by the map
$$ m_\gamma:\Lambda_\gamma  \otimes  \Lambda_{(1,1)} \ \to \
\Lambda_\gamma  \otimes  \Lambda_{(1,1)} \ \text {modulo the relations of the
i-th exterior power} \
\Lambda_i\ .
$$

Consider once more the tensor product ( see also Lemma 2.13 )
$$
\gather
\Lambda_\gamma  \otimes  \Lambda_{(1,1)} \\
\overset r^{-1} \to \cong (\Sigma^{(\alpha_1,
\alpha_2)} V_2  \boxtimes   \Sigma_{(\alpha_1, \alpha_2)^*} ) \\
\oplus 
\\
\Sigma^{(\alpha_1,
\alpha_2)} V_2  \boxtimes  (\Sigma_{(\alpha_1, \alpha_2 - 1,1)^*} \oplus
\Sigma_{(\alpha_1-1, \alpha_2 , 1)^*} 
\oplus
\Sigma_{(\alpha_1-1, \alpha_2-1, 1 , 1)^*} )
\endgather
$$ Since 
$
\Sigma^{(\alpha_1,
\alpha_2)} V_2  \boxtimes  (\Sigma_{(\alpha_1+1, \alpha_2-1)^*}  \oplus
\Sigma_{(\alpha_1-1, \alpha_2+1)^*} )
$ modulo the relations of the i-th exterior product is zero, it follows that
$m_\gamma \circ r $ is given by the projection
$$
\gather 
(\Sigma^{(\alpha_1,
\alpha_2)} V_2  \boxtimes   \Sigma_{(\alpha_1, \alpha_2)^*} ) \\
\oplus 
\\
\Sigma^{(\alpha_1,
\alpha_2)} V_2  \boxtimes  (\Sigma_{(\alpha_1, \alpha_2 - 1,1)^*} \oplus
\Sigma_{(\alpha_1-1, \alpha_2 , 1)^*} 
\oplus
\Sigma_{(\alpha_1-1, \alpha_2-1, 1 , 1)^*} )  
\\ \to (\Sigma^{(\alpha_1,
\alpha_2)} V_2  \boxtimes   \Sigma_{(\alpha_1, \alpha_2)^*} ) \ .
\endgather
$$
This projection is determined up to a constant $k_\alpha$, thus
$$
 m_\gamma (((v \wedge w) \boxtimes f \cdot g ) \otimes ( v_\gamma \boxtimes
f_\gamma)) = k_\alpha
r (( v \wedge w \otimes v_\gamma ) \boxtimes ( f_\gamma * ( f \cdot
g )))\ .
$$
Set $i_\alpha = k_\alpha id|_{\Lambda_{\alpha}} $.

 q.e.d.

\enddemo

\proclaim{Remark}
We will now identify 
$$
\Lambda = \Lambda ( V_2 \boxtimes (V/V_2)^* )
\ \ \text{ with} \ \
\underset  {\underset
 N-2 \geq \alpha_1 \geq \alpha_2 \geq 0 \to {|\alpha|=k} }
\to
\bigoplus
\Sigma^\alpha V_2
\boxtimes
\Sigma^{\alpha^*} (V/V_2)^* \ .
$$ 

\endproclaim

\proclaim{Theorem 2.16}

For $ N-2 \geq \alpha_1 \geq \alpha_2 \geq 0 $,
$$
\gather
M_\alpha \cong \Sigma^\alpha V_2 \boxtimes W_{\alpha^*} \\
\Phi_\alpha \cong \Sigma^\alpha S \boxtimes \Psi_{\alpha^*}
\endgather
$$

\endproclaim

\demo{Proof}

Let $\alpha^*=\beta-2i, v_\alpha \boxtimes f_\alpha \in 
\Sigma^\alpha V_2 \boxtimes \Sigma_{\alpha^*} $ and $(g,h) \in 
\frak p \times \frak p $. Consider
$$
\gather  (g,h) \cdot ((v_\alpha \otimes  (v_1 \wedge v_2 )^{\otimes i})
\boxtimes (f_\alpha \otimes  \eta^i ) )
\\   = (g \cdot v_\alpha \otimes  (v_1 \wedge v_2 )^{\otimes i} + v_\alpha
\otimes g \cdot (v_1 \wedge v_2 )^{\otimes i})
\boxtimes (f_\alpha \otimes  \eta^i ) )
\\ +((v_\alpha \otimes  (v_1 \wedge v_2 )^{\otimes i})
\boxtimes ( h \cdot f_\alpha \otimes \eta^i )
\\ + (   (v_1 \wedge v_2 )^{\otimes i-1}) \otimes (v_\alpha
\boxtimes f_\alpha)  \cdot ( v_1 \wedge v_2  \boxtimes h \cdot \eta)
 \otimes
\eta^{i-1} \\  =(g \cdot v_\alpha \otimes  (v_1 \wedge v_2 )^{\otimes i} +
v_\alpha \otimes g
\cdot (v_1 \wedge v_2 )^{\otimes i})
\boxtimes (f_\alpha \otimes  \eta^i ) )
\\  +((v_\alpha \otimes  (v_1 \wedge v_2 )^{\otimes i})
\boxtimes ( h \cdot f_\alpha \otimes \eta^i + ( f_\alpha * h \cdot \eta )
\cdot 
\eta^{i-1} ) 
\\ =(g \cdot (v_\alpha \otimes  (v_1 \wedge v_2 )^{\otimes i})
\boxtimes (f_\alpha \otimes  \eta^i ) )
\\ + ((v_\alpha \otimes  (v_1 \wedge v_2 )^{\otimes i})
\boxtimes h \cdot (f_\alpha \otimes  \eta^i ) )
\endgather
$$   The second equality follows from Proposition 2.15. These computations
 imply that $W_\beta$ is a representation of $\frak p$  coming from the
 representation $M_{\beta^*} $ of $\frak p
\times \frak p$ and thus
$$
 M_{\beta^*} \cong \Sigma^{\beta^*} V_2 \boxtimes W_\beta \ .  
$$ 

The isomorphism of the bundles
 follows immediately from here, since
$$
\gather
\Phi_\alpha \cong ( G \times G)_{P \times P} M_\alpha \cong ( G \times_P
\Sigma^\alpha V_2 ) \boxtimes ( G \times_P W_{\alpha^*})
\cong \Sigma^\alpha S \boxtimes \Psi_{\alpha^*} \ .
\endgather
$$ q.e.d.
\enddemo

We complete the Tate construction by extending the 
connecting morphism $d$ from the Koszul
complex
$C_. \otimes \Cal O_{X \times X}$ to the new construction $D_.$,
 thus extending the Koszul complex
$C_.$ to the  complex $D_.$.

\head {The map d}
\endhead

Recall that 
$ \eta = x_1 \otimes x_{n+1} + x_2 \otimes x_{n+2} + \dots + x_n \otimes
x_{2n} 
$ for the standard symplectic basis  $v_1, \dots , v_{2n}$ and its dual basis 
$x_1, \dots , x_{2n}$.

\proclaim{Definition 2.17}

\roster

\item  For $g,h \in G , i > 0 $, set
$$
\gather
  T_{(g,h)} = (g, v_1 \wedge v_2 ) \boxtimes (h,\eta) \in \Phi_{2, (gP, hP)}\
\ 
\text
 {and} \\
 T^i_{(g,h)} = (g, (v_1 \wedge v_2)^{\otimes i } ) 
\boxtimes (h,\eta^i)  \in \Phi_{2i, (gP, hP)}\ .
\endgather
$$

\item  Let $\overset \wedge \to d: D. \to D. $ be the extension of the map $d$
of the  Koszul complex 
$C. \otimes \Cal O_{X \times X}$ satisfying
 over a point 
$(gP, hP ) \in G/P \times G/P \cong X \times X  $:
$$
\align
\overset \wedge \to d(T_{(g,h)}) &=  (g, h \cdot x_1 ( g \cdot v_1 ) v_2 -  h
\cdot x_1 ( g \cdot v_2 ) v_1)
\boxtimes (h, x_{n+1} )
\tag a
\\ &+ (g, h \cdot x_2 ( g \cdot v_1 ) v_2 -  h \cdot x_2 ( g \cdot v_2 ) v_1)
\boxtimes (h, x_{n+2} )
\\  &+ \dots + \\ &+ (g, h \cdot x_n ( g \cdot v_1 ) v_2 -  h \cdot x_n ( g
\cdot v_2 ) v_1)
\boxtimes (h, x_{2n} ) 
\\ &=: t_{(g,h)} \in S_{gP} \boxtimes \st_{ hP} \ ,
\\
\overset \wedge \to  d ( T^i_{(g,h)} ) &= i \ t_{(g,h)} \ T^{i-1}_{(g,h)} \ \
\text {and}
\tag b
\\
\overset \wedge \to  d ( w T^i_{(g,h)} ) &=  ( dw) T^i_{(g,h)} +  (-1)^{ | w|
} i  w
\cdot  t_{(g,h)}  T^{i-1}_{(g,h)}) 
\
\text{for} \ w \in (C_j \otimes \Cal O_{X \times X} )_{(gP,hP)} \ .
\tag c
\endalign
$$

\endroster

\endproclaim

\proclaim{Lemma 2.18}
$\overset \wedge \to  d$ is  well defined and uniquely determined and 
$\overset \wedge \to  d$ restricted to the
 Koszul complex $C. \otimes \Cal O_{X \times X} $ is equal to $d$.

\endproclaim

\demo{Proof}

\item{(a)}  Consider 
$\overset \wedge \to d(T_{(g,h)}) $. Note that 
$T_{(g,h)} $ is an element in
$(( G \times G)_{P \times P} M_2)|_{(gP, hP)}$ and note that the
$P \times P$-representation
$M_2 = \overset 2 \to \Lambda V_2 \boxtimes W$ is contained in 
$\overset 2 \to \Lambda V_2 \boxtimes V^* \otimes (V/V_2)^*$.  Define 
$$
\gather
\overset \sim \to d: ( G \times_P \overset 2 \to \Lambda V_2 ) \ 
\boxtimes \ ( G \times _P ( V^* \otimes (V/V_2)^* )) \to ( G \times_P  V_2 )
\ 
\boxtimes \ ( G \times _P  (V/V_2)^* )
\\ (g, v \wedge w) \boxtimes (h, F \otimes G ) \mapsto (g, h \cdot F (g \cdot
v ) - h
\cdot F (g \cdot w) ) \boxtimes (h,G) \ .
\endgather 
$$  Clearly  
$\overset \sim \to d =\overset \wedge \to d $ on 
$\Phi = ( G \times G)_{P \times P} M_2 $. 
$\overset \sim \to d$ is the evaluation map on the first factor and it
 is well defined, since:
$$
\gather
 \overset \sim \to d( (g p_1,v \wedge w) \boxtimes (h p_2 ,F \otimes G) ) \\ =
(g p_1,h p_2 F (g p_1 v) w - h p_2 F (g p_1 w) v) \times (h p_2 , G)
\\
\text {and} \\
\overset \sim \to  d( (g, p_1 ( v \wedge w )) \boxtimes (h ,p_2 ( F \otimes
G)) ) \\  = (g,p_1 ( h p_2 F (g p_1 v) w - h p_2 F (g p_1 w) v) ) \times (h,
p_2  G)
\endgather
$$ are the same for all $p_1, p_2 \in P$.

Thus $\overset \wedge \to d=\overset \sim \to  d|_\Phi $ is well defined on
$\Phi$.

\item{(b)} The condition in (2) (b)   assures that
$\overset \wedge \to  d $ respects the relation $T^i T^j = T^{i+j} $, that is
$$
\overset \wedge \to  d ( T^{i}_{(g,h)} T^{j}_{(g,h)}) =\overset \wedge \to d (
T^{i}_{(g,h)}) T^{j}_{(g,h)} +   T^{i}_{(g,h)} \overset \wedge \to d
(T^{j}_{(g,h)})  = (i+j) \ t_{(g,h)} \ T^{i+j-1}_{(g,h)}
$$

\item{(c)} 
$\Phi_k = ( G \times G)_{P \times P} M_k $ and elements in $M_k$ are sums over
elements of the form
$ \lambda_{k-2i} \otimes ((v_1 \wedge v_2)^{\otimes i} \boxtimes \eta^i) $.

Thus over the point $(gP, hP) \in X \times X $, the elements of $\Phi_k|_{(gP,
hP)}$ are sums over elements of the form
$$
 w_{k-2i} \otimes ((g,h),(v_1 \wedge v_2)^{\otimes i} \boxtimes \eta^i) = 
w_{k-2i} \otimes T^{i}_{(g,h)} 
$$

 This determines $\overset
\wedge \to d $ uniquely and it is well defined since $\overset \wedge \to
d|_\Phi
$ is well defined .

 q.e.d.
\enddemo

Set $d = \overset \wedge \to d$.

\proclaim{Remark 2.18}
 Note that
\roster
\item
$ d^2=0$ and
\item
$d$ maps $\Phi_\alpha$ to $\Phi_{\alpha-1}$, where
$\Phi_{\alpha-1} $ denotes the sum of $\Phi_\gamma$ over all Young diagrams
$\gamma$ that are contained in
$\alpha$, i.e. $\gamma_i \leq \alpha_i $, and of length $|\alpha|-1$.
\item There is a short exact sequence of homogeneous vector bundles
$$ 0 \to \Sigma^\beta \st \to \Psi_\beta \to \Psi_{\beta-2} \to 0 \ ,
$$ for $2 \geq \beta_1 \geq \beta_2 \geq \dots \geq \beta_{N-2} \geq 0 $.
\item The Tate complex $D_.$  in degree $k$ is given by:
$$ D_k \cong \underset 
{\underset
 N-2 \geq \alpha_1 \geq \alpha_2 \geq 0 \to {|\alpha|=k} }
\to
\bigoplus
\Sigma^\alpha S
\boxtimes
\Psi_{\alpha^*} \ .
$$
\endroster

\endproclaim

\demo{Proof}
\roster
\item
 This follows from  $t \cdot t= t \wedge t =0$ and
$$
\gather dd(wT^i)=d((dw)T^i + (-1)^{|w|}i w \cdot t T^{i-1}) \\ =ddw T^i +
(-1)^{|dw|} i dw \cdot t T^{i-1} \\ + (-1)^{|dw|-1} i dw \cdot t T^{i-1} + 
(-1)^{|w|} (-1)^{|w|} i w d( t T^{i-1}) \\ = 0 + 0 + 0 + i w d( t T^{i-1}) \\
= i w dt T^{i-1} - (i-1) t \cdot t \ T^{i-2} = 0 \ .
\endgather
$$
\item Consider the map $d: M_\alpha \to M_{k-1}$ of $\frak p$ representations.
Let $|\alpha|=k w \in
\Lambda_\gamma, \gamma=(\alpha_1-i, \alpha_2 -i)$, then
$$
\gather d( w \eta^i ) = (dw) \eta^i + (-1)^{|w|} (w \cdot d \eta) T^{i-1} \\
\in \Lambda_\gamma \eta^i \oplus \Lambda_{(\gamma_1 +1, \gamma_2 +1)}
\eta^{i-1},   
\\
\text{because} \  w \cdot t \in \Lambda_\gamma \cdot \Lambda_{(1,0)} 
\subset \Lambda_{(\gamma_1 +1, \gamma_2) } \oplus \Lambda_{(\gamma_1 ,
\gamma_2+1)}\ .
\endgather
$$
 Note that
$d$ is not a homogeneous map.

\item The inclusion $\Sigma^\beta (V/V_2)^* \subset W_\beta $ defines a short
exact sequence of $P$ representations:
$$
 0 \to \Sigma^\beta (V/V_2)^* \to W_\beta \to W_{\beta-2} \to 0 \ ,
$$ which implies the  exactness of the corresponding sequence of homogeneous
vector bundles.
\item  This is obvious, since 
$$ D_. = \underset k \geq 0 \to
\bigoplus
\Phi_k \cong
\underset {\underset
 N-2 \geq \alpha_1 \geq \alpha_2 \geq 0 \to {|\alpha|=k} }
 \to
\bigoplus
\Phi_{\alpha^*} \ . \ 
$$

\endroster q.e.d.

\enddemo

We will show next that the complex $D_.$, that we have defined, is exact and thus
does resolve  $\Cal O_\Delta$.

\specialhead{The resolution}
\endspecialhead

\proclaim{Theorem 2.19}
$$ D.: \ \  \dots \to \Phi_k @> d>> \Phi_{k-1} @> d>> \dots @> d>>
\Phi @> d>>S \boxtimes \st @> d>> \Cal O_{X \times X} 
$$ is a left resolution of the structure sheaf $\Cal O_\Delta$ of the diagonal
$\Delta$ inside of $X \times X$ .

\endproclaim

\demo{Proof}  This is a local question.  Let $U_1, U_2$ be  open sets in $X$.
We may assume that all of the  vector bundles ${D_i}|_{U_1 \times U_2}$ are
trivial.   Thus
$$
\gather
\Psi|_ {U_2} \cong Sym_2 S|_ {U_2} \oplus \Cal T \ \ \text {and}\\
\Phi_k|_ {U_1 \times U_2} \cong 
\underset i\geq 0 \to \bigoplus
\overset k-2i \to \Lambda ( S|_{U_1} \boxtimes \st|_{U_2} ) 
 \ \Cal T^i
\endgather
$$ Therefore ${D_i}|_{U_1 \times U_2}$ is isomorphic to the  polynomial ring
in one variable $\Cal T$ of degree 2 over 
$(C. \otimes \Cal O_{X \times X})|_{U_1 \times U_2} $. Thus
$$ {D_i}|_{U_1 \times U_2} \cong (C. \otimes \Cal O_{X \times X})|_{U_1 \times
U_2}[\Cal T]  \ .
$$ We will show that  this is the same construction as in Tate [T].  Since
$D.|_{U_1
\times U_2}$ is gotten by adjoining the variable $\Cal T$ in degree 2 to $(C.
\otimes
\Cal O_{X \times X})|_{U_1 \times U_2} $, this complex will kill the class of
the cycle
$d\Cal T$, see [T],Theorem 2 or Lemma 2.20 below.

We will prove this theorem in two parts:

First we discuss the Tate construction of adjoining a variable in even degree
and then we will show that the homology of the complex
$C. $ restricted to $X \times X$ is given by the cycle $d\Cal T $. 

\enddemo

\head{Adjoining a variable in even degree}
\endhead

Let $\varLambda$ be an $R$-algebra, that is a graded algebra over a
commutative noetherian ring  $R$ with unit  element and an $R$-linear mapping 
$d: \varLambda_i \to
\varLambda_{i-1}$ satisfying 
\item{(1)}
$\varLambda_i=0 $ for $i<0$, $\varLambda_0 = R\cdot 1 $ and $\varLambda_i$ is
a  finitely generated $R-$module for all $i>0$.
\item{(2)}
$\varLambda$ is skew-commutative and
\item{(3)}
$d$ is a skew derivation, that is
$$ d ( x y ) = (dx) y + (-1)^ {deg(x)} x (dy) \ .
$$ 

Let $t$ be a cycle of odd degree $\rho-1$. We will adjoin a variable in even
degree
$\rho$ in order to kill the cycle $t$.

Let $Y=\varLambda[T]$ be the ring of polynomials in one commuting variable
$T$ of degree
$\rho$ over $\varLambda$. Extend the derivation $d$ in $\varLambda$ to $Y$ in
the unique way such that $dT=t$. Then $Y$ is an $R$-algebra.

\proclaim{Lemma 2.20}  Assume that  the homology class $[t] \in
H(\varLambda)$ is a skew non-zero divisor, that is
$$ [t] \xi = 0 \Rightarrow \xi \in [t] H(\varLambda)  \ 
\text{ for all} \ \xi \in H(\varLambda)
\ .
$$ Then $H(Y) \cong H(\varLambda) / \ [t]  H(\varLambda) \ $.
\endproclaim

\demo{Proof} See [T], Theorem 2, page 17.

\enddemo

We continue with the proof of the theorem.

Set $\varLambda=(C. \otimes \Cal O_{X \times X})|_{U_1 \times U_2} $, 
$Y=D.|_{U_1 \times U_2}$ and $t = d \Cal T $. Then $Y=\varLambda[\Cal T]$.
Suppose that the homology of $\varLambda$ is generated by $[t]$. Then the
condition of the previous Lemma is satisfied, which implies that $Y$ has no
homology, thus that $D.$ is exact.

Thus left to show that the homology of
 the complex 
$(C. \otimes \Cal O_{X \times X})|_{U_1 \times U_2}$ is generated by $[t]$. 

\head{The torsion}
\endhead

\proclaim{Lemma 2.21}
$$
\tor_i ^{\Cal O_{Gr \times Gr}} ( \Cal O_{\Delta_{Gr}} , \Cal O_{X \times X})
\cong
\cases
  \Cal O_\Delta     & i=0 \\
  \Cal O_\Delta (-1)   &   i=1 \\
  0    &  i>1
\endcases
$$
\endproclaim
\demo{Proof}

Recall that $X$ in $ Gr(2,V) $ is given by the hyperplane section $ z=0$ with
$$ z=x_1 \wedge x_{n+1} + x_2 \wedge x_{n+2} + ... + x_n \wedge x_{2n}
$$ where $ x_1 ... x_{2n} $  is the symplectic basis of $V^*$ . Let $Gr(2,V) =
Proj(A)$, then $X = Proj(A/z)$. Consider the short exact sequence:
$$ 0 \to \Cal O_{Gr } (-1)  @>  z >> \Cal O_{Gr } \to \Cal O_X \to 0 \ .
$$ 

Define $z_1= z \boxtimes 1 , z_2 =  1 \boxtimes z $. Now $z_1$ and $z_2$ form
a regular sequence in $A \boxtimes A $, hence the  Koszul complex of $z_1$ and
$z_2$ is exact. Since $X \times X = Proj ( A \boxtimes A / ( z_1, z_2 ) ) $
this gives a resolution of 
$\Cal O_{X \times X} $:
$$
\gather 0 \to \Cal O_{Gr\times Gr}(-1,-1) \to \Cal O_{Gr\times Gr} (-1,0) 
\oplus \Cal O_{Gr\times Gr}(0,-1) \\ \to
\Cal O_{Gr \times Gr} \to \Cal O_{X \times X} \to 0 \ .
\endgather
$$  To compute the $Tor$'s we tensor this complex by $ \Cal O_{\Delta_{Gr}}$
and compute cohomology of the thus obtained complex. Since $\Delta_{Gr} \cong
Gr$ and
$\Delta \cong X $ we get
$$ 0 \to \Cal O_{Gr} (-2) @> \beta > \left[ \matrix z & -z \endmatrix \right]
>
\Cal O_{Gr} (-1)
\oplus
\Cal O_{Gr}(-1)  @> \alpha > 
\left[ \matrix     z \\  z \endmatrix \right]  > \Cal O_{Gr} \to
\Cal O_{X } \to 0
$$ This complex has only cohomology in degree 0 and 1. In degree 0, $\oh^0
\cong
\Cal O_X $ and in degree 1, we have  
 
 $\alpha ( \underset i \to \sum (f_i \oplus g_i) ) = z \cdot (\underset i \to
\sum (f_i + g_i) ) = 0 
\Leftrightarrow \underset i \to \sum( f_i + g_i ) =0 $, thus 
$Ker \alpha \cong (A[-1])^\sim \cong \Cal O_{Gr} (-1)$,  and since $\beta(f)=z
\cdot f
\oplus ( - z \cdot f) $, we get
$Im \beta \cong z \cdot (A [-2])^\sim$.

 Therefore 
$ \oh^1 \cong (A[-1]/z A[-2])^\sim \cong
\Cal O_X(-1) \cong \Cal O_{\Delta}(-1) \ $. \qquad q.e.d. 
\enddemo

\proclaim {Corollary 2.22}

Let $y_1, ..,y_N$ be a symplectic basis of $V^*$, such that over a point $W_1
\times W_2 \in X \times X ,( V / W_2 )^*$ is generated by $ y_3, .. , y_N $.
Let
$ W_1 =v \wedge w$. Then the torsion over the point $W_1
\times W_2$ is generated by  the cycle $\gamma_{W_1 \times W_2}$ :
$$
\gamma_{W_1 \times W_2} = ( y_1(v) w - y_1 (w) v) \boxtimes y_{n+1} + ... + (
y_n(v) w - y_n (w) v) \boxtimes y_{2n} .
$$
\endproclaim
\demo {Proof}
$W_1$ is  an isotropic subspace of $V$ . Thus for $ v,w \in W_1$:
$$
\gather d \gamma_{W_1 \times W_2} =  y_1(v) y_{n+1} ( w) - y_1 (w) y_{n+1}(
v)  + ... +
 y_n(v)  y_{2n} ( w ) - y_n (w)  y_{2n} (v) \\
 = < v,w >=0
\ .
\endgather  
$$
 q.e.d.
\enddemo

\head{Proof of Theorem 2.19}
\endhead In order to finish the proof we have to show that 
$d \Cal T$ generates the homology of the complex
$C. \otimes \Cal O_{X \times X} $.  This follows immediately from the last
corollary, since
$d \Cal T_{(gP,hP)} = d T_{(g,h)} = t_{(g,h)} = \gamma_{(gP,hP)}  =
\gamma_{W_1
\times W_2} $, for the parabolic $gP$ and $hP$ that are preserving $W_1$
respectively $W_2$.

 q.e.d.

\specialhead The periodic resolution
\endspecialhead

We started with the Koszul resolution of $\Cal O_{\Delta_{Gr} }$ and we used
the Tate construction to define the  resolution $D.$ of $\Cal O_{\Delta} $.
The first thing we note is that this complex is not finite anymore, but
instead it becomes periodic in large degrees of period two. To be more
precise:

Recall from remark 2.18 that
 $$  0 \to \overset k \to \Lambda ( S \otimes \st ) \to \Phi_k \to \Phi_{k-2}
\otimes ( \Cal O_X(-1) \boxtimes \Cal O_X ) \to 0  
$$ { is exact }. In particular if $k>$ rank $ ( S \boxtimes \st ) = 2 \cdot (
N-2 ) $, then 
$$
\Phi_k \cong \Phi_{k-2} \otimes ( \Cal O_X(-1) \boxtimes \Cal O_X ) = \Cal
O_X(-1) \otimes \Phi_{k-2}
$$

Set $$\Phi_+ = \Phi_{2N-4},\ \Phi_- = \Phi_{2N-5} $$ and denote by $\Cal F
(-k) = \Cal F \otimes ( \Cal O_X(-k) \boxtimes \Cal O_X )$ 

Then the complex $D.$ becomes:

\proclaim {Corollary 2.23}

The complex $D.$ in large degrees is periodic of degree 2, that is $D.$ is equal to:
$$
\gather ... \to \Phi_+ (-i)  \to 
 \Phi_- (-i)  \to
\Phi_+ (-i+1) 
\to 
 \Phi_- (-i+1)  \\  \to...
\to  \\ \Phi_+ \to  \Phi_- 
\to \Phi_{2N-6} \to ... \to S \boxtimes \st \to 
\Cal O_{X \times X} .
\endgather
$$
\endproclaim

This completes the Tate construction of the resolution $D_.$.

\addline

\proclaim{Remark}
Although we will not use it in this work, we will  discuss one more issue.
The extensions $\Psi_\beta$ are extensions of $\Sigma^\beta \st$ and
of $\Psi_{\beta - 2}$, see Remark 2.18. 
It turns out that these are non-trivial extensions.

This is a very essential comment for the K-theory of $X$. The goal is to show algebraically
that
$$
K_* (X) \cong \underset \Cal G \in \Cal X \to \bigoplus K_*(\Bbbk) ,
$$
for a generating system $\Cal X$. Particularly the elements of a generating
system should not have any cohomology. 
We note without proof that, for example $ Sym_2 \st$  has higher cohomology
and that is why we chose to extend $Sym_2 \st$ to $\Psi$.
If $\Psi$ were a trivial extension, it would still have higher cohomology.
\endproclaim

\specialhead {Non-triviality of the extensions $\Psi_{\beta}$ }
\endspecialhead

\proclaim {Theorem 2.24}
$\Psi$ is a non-trivial extension of $Sym_2 \st$ and $\Cal O_X$, i.e.
$$ 0 \to Sym_2 \st \to \Psi \to \Cal O_X \to 0 \ \text { does not split.}
$$

\endproclaim

\demo{Proof} Consider the short exact sequence
$$ 0 \to \Cal O_{Gr} (-1) \to \Cal O_{Gr} \to \Cal O_X \to 0 \ .
$$

Tensor this sequence by $V^* \otimes \st_{Gr}$:
$$ 0 \to V^* \otimes \st_{Gr} (-1) \to V^* \otimes \st_{Gr} \to V^* \otimes
\st
\to 0 \ .
$$ Consider the long exact cohomology sequence:
$$
\gather 0 \to H^0 (V^* \otimes \st_{Gr} (-1)) \to H^0 (V^* \otimes \st_{Gr} ) 
\to H^0 (V^* \otimes \st ) \\
\to H^1( V^* \otimes \st_{Gr} (-1)) \to H^1 (V^* \otimes \st_{Gr} ) \to  H^1
(V^* \otimes \st ) 
\to \dots
\endgather
$$

Note that
 $$
\gather
\st_{Gr} (-1) \cong \st_{Gr} \otimes \overset 2 \to \Lambda \st_{Gr} 
\cong \young{1} \otimes \yng{1,1} \cong \yng {2,1} \oplus \yng{1,1,1} 
\cong 
 \Sigma^{(2,1)} \st_{Gr} \oplus \overset 3 \to \Lambda \st_{Gr}  
\endgather
$$  ( see appendix for the tensor product of two $ Gl_M$ representations ).

Since $\Sigma^{(2,1)} \st_{Gr}$ and  
$\overset 3 \to \Lambda \st_{Gr}$ have no cohomology,  see Lemma 1.9, chapter
1, it follows that
$V^* \otimes \st (-1)$ has no cohomology. Also  $V^* \otimes
\st_{Gr} $ has no cohomology. From the long exact cohomology sequence it
follows that 
 $V^* \otimes \st $ has no cohomology, particularly $V^* \otimes \st $ has no
non-trivial sections.

Suppose the sequence
$0 \to Sym_2 \st \to \Psi \to \Cal O_X \to 0
 $ splits, then $\Psi$ has a non-trivial section. Since $\Psi$ is inside of
$V^*
\otimes \st$,
 it follows that $V^* \otimes \st$ has a non-trivial section, which is a
contradiction to the
 cohomology computations above. q.e.d.
\enddemo

\proclaim{Proposition 2.25}
 The short exact sequence
$$ 0 \to \Sigma^\beta \st \to \Psi_{\beta} \to \Psi_{\beta-2} \to 0
$$ does not split, i.e. 
$\Psi_\beta$ is a non-trivial extension of 
$\Sigma^\beta \st $ and $\Psi_{\beta-2}$.
\endproclaim

\demo {Proof} We have to show that $\Psi_\beta \in  Ext^1 ( \Psi_{\beta-2} ,
\Sigma^\beta \st ) $ is not zero.  Tensor the short exact sequence above by 
$ \Sigma^{\beta-2} ( \st )^*$, then it is enough to show that 
$$
\Psi_\beta \otimes \Sigma^{\beta-2} ( \st )^* \in  Ext^1 ( \Psi_{\beta-2}
\otimes \Sigma^{\beta-2} ( \st )^* , 
\Sigma^\beta \st \otimes \Sigma^{\beta-2} ( \st )^* ) 
$$ is not zero. 

Consider
\item{(a)}
 the identity inclusion
$$ i: \Cal O_X \hookrightarrow \Sigma^{\beta-2} \st \otimes \Sigma^{\beta-2} (
\st )^*
\subset \Psi_\beta \otimes \Sigma^{\beta-2} ( \st )^* \ \text{and}
$$
\item{(b)} the evaluation map:
$$ ev: \Sigma^\beta \st \otimes \Sigma^{\beta-2} ( \st )^* \to Sym_2 \st \ , 
$$ defined by
$$
\gather
\Sigma^\beta \st \otimes \Sigma^{\beta-2} ( \st )^* 
\subset Sym_2 \st \otimes \Sigma^{\beta-2} \st \otimes \Sigma^{\beta-2} ( \st
)^* 
\to Sym_2 \st \\ f \cdot g \otimes f_{\beta-2} \otimes v_{\beta-2} \mapsto 
f_{\beta-2} ( v_{\beta-2} ) \ f \cdot g \ .
\endgather
$$ These maps induce maps on the extensions:
$$
\alignat 4
\Psi_\beta &\otimes \Sigma^{\beta-2} ( \st )^* &\in  Ext^1 ( \Psi_{\beta-2}
&\otimes
\Sigma^{\beta-2} ( \st )^* , 
\Sigma^\beta \st \otimes \Sigma^{\beta-2} ( \st )^* ) \\ &\downarrow & &
\downarrow \\ &\Cal F &\in Ext^1 ( \Psi_{\beta-2} &\otimes \Sigma^{\beta-2} (
\st )^* , 
\Cal O_X ) \\ &\downarrow && \downarrow \\ &\Cal G &\in \quad \quad \ Ext^1
(&Sym_2 \st , 
\Cal O_X ) 
\endalignat 
$$ We will show that $\Cal G$ is a non-trivial extension, thus 
 $\Psi_{\beta} \otimes \Sigma^{\beta-2} ( \st )^* $ and 
$\Psi_\beta$ are non-trivial extensions.

Let us first  discuss $\Cal F$, set
$$ W_\Cal F= \Sigma^\beta (V/V_2)^* \otimes \Sigma^{\beta-2} ( V/V_2)
\oplus \eta \cdot 1 \subset W_\beta \otimes \Sigma^{\beta-2} ( V/V_2),
$$ where $1 \in \Sigma^{\beta-2} (V/V_2)^* \otimes \Sigma^{\beta-2} ( V/V_2) $
is the element corresponding to the identity element in
$ End (\Sigma^{\beta-2} (V/V_2)^* ) 
\cong \Sigma^{\beta-2} (V/V_2)^* \otimes \Sigma^{\beta-2} ( V/V_2)$.					

We claim that the  action of $P$ leaves $W_\Cal F$ invariant.  Again we  do these
computations on the algebra level:
$$
\gather
\text {For} \ p \in \frak p: p \cdot ( \eta \cdot 1 ) = \pi ((p \cdot  \eta )
\cdot 1 ) + \eta \cdot (p \cdot 1) = \pi ((p \cdot  \eta ) \cdot 1 ) + 0
\\
\in \pi (Sym_2 (V/V_2)^* \otimes \Sigma^{\beta-2} (V/V_2)^* \otimes
\Sigma^{\beta-2} ( V/V_2) ) \\
\subset \Sigma^{\beta} (V/V_2)^* \otimes \Sigma^{\beta-2} ( V/V_2) \\
\subset W_\Cal F \ .
\endgather
$$ 

Claim: $\Cal F \cong G \times_P W_\Cal F $. This follows from the commutative
diagram of
representations of $P$:

$$
\alignat 8 0 \to \Sigma^{\beta} (V/V_2)^* &\otimes \Sigma^{\beta-2} ( V/V_2) &
\to  & W_\beta  &\otimes &\ \Sigma^{\beta-2} ( V/V_2) & \to  & W_{\beta-2} &
\otimes & \Sigma^{\beta-2} ( V/V_2) &\to 0\\ &\ \Vert & &&\cup & & & &\uparrow
\\ 0 \to 
\Sigma^{\beta} (V/V_2)^* &\otimes \Sigma^{\beta-2} ( V/V_2) & \to & &  W_\Cal
F && \to & &\Bbbk  &\to 0
\endalignat
$$

Next consider $\Cal G$:

We will see that $\Cal G \cong \Psi = G \times_P W$,  which by Theorem 2.4 is
a non-trivial extension  of $Sym_2 \st$ and $\Cal O_X$.

Let
$W_\Cal F \to W $ be the map of representations of $P$:
$$
\gather
\Sigma^\beta (V/V_2)^* \otimes \Sigma^{\beta-2} ( V/V_2)
\oplus \eta \cdot 1 \to Sym_2(V/V_2)^* \oplus \Bbbk \eta \\ f_\beta \otimes
v_{\beta-2} + c \eta
\cdot 1 \mapsto ev ( f_\beta \otimes v_{\beta-2})  + c
\eta  \ \ .
\endgather
$$ 

This map induces the commutative diagram of  representations of $P$:

$$
\matrix 0 &\to  &\Sigma^{\beta} (V/V_2)^* \otimes  \Sigma^{\beta-2} ( V/V_2) &
\to & &  W_\Cal F & \to & \Bbbk &\to 0 \\ &&  \downarrow ev \quad \quad & & &
\downarrow  & & \Vert \\ 0 &\to & Sym_2 (V/V_2)^* \quad \ & \to & &  W & \to &
\Bbbk &\to 0 & \ , \quad \quad \quad \quad
\endmatrix
$$ which in term induces the commutative diagram of homogeneous bundles:

$$
\matrix 0 &\to  &\Sigma^{\beta} \st \otimes  \Sigma^{\beta-2} ( \st)^* & \to &
&  \Cal F & \to &
\Cal O_X &\to 0 \\ &&  \downarrow ev \quad \quad & & &
\downarrow  & & \Vert \\ 0 &\to & Sym_2 (V/V_2)^* \quad \ & \to & &  \Psi &
\to & \Cal O_X &\to 0 &
\ . 
\endmatrix
$$

Thus $\Cal G \cong \Psi$ and therefore $\Cal G$ is a non-zero  element in 
$Ext^1(\Cal O_X, Sym_2 \st ) $.

 q.e.d.
\enddemo

\newpage

\topmatter
\title\chapter{3} The Main Theorem \endtitle
\endtopmatter

Let $X$ be the symplectic Grassmannian $X=SpGr(2,V)$. 
In the previous chapter we used Tate's techniques to define the  
complex $D.$, which turned out to be  a resolution of the structure sheaf of the diagonal 
$\Delta \subset
X \times X$.

In the case of  ordinary Grassmannians, chapter 1, we were able to define a finite resolution $C.$ of
$\Cal O_{\Delta_{Gr}}$ of the diagonal
$\Delta_{Gr} \subset
Gr \times Gr$. 
In the symplectic case of the
symplectic grassmannian $X$  we 
constructed the infinite resolution $D.$
  of $\Cal
O_\Delta$. 
This resolution  does not
terminate, instead it becomes periodic, see Corollary 2.23.

However it does contain a finite sub-complex $B.$.

\proclaim {Definition 3.1}
For $ i \leq 2N-6$, let
$$
B_i= \bigoplus \Sigma ^\alpha S 
\boxtimes \Psi_{\alpha^*}\ ,
$$
where the sum goes over all $\alpha$ with 
$l ( \alpha ) =i $ and $  N-2 \gvertneqq \alpha_1 \geq \alpha_2 \geq 0 $.
\endproclaim

Note that this sum excludes all $\alpha$ with $\alpha_1=N-2$.

Let us check that $B.$ is a sub-complex.
 Recall that
$\alpha -1 =\{ \gamma \subset \alpha \text{ and } l(\gamma)=l(\alpha) - 1 \}
$.
Let $\gamma \in \alpha -1 $, then $\alpha_1 \neq N-2$ implies $\gamma_1 \neq N-2$.
Since $d$ maps  $\Sigma ^\alpha S \boxtimes \Psi_{\alpha^*}$ to  
$\underset { \gamma \in \alpha-1 } \to \bigoplus \Sigma ^{\gamma} S \boxtimes
\Psi_{\gamma^*}$, Remark 2.18 (2),  $d$ maps $B.$ to itself.
(

We will devote  the next three chapters to the proof of the Main Theorem:

\proclaim {Theorem 3.2}
The sub-complex $B.$:
$$
B_{2N-6} \to ... \to B_2 \to S \boxtimes S^\perp \to \Cal O_{X \times X}
$$
is exact. 
\endproclaim

\specialhead
First part of the proof:
\endspecialhead

Let $Q.$ be the quotient complex of the total complex and the sub-complex.
In the following section on "exactness of sub- and quotient complexes", we
prove that it is enough to show that the quotient  complex is exact up to degree $2N-6$,
and in the section of "fiberwise exactness", we prove that it is enough to check
exactness of the quotient complex on fibers over the diagonal $\Delta \subset
X \times X .$

\heading {Exactness of sub-complexes and quotient complexes  }
\endheading

Let $D.$ be an arbitrary exact complex, and suppose $B. \subset D.$ is a sub-complex. Since
this sub-complex is not necessarily exact, we would like to find a way to  decide whether or
not $B.$ is exact. One possibility is by looking at the quotient-complex: $ Q.= D. / B. $.
 
\proclaim {Lemma 3.3}
$B.$ is exact up to degree $m$ if and only if $Q.$ is exact up to degree $m+1$.
\endproclaim

\demo {Proof} This can be easily proven by a diagram chase. Consider the following diagram:

$$
\alignat 9 & & 0 & & 0 & & 0 & & 0 & \quad \quad \quad \quad \quad \quad \quad \quad \quad
\quad \quad \quad
\quad \quad \quad\\
 & & \downarrow & & \downarrow & & \downarrow & & \downarrow &  \\ & ... \to & B_{i+2} & \to
& B_{i+1} & \to & B_i & \to & B_{i-1} &\to ... \\ & & \cap & & \cap & & \cap & & \cap &  \\
& ... \to & D_{i+2} & \to & D_{i+1} & \to & D_i & \to & D_{i-1} &\to ... \\ & & \downarrow &
& \downarrow & & \downarrow & & \downarrow &  \\ & ... \to & Q_{i+2} & \to & Q_{i+1} & \to &
Q_i & \to & Q_{i-1} &\to ... \\ & & \downarrow & & \downarrow & & \downarrow & & \downarrow
&  \\ & & 0 & & 0 & & 0 & & 0 &  
\endalignat
$$

Suppose $B.$ is exact up to degree m. If $ q \in ker \ Q_{i+1} \to Q_i$, for $ i \leq m$,
then we will show that $q=d \overset \sim \to q $, for some $\overset \sim \to q \in
Q_{i+2}. $ Denote by  $\overline c$ the image of $c \in D_j$ in $Q_j$. Then
$q= \overline c$, for some $c \in D_{i+1} $ . Since $ \overline {dc} = 0$, it follows that
$ dc=b \in B_i $ with $ db = 0 $. Using the fact that $B.$ is exact in degree $i$, there
exist an  element $b' \in B_{i+1}$ with $ b = db'$. Consider $c-b'$, then $d ( c-b' ) =0$
and since $C.$ is  exact, $c-b' = d \overset \sim \to c$ for some $\overset \sim \to c \in
C_{i+2}.$ Set $\overset \sim \to q = \overline {\overset \sim \to c}. $

Now suppose $Q.$ is exact up to degree $m+1$. If $b \in ker \ B_i \to B_{i-1}$, for $i \leq
m$, then we need to show that $b=d \overset \sim \to b $, for some $\overset \sim \to b \in
B_{i+1}$. Since $ d \overline c = \overline {dc}=0$, and since $Q.$ is exact up to degree
$m+1$, $\overline c=dq$ for some $q\in Q_{i+2}$. But $q=d \overset \sim \to c $,  for some $
\overset \sim \to c \in C_{i+2} $. Consider $c - d \overset \sim \to c $, then $ \overline
{c - d \overset \sim \to c }=0 $, therefore
$c - d \overset \sim \to c \in B_{i+1} $. Set $ \overset \sim \to b = c - d \overset \sim
\to c $, then $ d \overset \sim \to b = d (c - d \overset \sim \to c )=b-0=b $, as desired.
q.e.d.

\enddemo

\specialhead {Fiberwise exactness}
\endspecialhead
In this section we will show that it is enough to check the exactness of the 
quotient complex on the diagonal $\Delta \subset X \times X$. Since $Q.$ is a complex that
resolves a vector bundle, it is enough to check exactness fiberwise.
Suppose for  $ p \times q $, we can find an element $g$ in the isotropy group of
$p$,
$Iso_p \subset Sp(V)$,
 such that $g(q)=p$. Then the cohomology of the complex $ Q_{p \times q} . $ and the
cohomology of $Q_{g(p) \times g( q)} .=Q_{p \times p} .$, for $g \in Iso_p$, are just the
same, because $Q.$ is a complex of homogeneous vector bundles over $X \times X$ and 
$X \cong Sp(V) / P_X $ where $P_X$ is the parabolic $P_X = Iso_e $ and $e={e_1 \wedge e_2}$
for some fixed  symplectic basis
$\{ e_1, e_2 ... , e_{2n} \}$ of $X$. 

Suppose we can find, if not a $g$ , at least a family $g_\lambda$ of elements in $Iso_p$,
such that the limit of $g_\lambda ( q ) $ over $\lambda$ is equal to $p$. Then from the
Upper-Semi-Continuity Theorem for Cohomology,[H] III, Theorem 12.8, it follows, that the cohomology of $Q_{p
\times q} .$ and $Q_{p \times p} .$ are the same, that is $Q.$ is exact up to degree $m$ if it
is exact on the diagonal  up to degree $m$.

\proclaim {Lemma 3.4} If $ p \times q \in X \times X $, then we can use the isotropy group
of p to move q to p, i.e. there is a family $g_\lambda \in Iso_p$, such that $\underset
\lambda
\to {lim} \ g_\lambda (q) = p$.
\endproclaim

\demo {Proof} Let $p \times q $ be in $X \times X$. Then we may assume that $p=e_1 \wedge
e_2$, for some symplectic basis
$e_1, e_2, ..., e_{2n}$. Without loss of generality, we may even assume that
$p$ and $q$ fit into one of the following five cases:

If $p \cap q = \{ 0 \} $, then
\item {a)} $ q= e_{n+1} \wedge e_{n+2}$, or
\item {b)} $ q= e_{n+1} \wedge e_{n+3}$, or
\item {c)} $ q= e_{n+3} \wedge e_{n+4}$.

If $p \cap q = \Bbbk v $, for some vector $v \in V$, then we may assume that $v=e_1$, and
\item {d)} $ q= e_1 \wedge e_{n+2}$, or
\item {e)} $ q= e_1 \wedge e_{n+3}$

Let $\gamma = (1-\lambda) \lambda^{-1} $
\item {a)} Define
$$ g_\lambda = 
\left( \matrix
\lambda^{-1} & & & & &1-\lambda & & & & & \\
 & \lambda^{-1}  & & & &  & 1-\lambda  & & &  & \\ & & 1 & & & & & & & & \\ & & & \ddots & &
& & & & & & \\ & &  & & 1 & & & & &  & \\ & &  & & & \lambda & & &  & \\ & &  & & & &
\lambda  & & & \\ & &  & & &  & & 1 & & \\ & &  & & &  & & & \ddots & \\ & &  & & & & & & &
1 \\
\endmatrix \right).
$$  Then $g_1 (q) = e_{n+1} \wedge e_{n+2} = q$ and
$g_\lambda (q)=g_\lambda (e_{n+1} \wedge e_{n+2}) = ((1- \lambda ) e_1 + \lambda e_{n+1})
\wedge ((1- \lambda ) e_2 + \lambda e_{n+2})  $. Therefore the limit $ \underset { \lambda
\to 0 } \to {lim} \ g_\lambda (q) = e_1 \wedge e_2 $.

\newpage

\item {b)} Define $ g_\lambda = $
$$  
\left( \matrix
\lambda^{-1} & & & & & &1-\lambda & & & & & & \\
 & 1  & & & & & & & 1-\lambda & & &  & \\
 & & \lambda^{-1} & & & & & \gamma & & & & & \\
 & & & 1 & & & & & & & & & \\
 & & & & \ddots & & & & & & & & & \\
 & & & & & 1 & & & & & &  & \\  & &  & & & & \lambda & & & &  & \\
 & &  & & & && 1  & & & & \\
 & &  & & & & & & \lambda & & & \\
 & &  & & &  & & & & 1 &  &  \\ & &  & & &  & & & & & \ddots &  \\  & &  & & & & & & & & & 1
\\
\endmatrix \right).
$$  Then $g_1 (q)=e_{n+1} \wedge e_{n+3}=q$ and
$g_\lambda (q)=g_\lambda (e_{n+1} \wedge e_{n+3}) = ((1- \lambda ) e_1 + \lambda e_{n+1})
\wedge ((1- \lambda ) e_2 + \lambda e_{n+3})  $. Therefore the limit $ \underset { \lambda
\to 0 }
\to {lim} \ g_\lambda (q) = e_1 \wedge e_2 $.

\item {c)} Define $ g_\lambda = $
$$  
\left( \matrix 1 & & & & & & & &1-\lambda & & & &  \\
 & 1  & & & & & & & & 1-\lambda  & & &   \\
 & & \lambda^{-1} & & & &  \gamma & & & & & & \\
 & & & \lambda^{-1} & & & &  \gamma & & & & & \\
 & & & & \ddots & & & & & & & & & \\
 & & & & & 1 & & & & & &  & \\
  & &  & & & & 1 & & & &  & \\
 & &  & & & && 1  & & & & \\
 & &  & & & & & & \lambda & & & \\
 & &  & & &  & & & & \lambda &  &  \\  & &  & & &  & & & & & \ddots &  \\ 
 & &  & & & & & & & & & 1 \\
\endmatrix \right).
$$  Then $g_1 (q)=e_{n+3} \wedge e_{n+4}=q$ and
$g_\lambda (q)=g_\lambda (e_{n+3} \wedge e_{n+4}) = ((1- \lambda ) e_1 + \lambda e_{n+3})
\wedge ((1- \lambda ) e_2 + \lambda e_{n+4})  $. Therefore the limit $ \underset { \lambda
\to 0 }
\to {lim} \ g_\lambda (q) = e_1 \wedge e_2 $.

\newpage

\item {d)} Define
$$ g_\lambda = 
\left( \matrix 1 & & & & & & & & & & \\
 & \lambda^{-1}  & & & &  & 1-\lambda  & & &  & \\
 & & 1 & & & & & & & & \\
 & & & \ddots & & & & & & & & \\
 & &  & & 1 & & & & &  & \\  & &  & & & 1 & & &  & \\
 & &  & & & & \lambda  & & & \\  & &  & & &  & & 1 & & \\  & &  & & &  & & & \ddots & \\  &
&  & & & & & & & 1 \\
\endmatrix \right).
$$  Then $g_1 (q) = e_1 \wedge e_{n+2} = q$ and
$g_\lambda (q)=g_\lambda (e_1 \wedge e_{n+2}) =  e_1 
\wedge ((1- \lambda ) e_2 + \lambda e_{n+2})  $. Therefore the limit $ \underset {
\lambda \to 0 } \to {lim} \ g_\lambda (q) = e_1 \wedge e_2 $.

\item {e)} Define $ g_\lambda = $
$$  
\left( \matrix 1 & & & & & & & & & & & & \\
 & 1  & & & & & & & 1-\lambda & & &  & \\
 & & \lambda^{-1} & & & & & \gamma & & & & & \\
 & & & 1 & & & & & & & & & \\
 & & & & \ddots & & & & & & & & & \\
 & & & & & 1 & & & & & &  & \\  & &  & & & & 1 & & & &  & \\
 & &  & & & && 1  & & & & \\
 & &  & & & & & & \lambda & & & \\
 & &  & & &  & & & & 1 &  &  \\ & &  & & &  & & & & & \ddots &  \\  & &  & & & & & & & & & 1
\\
\endmatrix \right).
$$  Then $g_1 (q)=e_{n+1} \wedge e_{n+3}=q$ and
$g_\lambda (q)=g_\lambda (e_{n+1} \wedge e_{n+3}) = ((1- \lambda ) e_1 + \lambda e_{n+1})
\wedge ((1- \lambda ) e_2 + \lambda e_{n+3})  $. Therefore the limit $ \underset {
\lambda \to 0 }
\to {lim} \ g_\lambda (q) = e_1 \wedge e_2 $.

Left to check that all the above matrices are indeed symplectic. The symplectic form that we
are using, is given by:
$$
 <v,w>  = \ ^\tau v 
\left( \matrix  0 & I \\ -I & 0
\endmatrix
\right) w .
$$ It is easy to verify that the columns of all the above defined matrices are in fact  a
symplectic basis. q.e.d.
\enddemo

Before we finish the proof of the Main Theorem, we will discuss some exact sequences, chapter 4, and
apply these  to the K-theory of the symplectic grassmannian of isotropic two planes in four space,
chapter 5, and finally proof the Main Theorem, chapter 6.

\newpage

\topmatter
\title\chapter{4}{ Some exact sequences } \endtitle
\endtopmatter

In this chapter we will discuss some special exact sequences which play an important role
in proving  that the sub-complex
$B_.$ is exact. These sequences are a special case of  Lascoux Resolutions.
The techniques that we are applying in this chapter are the geometric techniques
of calculating syzygies, which go back to G. Kempf [K] and which were developed by 
P.Pragacz and J. Weyman [PW]. Similar methods can also be found in  
 J. Weyman's  book [W], under preparation, On the cohomology of vector bundles and syzygies, chapter
5 and 6.

The idea of using these complexes to proof the main theorem 
and the construction of these complexes is entirely due to J. Weyman. 
Nevertheless
the proof
here is different from  J. Weyman's proof,  since we are  interested only in a very special
case of  these sequences and for this special case we can proof the results more directly.

\head
Notations
\endhead

Let $T,T'$ be vector spaces of dimension $n$ and $m$ and set 
$$
\gather
X=Hom(T,T') \cong T^* \otimes T'\ , \\
  X_r \subset X , X_r = \{ \phi \in X | \text {rank} ( \phi ) \leq r\} \quad \text {and}\\ 
\overset \wedge \to X = X - X_r.
\endgather
$$

Note, that $X$ can also be described as the set of all $m \times n$ matrices, 
thus $X_r$ is the subset of
matrices with all $(r+1) \times (r+1)$ minors vanishing. This gives a description of $X_r$ 
as a
determinantel variety. A more general case of determinantel varieties can be found in
Chapter 6 [W].

Set $V= Gr(n - r,T )$, let $\Cal R $ be the tautological bundle and consider the product 
$V \times X$. Denote the projections of $ V \times X $ to $V$ and $X$ by
$p$ respectively $q$, i.e.
$p: V \times X \to V $ and $ q: V \times X \to X .$ Let 
$$
Z_r= \{ (l,\phi) \in V \times X| \  l \subset \text {ker} \phi \}. 
$$
Note, that $Z_r$ is just the pullback of $X_r$ under $q$.
Set $\overset \wedge \to Z = V \times X - Z_r $. 
Denote the projections of $ \overset \wedge \to Z $ to $V$ and $\overset \wedge \to X$ by
$\overset \wedge \to p$ respectively $\overset \wedge \to q$, that is:
$$
\alignat 5 & & & \ \ \   \overset \wedge \to Z \ \ \ & & \quad \quad \quad \quad \quad \quad
\quad
\quad \quad 
\quad 
\quad \quad \quad \quad \quad \quad
 \\ &  & \overset {\overset \wedge \to p} \to \swarrow & & 
\overset {\overset \wedge \to q} \to \searrow & \\ & V & & & &\overset \wedge \to X 
\endalignat
$$

 We will now apply the geometric techniques of calculating syzygies to define an exact
Koszul complex of $ \Cal O_{\overset \wedge \to Z} $~-modules,
which after twisting we will then push forward under the map $\overset \wedge \to q$ to get 
an exact complex of
 $\Cal O_{\overset \wedge \to X} $-modules.

\specialhead {A special Koszul complex}
\endspecialhead

View ${T'}^*$ as trivial bundle over $V$ and consider the bundle 
$\Cal R \otimes {T'}^*$ over $V$.
Define a map of vector bundles: 
$$
\gather
 \Psi : \overset \wedge \to p^* (\Cal R \otimes {T'}^* ) \to \Cal O_{\overset \wedge \to Z
},
\\
 \text { over a fiber} \ l \times
\psi :
\Psi ( r \otimes T ) = T ( \psi ( r) ). 
\endgather
$$

 $\Psi$ is surjective, because 
 for a point $( l , \psi ) \in
\overset \wedge \to Z, l$ is not contained in the kernel of $\psi$, thus
the image of $\Psi_{l \times \psi} = \Bbb C $.

Let $K_. =K_.(\overset \wedge \to p^* (\Cal R \otimes {T'}^* ) , \Psi )$ be the Koszul
complex of
$\overset \wedge \to p^* (\Cal R \otimes {T'}^* )$ and 
$\Psi :\overset \wedge \to p^* (\Cal R \otimes {T'}^* )
\to \Cal O_{\overset \wedge \to Z } $, see also chapter 1, Koszul complexes.

\proclaim {Proposition 4.1}
$$
K_.: 
0 \to \overset {top} \to \Lambda \overset \wedge \to  p^* (\Cal R \otimes {T'}^* ) 
\to ... 
\to \overset 2 \to \Lambda \overset \wedge \to p^* (\Cal R \otimes {T'}^* ) 
\to \overset \wedge \to  p^* (\Cal R \otimes {T'}^* )
\to   \Cal O_{\overset \wedge \to Z } \to 0
$$
is exact.
\endproclaim

\demo {Proof}
Lemma 1.2 gives a criteria for a Koszul complex to be exact.
Since  $ \Psi : \overset \wedge \to p^* (\Cal R \otimes {T'}^* ) \to \Cal O_{\overset \wedge
\to Z }  $ is surjective, we have to
check that $\overset \wedge \to p^* (\Cal R \otimes {T'}^* )$ has the right rank,  that
is:
$$
 \text { rank} (\overset \wedge \to p^* (\Cal R \otimes {T'}^* ) ) =
\text {dim} \overset \wedge \to Z = \text {dim} Z - \text {dim} Z_r \ .
$$

Fix $l \in V$ and consider the subspace $X_{l,r}$ of $X$, consisting of all 
$\varphi \in X=Hom(T,T') $ with $ l \subset ker \varphi $.
Recall that dim $l=n-r$, dim $T=n$ and dim $T' =m$. Thus
$X_{l,r}$ is a subspace of $X$ of codimension $ r  m$, and therefore
$\text {codim}_{V \times X}  Z_r = ( \text {dim} V ) \  rm = \text {dim} \overset \wedge \to
Z $. 

The rank of $\overset \wedge \to p^* (\Cal R \otimes {T'}^* )$
is equal to $ ( \text {dim} V ) \text {rank}  (\Cal R \otimes {T'}^*) = (\text {dim} V)
\ rm$ as desired.  q.e.d.
\enddemo

We will twist $ K_. $ and push this twisted Koszul complex forward under
the projection
$\overset \wedge \to  q: 
\overset \wedge \to Z \to \overset \wedge \to X$ .
Since there might be higher direct images we can not always do this. For our purposes it
will be enough to look at special cases where the higher direct images vanish
automatically. For a broader approach, see appendix or [W]. 

In Chapter 1, Lemma 1.5, we got a criteria to when the pushforward of a complex is exact.
In the above situation the hypothesis of Lemma 1.5 that all the higher direct images under
the  projection $\overset \wedge \to q:  \overset \wedge \to Z \to
\overset \wedge \to X$  of all
the twisted
$K_j(\Cal V)  $'s, $\Cal V $ a vector bundle on $V$, are
zero, are not satisfied. Here only certain direct images vanish,
but as we will see next, this weaker condition is enough to imply the exactness of the 
complex
$\overset \wedge \to q_* K_.(\Cal V)  $ .

\specialhead $ R^i q_*$ of a complex of coherent sheaves \endspecialhead

\proclaim {Lemma 4.2} Let $q: X \to Y$ be a projective morphism of noetherian schemes. Let 
$ 0 \to \Cal F_{top} @> d_{top} >> \to ...\to \Cal F_m @> d_m >>  ...
 @> d_1 >> \Cal F_0 \to 0
$ be an exact sequence of coherent $\Cal O_X$ - modules.  Suppose $R^i q_* \Cal F_{l+i} =0 $ for all
$i>0$ and all $0 \leq l \leq m$, then
$$
 q_* \Cal F_m \to ... \to q_* \Cal F_0 \to 0
$$ is an exact complex of coherent $\Cal O_Y$ - modules.
\endproclaim

\demo {Proof} We will follow the same lines as in the proof of Lemma 1.3. First we split the
complex $\Cal F.$ into short exact sequences,

$ 0 \to ker \ d_{m-k} \to \Cal F_{m-k} \to Im \ d_{m-k} \to 0 $. Next, in order to show that $q_*
\Cal F.$ is exact, we also split this complex into  sequences,

$ 0 \to ker \ q_* d_{m-k} \to q_* \Cal F_{m-k} \to q_* Im \ d_{m-k} \to 0 $, and show that  these are
exact for all $k$.

Again it suffices to show that $ R^i q_* ker \ d_{(m-k)} =0 $ for all $i>0, 0 \leq k \leq m$.. As
before we will use  induction over $k$, but now the  case of
$k=0$ is not obvious anymore.
 
\item {$k=0:$} Set $\Cal A = ker \ d_m$. We will show
$R^i q_* ker \ d_{(m+l)} \cong R^{i+1} q_* ker \ d_{(m+l+1)}$ for all $ l \geq 0 $,  which in turn
implies 
$$
\gather R^i q_* \Cal A \cong R^{i+1} q_* ker \ d_{m+1} \cong R^{i+2} q_* ker \ d_{m+2} 
\cong R^{i+3} q_* ker \ d_{m+3} \\ \cong ... \cong R^{top} q_* ker \ d_{top} =0.
\endgather
$$

\item {l=0:}
 $\Cal A \cong Im \ d_{m+1}$ and $\Cal A$ fits into the  short exact sequence:
$ 0 \to ker \ d_{m+1} \to \Cal F _{m+1} \to \Cal A \to 0 $. Consider the corresponding long exact
sequence of $Rq_*$:
$$ ...\to R^i q_* \Cal F_{m+1} \to R^i q_* \Cal A \to R^{i+1} q_* ker \ d_{m+1} \to R^{i+1} q_* \Cal
F_{m+1} \to ...
$$

Since the higher direct images of $ \Cal F_{m+1}$ are all zero, it follows that
$ R^i q_* \Cal A \cong R^{i+1} q_* ker \ d_{m+1}$.

\item{$l>0:$}  Consider  the short exact sequence 
$0 \to ker \ d_{m+l+1} \to \Cal F _{m+l+1} \to ker \ d_{m+l} \to 0 $, and consider the resulting long
exact sequence of $R q_*$:
$$
\gather ...\to R^i q_* \Cal F_{m+l+1} \to R^i q_*  ker \ d_{m+l} \\ \to R^{i+1} q_* ker \ d_{m+l+1}
\to R^{i+1} q_* \Cal F_{m+l+1} \to ...
\endgather
$$ Here the additional assumption of $R^i q_* \Cal F_{m+l+i}=0, 0 \leq l \leq m$ is needed. These
imply 
$R^i q_* \Cal F_{m+l+1}=0$ for
$ i \geq l+1$,  hence

$R^i q_*  ker \ d_{m+l} \cong R^{i+1} q_* ker \ d_{m+l+1}$ for $i \geq l+1 $. q.e.d.

\item{$k \rightsquigarrow k+1:$} This is just the same proof as in Lemma 1.5 , since the
only facts used here are
$R^i q_* \Cal A =0$ and $R^i q_* \Cal F_{j} =0$ for $0 \leq j \leq m$. q.e.d.
\enddemo

\specialhead Some exact sequences \endspecialhead

Recall  the notations. 
$\overset \wedge \to X = X - X_r, X=Hom(T,T')$, $V=Gr(n-r,T)$ and 
$\overset \wedge \to q$ is the projection
 $\overset \wedge \to q: \overset \wedge \to Z \to \overset \wedge \to X $.
 
Let $\Cal V $
be a vector bundle on $V $.  Denote by $K_.(\Cal V) $ the tensor product of the Koszul
complex
$K_. $ by the vector bundle $\Cal V$.
 
We will use Lemma 4.2 to show that the pushforward
  $\overset \wedge \to q_* K_.(\Cal V)$   under the projection $\overset \wedge \to q$ is
exact up to a certain degree, depending on the cohomology of the vector bundle $\Cal V$:
 
\proclaim {Corollary 4.3}
Let $\Cal V$ be a vector bundle on $V$ and suppose that 
$$
H^i(V, \overset i+l \to \Lambda ( \Cal R \otimes {T'}^*) \otimes \Cal V ) =0 \  \text {for } i > 0 
\text{ and all } 0 \leq l \leq m ,
$$
then $\overset \wedge \to q_* K_.(\Cal V)$ is exact up to degree $m-1$, i.e.
$$
\overset \wedge \to q_*(\overset m \to \Lambda p^* (\Cal R \otimes {T'}^* ) \otimes \Cal V) \to
 ... \to
\overset \wedge \to q_*
(\overset \wedge \to  p^* (\Cal R \otimes {T'}^* ) \otimes \Cal V) \to  \Gamma ( V, \Cal V)
\otimes \Cal O_{\overset \wedge \to X}
\to 0 
$$
is exact.
\endproclaim

\addline

\demo {Proof}
This is a special case of Lemma 4.2, thus we only need to check that
 $R^i \overset \wedge \to q_* K_{l+i} (\Cal V)=0$ for all $ i>0$ and all $ 0 \leq l \leq m$.
 Consider the commutative
diagram:
$$
\alignat 5
& & & \ \ \ \overset \wedge \to Z & & \quad \quad \quad \quad \quad \quad \quad \quad \quad 
\quad 
\quad \quad \quad \quad \quad \quad
 \\
&  & \overset {\overset \wedge \to p} \to \swarrow & & 
\overset {\overset \wedge \to q} \to \searrow & \\
& V & & & & \overset \wedge \to X \\
&  & \overset {\overset \wedge \to p'} \to \searrow & & 
\overset {\overset \wedge \to q'} \to \swarrow & \\
& & & \ Spec \ \Bbbk & &
\endalignat
$$

For all $0 \leq l \leq m$, the higher direct images of
$K_* ( \Cal V ) $ are:

$$
\gather
 R^i \overset \wedge \to q_* 
(\overset \wedge \to  p^* ( \overset i+l \to \Lambda  (\Cal R \otimes {T'}^* )
\otimes \Cal V) ) \cong  {\overset \wedge \to q'}^* 
R^i \overset \wedge \to p'_* (\overset i+l \to \Lambda( \Cal
R \otimes {T'}^* )
\otimes \Cal V) \cong \\
 H^i ( V, \overset i+l \to \Lambda  (\Cal R \otimes {T'}^* ) \otimes \Cal V) \otimes 
\Cal O_{\overset \wedge \to X} = 
\cases
 H^0 ( V, \overset l \to \Lambda  (\Cal R \otimes {T'}^* ) \otimes \Cal V) \otimes 
\Cal O_{\overset \wedge \to X}
& \text {if} \ i=0 
\\  0 & \text { otherwise} \ .
\endcases 
\endgather 
$$
Thus the higher direct images that are supposed to vanish, are zero, and
by Lemma 4.2 $\overset \wedge \to q_* K_.(\Cal V)$ is exact up to degree $m-1$. q.e.d.
\enddemo

Under the same assumption of the Corollary, note that
 $ \overset \wedge \to  q_*(\overset l \to \Lambda 
\overset \wedge \to p^* (\Cal R \otimes {T'}^* ) \otimes
\Cal V) 
\cong  H^0 ( V, \overset l \to \Lambda  (\Cal R \otimes {T'}^* ) \otimes \Cal V) \otimes 
\Cal O_{\overset \wedge \to X} $.
Thus
$\overset \wedge \to q_* K_.(\Cal V)$ is given by:

$$
\gather
 H^0 ( V, \overset m \to \Lambda  (\Cal R \otimes {T'}^* ) \otimes \Cal V) \otimes 
\Cal O_{\overset \wedge \to X} \to
... \to  H^0 ( V, \overset 2\to \Lambda  (\Cal R \otimes {T'}^* ) \otimes \Cal V) \otimes 
\Cal O_{\overset \wedge \to X} \\
\to  H^0 ( V,   \Cal R \otimes {T'}^*  \otimes \Cal V) \otimes 
\Cal O_{\overset \wedge \to X}
\to  H^0 ( V,  \Cal V) \otimes \Cal O_{\overset \wedge \to X} \to 0,
\endgather
$$
which is exact by the previous Corollary.

Next we consider the case, when $V$ is a projective space, to show that the sub-complex is 
exact. In particular, in this case $\Cal V = \Cal O_{\Bbb P^{n-1}}(j)$ for some $j$.

\specialhead {Two special sequences}
\endspecialhead

\proclaim {Corollary 4.4}
Suppose $V\cong \Bbb P ( T ) \cong \Bbb P^{n-1}$ and let $\Cal V = \Cal O_{\Bbb P^{n-1}}(j)$
for $ j \geq 0 $.
Then
$$
\gather
 \overset j \to \Lambda {T'}^* \otimes \Cal O_{\overset \wedge \to X}
\to T^* \otimes \overset j-1 \to \Lambda {T'}^* \otimes 
\Cal O_{\overset \wedge \to X}
\to Sym_2(T^*)  \otimes \overset j-2 \to \Lambda {T'}^* \otimes 
\Cal O_{\overset \wedge \to X} \to ...\tag {a} \\ \to
Sym_{j-1}(T^*)  \otimes {T'}^* \otimes \Cal O_{\overset \wedge \to X}  \to
Sym_j(T^*)  \otimes \Cal O_{\overset \wedge \to X} \to 0 \\
\text {is exact.}
\endgather 
$$
\item{(b)}
Over a point $\varphi = s \otimes f \in {\overset \wedge \to X}
\subset X = Hom(T,T') \cong T^*
\otimes T', $ the maps of  this
complex are given by:
$$
\gather
d_\varphi : Sym_k (T^*) \otimes \overset j-k \to \Lambda {T'}^* \otimes \Bbbk (\varphi) \to
 Sym_{k+1 }(T^*) \otimes \overset j-(k +1) \to \Lambda {T'}^* \otimes \Bbbk (\varphi) \\
d_\varphi ( ( s_1 ... s_k ) \otimes ( t_1 \wedge ... \wedge t_{j-k} ) )
= s_1 ... s_k s \otimes \underset {i \geq 0} \to \Sigma (-1)^i f(t_i) \
t_1 \wedge ... \wedge \overset \wedge \to t_i \wedge ... \wedge t_{j-k}
\endgather
$$ 

\endproclaim

\addline

\demo {Proof}
\item{(a)}
This is a very special case of the previous Corollary. We need to check that the
hypothesis are satisfied.
The tautological bundle $\Cal R$ in this case is equal to $\Cal O_{\Bbb P^{n-1}} (-1)$.
The higher cohomology,
  $$
\gather
 H^i( \overset i+l \to \Lambda ( \Cal R \otimes {T'}^*) \otimes \Cal V ) 
=H^i( \overset i+l \to \Lambda ( \Cal O_{\Bbb P^{n-1}} (-1) \otimes {T'}^*) \otimes 
\Cal O_{\Bbb P^{n-1}}(j) ) \\
\cong
H^i(  \Cal O_{\Bbb P^{n-1}} (j-(i+l)) \otimes  \overset i+l \to \Lambda {T'}^* 
=0 \quad \   \text {for } i > 0 
\text{ and all } 0 \leq l \leq m .
\endgather
$$
Thus the hypothesis of Corollary 4.2 are satisfied and 
for $i=0$, the sections are given by:
$$
\gather
H^0( \overset l \to \Lambda ( \Cal O_{\Bbb P^{n-1}} (-1) \otimes {T'}^*) \otimes 
\Cal O_{\Bbb P^{n-1}}(j) ) \cong
H^0(  \Cal O_{\Bbb P^{n-1}} (j-l) \otimes  \overset l \to \Lambda {T'}^* ) \\ \cong
 Sym_{j-l}(T^*)  \otimes \overset l \to \Lambda {T'}^* 
\endgather
$$
\item{(b)}
We need to trace the map through the isomorphism of
$$
\gather
H^0( \overset {j-k} \to \Lambda ( \Cal O_{\Bbb P^{n-1}} (-1) \otimes {T'}^*) \otimes 
\Cal O_{\Bbb P^{n-1}}(j) ) \\
\text{with} \ \
Sym_{k}(T^*)  \otimes \overset {j-k} \to \Lambda {T'}^* 
\endgather
$$
Set $ m=j-k$. Let us start with the map $d$ of the Koszul complex $K_.(\Cal V ) $:

$$
\gather
d: \overset m \to \Lambda ( \Cal O_{\Bbb P^{n-1}} (-1) \otimes {T'}^*) \to
\overset {m-1} \to \Lambda ( \Cal O_{\Bbb P^{n-1}} (-1) \otimes {T'}^*) \\
d(  (r_1 \otimes T_1 ) \wedge ... \wedge (r_m \otimes T_m ) )  \\ =
\underset {i \geq 0} \to \Sigma (-1)^i \ s(r_i) \ f(T_i) 
(r_1 \otimes T_1 ) \wedge ... \wedge \overset \wedge \to {(r_i \otimes T_i ) } \wedge ... 
\wedge (r_m \otimes T_m ).
\endgather
$$

Use the isomorphism of
$
\overset l \to \Lambda ( \Cal O_{\Bbb P^{n-1}} (-1) \otimes {T'}^*) $ with
$
\Cal O_{\Bbb P^{n-1}} (-l) \otimes \overset l \to \Lambda (  {T'}^*).
$
Then the map $d$ is given by:
$$
\gather
 d:
\Cal O_{\Bbb P^{n-1}} (-m) \otimes \overset m \to \Lambda (  {T'}^*) \otimes \Bbbk(\varphi) \to
\Cal O_{\Bbb P^{n-1}} (-(m-1)) \otimes \overset {m-1} \to \Lambda (  {T'}^*) \otimes \Bbbk(\varphi) 
\\
 \ d(  (r_1 ... r_m) \otimes (T_1  \wedge ... \wedge  T_m ) )  \\ =
\underset {i,j \geq 0} \to \Sigma (-1)^{i+j} 
\ s(r_j) \ f(T_i) \ ( r_1 ...\overset \wedge \to r_j ... r_m ) \otimes (T_1  \wedge ... \wedge
\overset \wedge \to  T_i   \wedge ... 
\wedge  T_m ) \\
= r_1 ... r_m s \ \otimes \  
\underset {i \geq 0} \to \Sigma (-1)^{i} 
 \ f(T_i) \  (T_1  \wedge ... \wedge
\overset \wedge \to  T_i   \wedge ... 
\wedge  T_m )
.
\endgather
$$
The last equality follows from the isomorphism: 
$
 \Cal O_{\Bbb P^{n-1}} \cong  \Cal O_{\Bbb P^{n-1}} (-1) \otimes \Cal O_{\Bbb P^{n-1}} (1) 
$,
which sends $s$ to $\Sigma s(t_i) \otimes t^i  $ for a basis $t_1 ... t_n$ of $T$ and its
dual basis $t^1 ... t^n$.

This  implies that the map of $\overset \wedge \to q_* K_. ( \Cal V ) $ is given by:
$$
\gather
\text  d:
\Cal O_{\Bbb P^{n-1}} (k) \otimes \overset m \to \Lambda (  {T'}^*) \otimes \Bbbk(\varphi) \to
\Cal O_{\Bbb P^{n-1}} (k+1)) \otimes \overset {m-1} \to \Lambda (  {T'}^*) \otimes \Bbbk(\varphi) 
\\
 d(  (s_1 ... s_k) \otimes (T_1  \wedge ... \wedge  T_m ) )  \\ =
 s_1 ... s_k s \ \otimes \  
\underset {i \geq 0} \to \Sigma (-1)^{i} 
 \ f(T_i) \  (T_1  \wedge ... \wedge
\overset \wedge \to  T_i   \wedge ... 
\wedge  T_m ) \ .
\endgather
$$
q.e.d.
\enddemo

Tensor this  complex of the  Corollary 4.4 by $\overset {top} \to \Lambda (T') $, 
then this exact sequence is 
the  sequence that we will use in the proof of the Main Theorem 3.2.

\proclaim {Corollary 4.5} 
Let $V\cong \Bbb  P ( T ) \cong \Bbb P^{n-1}$ and set $\Cal V = \Cal O_{\Bbb P^{n-1}}(j)$
for $ j
\geq 0 $.
$$
\gather
 \overset m-j \to \Lambda {T'} \otimes 
\Cal O_{\overset \wedge \to X}
\to 
T^* \otimes \overset m-j+1 \to
   \Lambda {T'} \otimes 
\Cal O_{\overset \wedge \to X}
\to Sym_2(T^*)  \otimes \overset m-j+2 \to \Lambda {T'} \otimes 
\Cal O_{\overset \wedge \to X} \to ...\tag{a} \\ \to
Sym_{j-1}(T^*)  \otimes \overset m-1 \to \Lambda{T'} \otimes 
\Cal O_{\overset \wedge \to X}  \to 
Sym_j(T^*)  \otimes \overset m \to \Lambda T' \otimes
\Cal O_{\overset \wedge \to X} \to 0
\endgather 
$$
is exact.
\item{(b)} Over a point $\varphi = s \otimes f \in 
\overset \wedge \to X \subset X = Hom(T,T') \cong T^* \otimes T', $ the maps of
this complex are given by:
$$
\gather d_\varphi : Sym_k (T^*) \otimes \overset m-j+k \to \Lambda {T'}^* \otimes \Bbbk (\varphi)
\to
 Sym_{k+1 }(T^*) \otimes \overset m-j+{k +1} \to \Lambda {T'}^* \otimes \Bbbk (\varphi) \\ d_\varphi (
( s_1 ... s_k ) \otimes ( f_1 \wedge ... \wedge f_{m-j+k} ) ) 
= s_1 ... s_k s \otimes 
f \wedge f_1 \wedge ...  \wedge f_{m-j+k}
\endgather
$$ 

\endproclaim

\demo {Proof}
This follows from:
$$
\overset l \to \Lambda {T'}^* \otimes \overset top \to \Lambda T'
\cong \overset top - l \to \Lambda T'
$$ q.e.d.
\enddemo
We can now  proof the Main Theorem, but first
we discuss the  example of the symplectic grassmannian of 2-planes in four space,
$$
X=SpGr(2,4) \ .
$$

\newpage

\topmatter
\title\chapter{5} The  Symplectic Grassmannian X of two
dimensional planes in four dimensional space \endtitle
\endtopmatter

Denote by $X$ the symplectic grassmannian of two-planes in the
four dimensional vector space $V$ over an algebraically closed
field $\Bbbk$ of characteristic 0. We keep the notations of
chapter 2 and 3.

In this chapter we will study the K-theory of $X$. Since $X$ is
isomorphic to a quadric in $\Bbb P^4 $ this case has been proven
by Kapranov in [KII] and by Swan in [S] using the Clifford algebra
of $ \overset 2 \to\Lambda V $.

Our approach here is different, since we discuss $X$  as a
symplectic  grassmannian inside a grassmannian, and not as a
quadric inside a projective space. Therefore we also get a
different generating system for the K-theory of $X$, but overall
all of these approaches give similar algebraic descriptions of the K-theory.

Similar to the case of ordinary grassmannians, chapter 1, we will
construct a finite resolution of 
$\Cal O_{\Delta} $ and use this to show that 
$$
\Cal X= \{ \Cal O_X (-2),\Cal O_X (-1), S, \Cal O_X \}
$$ is a generating system  for the K-theory of $X$, that is: 

\proclaim{Theorem 5.1}
$$ K_* ( SpGr(2,4) ) \cong \underset {\Cal G \in X} \to \bigoplus
K_*( \Bbbk )
$$
\endproclaim

\demo {Proof} The idea of the proof is  similar to the proof in
the  case of   ordinary grassmannians. We proof this in two
steps. First we construct the resolution, then we use some
cohomology computations in order to show that $\Cal X$ is a 
generating system.
\enddemo

\specialhead The resolution
\endspecialhead

Consider the Tate construction in this case:
$$
\matrix & &  & Sym_2 S \boxtimes \overset 2 \to \Lambda S^\perp & 
\\ & & \nearrow & & \searrow \\ &  \to  \Sigma^{(2,2)} S \boxtimes
\Phi_+ \to \Sigma^{(2,1)} S \boxtimes \Phi_- & & \oplus & & S
\boxtimes S^\perp \to \Cal O_{X \times X} \\ & & \searrow  & &
\nearrow \\ & &  & \overset 2 \to \Lambda S \boxtimes \Psi  &  \\
\endmatrix
$$ and consider the sub-complex:
$$
\overset 2 \to \Lambda S \boxtimes \Psi \to S \boxtimes S^\perp
\to \Cal O_{X \times X}
$$ In this section  we show that this sub-complex is exact and
furthermore that the kernel is given by
$ ( \overset 2 \to \Lambda S )^{\otimes 2} \boxtimes \overset 2
\to \Lambda S^\perp 
\cong \Cal O_{X \times X} (-2,-1) $, that is:

\proclaim {Theorem 5.2}
$$ 0 \to ( \overset 2 \to \Lambda S )^{\otimes 2} \boxtimes
\overset 2 \to \Lambda S^\perp \to
\overset 2 \to \Lambda S \boxtimes \Psi \to S \boxtimes S^\perp
\to \Cal O_{X \times X}
$$ is a finite resolution of $\Cal O_\Delta$.
\endproclaim

\demo {Proof} We split this proof into two parts. We show first
that the sub-complex is exact, then we show that the kernel is
equal to
$( \overset 2 \to \Lambda S )^{\otimes 2} \boxtimes \overset 2 \to
\Lambda S^\perp$.
\enddemo

 \head{The sub-complex}
\endhead

Let us represent the representations of $S$ and $S^\perp$ by Young
diagrams, i.e. let us represent $\Sigma^\alpha S \boxtimes
\Sigma^\beta S^\perp $ by the external tensor product of the Young
diagrams corresponding to the weights $\alpha$ and $\beta$. For
example:
$$
\gather
\young{2} \boxtimes \yng{1,1} = Sym_2 S \boxtimes \overset 2 \to
\Lambda S^\perp 
\endgather
$$ Over an open set $U$, where all the $\Psi_\alpha$'s are trivial
extensions,  that is $\Psi_\alpha$ is a direct sum over all
$\Sigma^{\alpha-2k} S^\perp |_{U}, k \geq 0 $,  consider the
resolution $D.$ of $\Cal O_{\Delta}$ from chapter 2:

$$
  \alignat 8  ...\to \yng{2,2} & \boxtimes \yng{2,2} & &
\longrightarrow & \yng{2,1} & \boxtimes 
\yng{2,1} & &\longrightarrow &
  \young{2} & \boxtimes \yng{1,1} & &  & 
\\ & \oplus & & \nearrow &   & \oplus && \nearrow && \oplus  &
\searrow & &
\\ ... \to \yng{2,2} & \boxtimes \young{2}  & & \longrightarrow &
\yng{2,1} & \boxtimes 
\young{1} & & \longrightarrow &
\yng{1,1} & \boxtimes \young{2} & \longrightarrow & \young{1}
\boxtimes 
\young{1} \longrightarrow ^.
 \boxtimes ^.  
\\ & \oplus && \nearrow && &&  \searrow && \oplus & \nearrow &
\\ ...\to \yng{2,2} & \boxtimes \ ^. &{}& &{}& & {} & & \yng{1,1}
& \boxtimes \
^.                                                                                         
\
\\ & \oplus && \nearrow
\\ ...\to \yng{3,1} & \boxtimes \yng{1,1} 
\endalignat
$$ Recall that in order to show that the sub-complex is exact up
to degree 1, it is enough to show that the quotient complex is
exact up to degree 2, Lemma 3.3 . In addition it is enough to
show that the quotient complex is exact over points on the
diagonal $\Delta \subset X \times X$, Lemma 3.4.

Consider the quotient complex over the open set $U$:
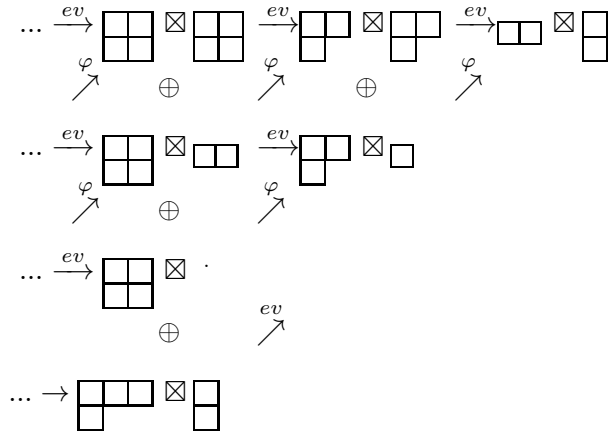
$$
  \alignat 8  ...@> ev >> \yng{2,2} & \boxtimes \yng{2,2} & & @>
ev >> & \yng{2,1} & \boxtimes 
\yng{2,1} & & @> ev >> & 
  \young{2} & \boxtimes \yng{1,1} & &  & \quad \quad \quad \quad
\quad \quad \quad
\\ \overset \varphi \to \nearrow \quad \quad  & \oplus & & 
\overset \varphi \to \nearrow &   & \oplus &&  \overset \varphi
\to \nearrow &&  & &
\\ ... @> ev >> \yng{2,2} & \boxtimes \young{2}  & & @> ev >> &
\yng{2,1} & \boxtimes 
\young{1} & &  
\\ \overset \varphi \to \nearrow \quad \quad & \oplus && 
\overset \varphi \to \nearrow && && &&
\\ ...@> ev >> \yng{2,2} & \boxtimes \ ^. &{}& &{}& & {} &
&
\\ & \oplus && \overset ev
\to  \nearrow
\\ ...\to \yng{3,1} & \boxtimes \yng{1,1}                                       
\endalignat
$$

Over the open set $U$, the map $d$ of the  quotient complex,
splits into two kinds of maps, $ev$ and $\varphi$. Let us discuss
these maps.

\subhead
The maps $ev$ and $\varphi$
\endsubhead

Let $(W_1, W_2)$ be a point in $X \times X$ and let
$gP$ respectively $hP$ be the parabolic leaving $W_1$
respectively $W_2$ invariant.

Recall the definitions of $ T^i_{(g,h)} $ and $ t_{(g,h)} $.

For $g,h \in G , i > 0 $, let
$$
\gather
  T_{(g,h)} = (g, v_1 \wedge v_2 ) \boxtimes (h,\eta) \ ,
 T^i_{(g,h)} = (g, (v_1 \wedge v_2)^{\otimes i } ) 
\boxtimes (h,\eta^i)  \\ \text {and} \ t_{(g,h)}= d T_{(g,h)} \ .
\endgather
$$ Recall that an element in
$\Sigma^\alpha S \boxtimes \Psi_{\alpha^*} $ over the point
$(gP,hP)$  can be  written uniquely as a sum of elements of the
form
$$ w T^i_{(g,h)} \ , \ w \in (\Sigma^{(\alpha_1-i,\alpha_2-i)}
S  )_{gP} \boxtimes  (\Sigma^{(\alpha^*-2i)} S^\perp) _{hP} \
.
$$ The map $d$ is defined as:
$$ d ( w T^i_{(g,h)})=  ( dw) T^i_{(g,h)} + (-1)^{ | w| }i w \cdot
t_{(g,h)} T^{i-1}_{(g,h)} \ .
$$ This translates to the two maps $ev $ and $\varphi $ as
follows.

Let $\beta= ( \alpha_1 - i, \alpha_2 -i)^*$.

$$
\alignat 8  & & & \Sigma^{\alpha-1} W_1 \ \boxtimes \  
\Sigma^{\beta-1}(V/W_2)^* & 
\hskip 1 in
\\ & & \overset ev \to \nearrow \\ & d: \Sigma^\alpha W_1 \
\boxtimes \   \Sigma^\beta (V/W_2)^*  & & \hskip .4 in\oplus
\\ & & \overset \varphi \to \searrow \\ & & & \Sigma^{\alpha-1}
W_1
\ \boxtimes \   
\Sigma^{\beta+1 }(V/W_2)^*
\endalignat
$$ The  map $ev$ is induced by the  evaluation map or
$$
\gather  ev ( (w_1 \wedge w_2)^{\otimes i} \otimes w )  =  (w_1
\wedge w_2)^{\otimes i} \otimes d w
\endgather
$$ and 
$$
\varphi ((w_1 \wedge w_2)^{\otimes i} \otimes w) =(-1)^{ | w| }i 
(w_1 \wedge w_2)^{\otimes i-1 }( w \cdot t_{(g,h)}) \ .
$$

Over a point $ W \times W$ on the diagonal $\Delta$, all the
evaluation maps $ev$ are zero, since for $v \in W$ and $f \in
(V/W)^*$,
$f(v)=0$. Therefore the quotient complex splits into several
subsequences. 
$$
  \alignat 8
  &  \text {deg 4} & &  &  & \text {deg 3}
 & & &  & \text {deg 2}
   &  & &   \\
\\ 
 \yng{2,2} & \boxtimes \yng{2,2} & &  & \yng{2,1} & \boxtimes 
\yng{2,1} & & & 
  \young{2} & \boxtimes \yng{1,1} & &  & \quad \quad \quad \quad
\quad \quad \quad
\\ ...\quad \overset \varphi \to \nearrow \quad \quad \quad  &
\oplus & & \quad \overset \varphi \to
\nearrow
\quad &   &
\oplus && \quad \overset \varphi \to \nearrow \quad &&  & & \\
\\  \yng{2,2} & \boxtimes \young{2}  & &  & \yng{2,1} & \boxtimes 
\young{1} & &  
\\ ... \quad \overset \varphi \to \nearrow \quad \quad \quad &
\oplus && \quad \overset \varphi \to 
\nearrow
\quad && && && \\
\\  \yng{2,2} & \boxtimes \ ^. &{}& &{}& & {} &
&
\\ & \dots                                      
\endalignat
$$

Hence to show that the quotient complex is exact in degree 2, it
is enough to show that over a point 
$W \times W \in \Delta$:
$$
\alignat 1
\yng{2,1} \boxtimes \young{1} @> \varphi >> \young{2} \boxtimes
\yng{1,1} \quad \text{is surjective.}
\quad
\quad
\quad \quad \quad \quad \quad \quad
\endalignat
$$

This on the other hand is a special case of our sequences in
chapter 4, Corollary 4.5. Recall the definition of $Y$:

\subhead { Some exact sequences }
\endsubhead

Let $W \times W$ be a point on the diagonal inside of $X \times X$.
Let 
$Y= Hom(T,T')$, 
$Y_r=\{ x \in Y \ |\ \text{rank} (x) \leq r \} $ and
$\overset \wedge \to Y = Y - Y_r $. Set $V=Gr(n-r,T) $,
 $T=W^*$ and $T'=(V/W)^*$.

Then 
$Y=Hom(W^*,(V/W)^* ). $ Consider the  case when  $r=1$, that is
when $V$ is equal to the projective space 
 $ \Bbb P (W^*) \cong \Bbb P^1  $.

\proclaim {Lemma 5.3} Let  $\Cal V = \Cal O_{\Bbb P^1 } (2)$ , then
$$
\Cal O_{\overset \wedge \to Y} @> d >>  W \otimes (V/W)^* \otimes
\Cal O_{\overset \wedge \to Y} @> d >>  Sym_2 (W)
\otimes \overset 2 \to 
\Lambda (V/W)^* \otimes \Cal O_{\overset \wedge \to Y}  \to 0
$$ is exact and over a point $x \in {\overset \wedge \to Y}, d_x$
is given by
$$ d_x(w_1 w_2 ... w_m \otimes f_1 \wedge ... \wedge f_l)=w_1 w_2
... w_m x \wedge f_1 \wedge ...\wedge f_l
$$
\endproclaim

\demo{Proof}

Recall the statement of Corollary 4.5, using
 $\Cal V = \Cal O_\Bbb P (j)$:

$$
\gather
 \overset m-j \to \Lambda {T'} \otimes 
\Cal O_{\overset \wedge \to Y}\to T^* \otimes \overset m-j+1 \to
\Lambda {T'} \otimes 
\Cal O_{\overset \wedge \to Y}
\to Sym_2(T^*)  \otimes \overset m-j+2 \to \Lambda {T'} \otimes
\Cal O_{\overset \wedge \to Y} \to
\\ \dots \to Sym_{j-1}(T^*)  \otimes \overset m-1 \to \Lambda{T'}
\otimes \Cal O_{\overset \wedge
\to Y}  \to  Sym_j(T^*) 
\otimes
\overset m \to \Lambda T' \otimes
\Cal O_{\overset \wedge \to Y}  \to 0\\
\text {is exact}.
\endgather 
$$Over a point $x = s \otimes f \in \overset \wedge \to Y
\subset Y = Hom(T,T') \cong T^* \otimes T', $ the maps of this
complex are given by:
$$
\gather d_x : Sym_k (T^*) \otimes \overset m-j+k \to \Lambda
{T'}^* \otimes \Bbbk (\varphi)
\to
 Sym_{k+1 }(T^*) \otimes \overset m-j+{k +1} \to \Lambda {T'}^*
\otimes \Bbbk (\varphi) \\ d_x ( ( s_1 ... s_k ) \otimes ( f_1
\wedge ... \wedge f_{m-j+k} ) )  = s_1 ... s_k \cdot  s \otimes  f
\wedge f_1
\wedge ...  \wedge f_{m-j+k}
\endgather
$$

In particular for $\Cal V = \Cal O_{\Bbb P^1 } (2)$ , we get the
exact complex
$$
\overset 2 \to \Lambda (V/W) \otimes \Cal O_{\overset \wedge \to
Y} \to W \otimes V/W \otimes \Cal O_{\overset \wedge \to Y} \to
Sym_2 (W)
\otimes \Cal O_{\overset \wedge \to Y} \to 0 \ .
$$  Twist this last complex by $\overset 2 \to \Lambda (V/W)^* $,
then
$$
\Cal O_{\overset \wedge \to Y} @> d >>  W \otimes (V/W)^* \otimes
\Cal O_{\overset \wedge \to Y} @> d >>  Sym_2 (W)
\otimes \overset 2 \to \Lambda (V/W)^* \otimes \Cal O_{\overset
\wedge \to Y}  \to 0
$$ is exact.

 q.e.d.
\enddemo

\subhead {Exactness of the sub-complex }
\endsubhead

Choose $g \in G$ in such a way that the parabolic $gP$   preserves
the space $W$. We will set $x = t $ with $t=t_{(g,g)} $  in the
above sequence.

Let $e_1, ... ,e_4$ be a symplectic basis of $V$, and $y_1, ...
,y_4$ the corresponding dual basis of $V^*$, such that $W=e_1
\wedge e_2$ and $ g \cdot \eta = y_1 \otimes y_3 + y_2
\otimes y_4 $. 

Recall that $\eta$ is defined as
$\eta =x_1 \otimes x_3 + x_2
\otimes x_4$.

Denote by $d$ the map from the Tate construction. Then
$$
\gather t= dT_{(g,g)} \\ = d ((g, v_1 \wedge v_2) \boxtimes (g,
\eta) ) \\ =  d ( (e_1 \wedge e_2) \boxtimes ( y_1 \otimes y_3 +
y_2
\otimes y_4 ) \\  = y_1 (e_1)e_2 \boxtimes y_3 - y_1 (e_2)e_1
\boxtimes y_3 + y_2(e_1) e_2 \boxtimes y_4 - y_2(e_2) e_1
\boxtimes y_4 \\ =e_2 \boxtimes y_3 - e_1 \boxtimes y_4 \ .
\endgather
$$ Set $x=t=e_2 \boxtimes y_3 - e_1 \boxtimes y_4$. Then $x \in Y$
and the rank of $x$ is equal to $2$, because $x(y_1)=-y_4 ,
x(y_2)=y_3  $ and $y_3, y_4$ are linearly independent. Thus $ x
\in \overset \wedge \to Y $ , because $x \notin Y_1$.

Now evaluate the complex of Lemma 4.3 at $ x \in \overset \wedge
\to Y$:
$$
\gather
\Cal O_{\overset \wedge \to Y,x} @> d_x >>  W \otimes (V/W)^*
\otimes \Cal O_{\overset \wedge \to Y,x} \\ @> d_x >>  Sym_2 (W)
\otimes \overset 2 \to \Lambda (V/W)^* \otimes \Cal O_{\overset
\wedge \to Y,x} 
\to 0
\endgather
$$ In terms of Young diagrams:
$$
\gather ^. \boxtimes ^. @> d_x >>    \young{1} \boxtimes \young{1} @>
d_x >> \young{2} \boxtimes
\yng{1,1} \to 0 
\endgather
$$

This is almost the sequence we want, except for a twist. Recall
the definition of  $\varphi$:
$$
\gather
 \yng{2,1} \boxtimes \young{1} @> \varphi >> \young{2} \boxtimes
\yng{1,1} \\ \\ 
\varphi ((e_1 \wedge e_2) \otimes v \boxtimes f  ) = -(v \boxtimes
f ) \cdot t =  d_x ( v \boxtimes f )
\endgather 
$$ Hence the difference between $ \varphi $ and $ d_x $ is just a
twist by 
$\Cal O_{\overset \wedge \to Y} (-1)$, which is trivial over the
affine variety $\overset \wedge
\to Y$.

Therefore
$$
\gather
\yng{2,2} \boxtimes ^. @> \varphi >>    \yng{2,1} \boxtimes
\young{1} @> \varphi >> \young{2} \boxtimes
\yng{1,1} \to 0 
\endgather
$$ is exact, which is just what we needed in order to show that
the sub-complex is exact. 

q.e.d.

\head The Kernel:
\endhead

First we define a map $ \vartheta$ from 
$ K=  ( \overset 2 \to \Lambda S )^{\otimes 2} \boxtimes
\overset 2 \to \Lambda S^\perp 
\cong \Cal O_{X \times X} (-2,-1) $ to the sub-complex and 
 show that this map is well defined and nonzero.  Then we
discuss the quotient complex and show that the kernel of the
sub-complex is indeed isomorphic to $K$, which then implies the
exactness of 
$$ 0 \to ( \overset 2 \to \Lambda S )^{\otimes 2} \boxtimes
\overset 2 
\to \Lambda S^\perp @> \vartheta >> 
\overset 2 \to \Lambda S \boxtimes \Psi @> d >> S \boxtimes
S^\perp @> d >> 
\Cal O_{X \times X} \ .
$$

\subhead { Definition of the map $\vartheta$}
\endsubhead

Define
$$
\gather
\vartheta : ( \overset 2 \to \Lambda S )^{\otimes 2} \boxtimes
\overset 2 \to \Lambda S^\perp
\to \overset 2 \to \Lambda S \boxtimes (V^* \otimes S^\perp) \\
\vartheta (s_1 \wedge s_2 \ \otimes \ v \wedge w \ \boxtimes \ f
\wedge g) \\ = s_1 \wedge s_2 \ \boxtimes \ 
\left( \matrix
 \ [ g(w) <v,\ > - g(v) <w,\ >] \otimes f \\ - [ f(w) <v,\ > -
f(v) <w,\ >] \otimes g 
\endmatrix
\right) \ .
\endgather
$$

This is certainly well defined and since $v,w,s_1$ and $s_2$ are
all in $S_{W_1}=W_1$, it follows that
$d(\vartheta (s_1 \wedge s_2 \ \otimes \ v \wedge w \ \boxtimes \
f \wedge g)) = 0 $ in 
$S \boxtimes S^\perp $.  Set $$
\rho =\left( \matrix
 \ [ g(w) <v,\ > - g(v) <w,\ >] \otimes f \\ - [ f(w) <v,\ > -
f(v) <w,\ >] \otimes g 
\endmatrix
\right) \ 
$$

We  need to check that the Image of $\vartheta$ is  contained in
$\overset 2 \to \Lambda S \boxtimes \Psi . $

\proclaim{Remark} Let us discuss an alternative definition of the
extension $\Psi$. Recall that $\Psi \subset V^* \otimes S^\perp$,
 Lemma 2.3, Remark (1), 
is an extension of
$Sym_2 S^\perp $ and $\Cal O_X $. Consider the exact sequence:
$$ 0 \to Sym_2 S^\perp \to V^* \otimes S^\perp @> \psi >>
\overset 2 \to \Lambda V^*
$$ Set $ \delta= (x_1 \wedge x_3 + x_2 \wedge x_4)$. Then
 $\Psi$ is  the pullback of $ \delta \cdot \Cal O_X $ under the map
$\psi   $ .

\endproclaim
\addline

Thus we have to show that
$\psi (\rho ) \in \delta {}^. \Cal O_X $.

We can do this calculation fiberwise. Over a fiber $W_1 \times
W_2$, we choose a symplectic basis of $V$ : $ \{ e_1, ... , e_4 \}
$, and its dual basis $\{ y_1, ... , y_4 \}$, in such a way that
$W_2 = e_1 \wedge e_2 $. Without loss of generality, we may assume
$f=y_3$ and $g=y_4$. Let $ v =\Sigma a_i e_i $ and  $ w =\Sigma
b_i e_i $. Then
$$
\matrix
 & f(v) = a_3, & g(v) = a_4, & <v,\ >=a_1 y_3 + a_2 y_4 - a_3 y_1
- a_4 y_2   ,\\
 & f(w) = b_3, & g(w) = b_4  , & <w,\ >=b_1 y_3 + b_2 y_4 - b_3
y_1 - b_4 y_2 .
\endmatrix 
$$
 And $<v,w>=0$ is equivalent to $a_1 b_3 + a_2 b_4 - a_3 b_1 - a_4
b_2 = 0 $,  thus $a_1 b_3  - a_3 b_1 = -( a_2 b_4 - a_4 b_2 ) $.
Finally:
$$
\gather
\psi ( \rho) = [ g(w) <v,\ > - g(v) <w,\ >] \wedge f - [ f(w) <v,\
> - f(v) <w,\ >] \wedge g \\ =[ b_4 <v,\ > - a_4 <w,\ >] \wedge
y_3 - [ b_3 <v,\ > - a_3 <w,\ >] \wedge y_4 \\ =[ b_4 [ a_1 y_3 +
a_2 y_4 - a_3 y_1 - a_4 y_2 ]  -  a_4 [ b_1 y_3 + b_2 y_4 - b_3
y_1 - b_4 y_2 ] ] \wedge y_3 \\ - [ b_3 [ a_1 y_3 + a_2 y_4 - a_3
y_1 - a_4 y_2 ]  -  a_3 [ b_1 y_3 + b_2 y_4 - b_3 y_1 - b_4 y_2 ]]
\wedge y_4 \\ = (b_4 a_2 - a_4 b_2 ) y_4 \wedge y_3 - (b_3 a_1 -
a_3 b_1 ) y_3 \wedge y_4 \\ +(-b_4 a_4 + a_4 b_4 ) y_2 \wedge y_3
- (-b_3 a_3 + a_3 b_3 ) y_1 \wedge y_4 \\ + (-b_4 a_3 + a_4 b_3)
y_1 \wedge y_3 - (-b_3 a_4 + a_3 b_4) y_2 \wedge y_4 \\ =(a_4 b_3
- b_4 a_3 ) (y_1 \wedge y_3 + y_2 \wedge y_4  ) \\ = (a_4 b_3 -
b_4 a_3 ) {}^. \delta 
\endgather
$$

Thus the map from $K$ to the sub-complex is well defined and
nonzero.  Next we show that $K$ is indeed isomorphic to the
kernel. 

\subhead{K $ \cong  $ kernel}
\endsubhead

First we show that we can compute the kernel from the quotient
complex:

\proclaim{Lemma 5.4} Let $D.$ be an exact complex, $B. \subset D.$
a sub-complex, and $Q.=D./B.$ the quotient complex. Suppose $ Q_j
=D_j$ for $j \geq m+3, $ and $Q.$ exact up to degree $ m+1 $. Then
$$
\text{kernel}\ ( B_{m+1} \to B_{m} ) = H_{m+2}( Q.) .
$$
\endproclaim

\demo{Proof} This can be proven by a diagram chase. Consider the
following diagram:
$$
\alignat 9 & &  & &  & &  & &  & \quad \quad \quad \quad \quad
\quad \quad \quad \quad \quad
\\
  &  & 0 \ & \to & K & \to & B_{m+1} & \to & B_{m} &\to ... \\ &
&  & &  & & \cap \quad & & \cap &  \\ & ...
\to & D_{m+3} \ & \to & D_{m+2} & \to & D_{m+1} & \to & D_{m} &\to
... 
\\ & & \| \quad & &
\| \quad & &
\downarrow \quad & & \downarrow &  
\\ & ... \to & Q_{m+3} \ & \to & Q_{m+2} & \to & Q_{m+1} & \to &
Q_{m} &\to ... \\ & & \uparrow \quad & & \uparrow \quad & &
\uparrow \quad & &  &  
\\ & & Q.\text{ exact} \ & & Q.\text{not exact} & & Q.\text{
exact} & &  &  
\endalignat
$$ Let $K=$ kernel $(B_{m+1}  \to  B_{m} )$ and $H=H_{m+2} ( Q.
)$. Then
$$
\gather H = ker (Q_{m+2}  \to  Q_{m+1} ) / Im (Q_{m+3}  \to 
Q_{m+2} ) \\
\cong  ker (D_{m+2}  \to  Q_{m+1} ) / ker (D_{m+2}  \to  D_{m+1} ).
\endgather
$$ For $ b\in K$, there exist a $ c \in D_{m+2}$, unique up to
$ker (D_{m+2}  \to  D_{m+1} ) $. This defines a map $f : K \to H
$.  On the other hand if $c \in ker (D_{m+2}  \to  Q_{m+1} )$,
then $dc \in K$, and if 
$ c \in  ker (D_{m+2}  \to  D_{m+1} ) $, then $dc=0$. This defines
a map 
$g: H \to K $. And since $ g \circ f = id_K $, it follows that $H
\cong K. $ 

q.e.d.
\enddemo

This implies for our computation:
$$
\gather ker(\overset 2 \to \Lambda S \boxtimes \Psi @> d >> S
\boxtimes S^\perp) \cong H_3(Q. )=  ker(Q_3 \to Q_2 )/  im(Q_4 \to
Q_3 )\ . 
\endgather
 $$ Thus all we need to show is that the rank of
$H_3(Q. )$ is equal to the rank of $ K $ which is equal to 1.
Since the quotient complex resolves the vector bundle 
$Sym_2 S \otimes \overset 2 \to \Lambda S^\perp $, the rank of  $
H_3(Q. ) $ is a constant. Therefore it is enough to compute the
rank of the kernel over a point $W \times W$ on the diagonal 
inside of $X \times X$:

Consider once more the quotient complex over a point $W \times W$:
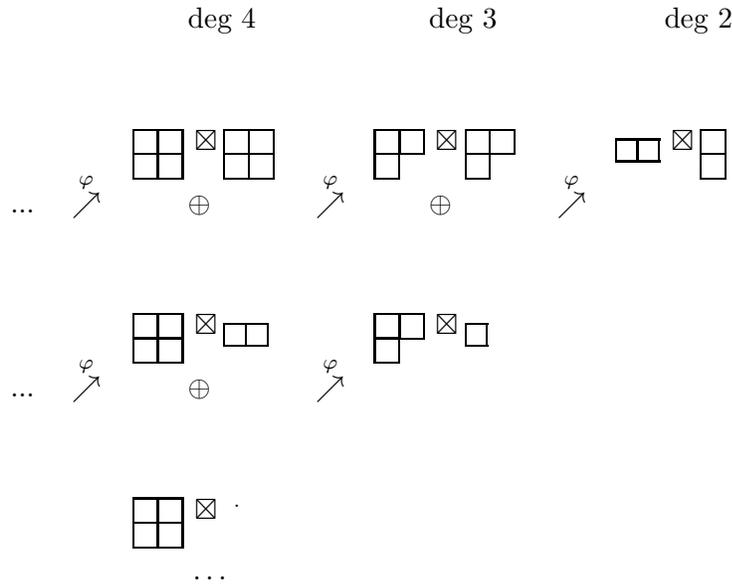
$$
  \alignat 8
  &  \text {deg 4} & &  &  & \text {deg 3}
 & & &  & \text {deg 2}
   &  & &   \\
\\ 
 \yng{2,2} & \boxtimes \yng{2,2} & &  & \yng{2,1} & \boxtimes 
\yng{2,1} & & & 
  \young{2} & \boxtimes \yng{1,1} & &  & \quad \quad \quad \quad
\quad \quad \quad
\\ ...\quad \overset \varphi \to \nearrow \quad \quad \quad  &
\oplus & & \quad \overset \varphi \to
\nearrow
\quad &   &
\oplus && \quad \overset \varphi \to \nearrow \quad &&  & & \\
\\  \yng{2,2} & \boxtimes \young{2}  & &  & \yng{2,1} & \boxtimes 
\young{1} & &  
\\ ... \quad \overset \varphi \to \nearrow \quad \quad \quad &
\oplus && \quad \overset \varphi \to 
\nearrow
\quad && && && \\
\\  \yng{2,2} & \boxtimes \ ^. &{}& &{}& & {} &
& 
\\ & \dots                                     
\endalignat
$$

$$
\gather
\text{Since} \ 
\yng{2,2}  \boxtimes \ ^. \to \yng{2,1}  \boxtimes \young{1} \to
\young{2}  \boxtimes \yng{1,1}
\to 0 \ \text{is exact, see Lemma 5.3, }
\endgather
$$ the rank of the kernel is equal to the rank of
$$
\gather ker (\yng{2,1} \boxtimes \yng{2,1} \to 0 ) /  im
(\yng{2,2} \boxtimes \young{2} \to \yng{2,1}
\boxtimes
\yng{2,1})
\endgather
$$

We again use Corollary 4.5: Set $\Cal V= \Cal O_{\Bbb P} (1),
x=t=t_{(g,g)} $ as in the section on the sub-complex, then 
$$
\gather   d_x: \   \yng{1,1} \boxtimes \young{1}  \cong \yng{2,1}
\boxtimes \yng{1,1}\ \text {is an isomorphism. }
\endgather  
$$ Tensor this on the right by $\young{1}$, then
$$
\alignat 8  
  &\yng{1,1} \boxtimes \young{2} 
\\
 d_x ( 1) : &\hskip .2in \oplus  & \ \cong \ & \yng{2,1} \boxtimes
\yng{2,1}
\\ &\yng{1,1} \boxtimes \yng{1,1}  & & & &\ \hskip .4in \text {is
an isomorphism. } \hskip 1in
\endalignat  
$$

 Recall the definition of  $\varphi$:
$$
\gather
 \yng{2,2} \boxtimes \young{2} @> \varphi >> \yng{2,1} \boxtimes
\yng{2,1} \\ \\ 
\varphi ((s_1 \wedge s_2) \otimes (t_1 \wedge t_2) \boxtimes f
\cdot g  )  =(s_1 \wedge s_2) \otimes  (t_1 \wedge t_2) \cdot t
\otimes f \cdot g \\ =  d_x ( (t_1 \wedge t_2) \boxtimes f \cdot g
)
\endgather 
$$ Hence the difference between $ \varphi $ and $ d_x $ is just a
twist by 
$\Cal O_{\overset \wedge \to Y} (-1)$, which is trivial over the
affine variety $\overset \wedge
\to Y$.

Therefore
$$
\alignat 8  
  &\yng{2,2} \boxtimes \young{2} 
\\
 &\hskip .4in \oplus  & \ \cong \ & \yng{2,1} \boxtimes \yng{2,1}
\\ &\yng{2,2} \boxtimes \yng{1,1} \  & & & &\ \hskip .4in  \hskip
1in
\endalignat  
$$ and thus
$$
\gather H_3(Q.)|_{W \times W} \\ \\
\cong
\yng{2,1} \boxtimes \yng{2,1}\  / \  im ( \yng{2,2} \boxtimes
\young{2} \to \yng{2,1} \boxtimes
\yng{2,1} )
\\ \\
 \cong ( \yng{2,2} \boxtimes \young{2}  \oplus \yng{2,2} \boxtimes
\yng{1,1} )
\ / \ \yng{2,2} \boxtimes \young{2}  \\ \\
\cong \yng{2,2} \boxtimes \yng{1,1} \ .
\endgather 
$$ Note that this has rank one and  and it also is clearly
isomorphic to
 $K|_{W \times W}$.
 This completes the proof of Theorem 5.2.

 q.e.d.

Next we discuss the K-theory of $X$.

\specialhead K-theory of $X$
\endspecialhead

Let $\Cal P$ be the category of vector bundles over $\Bbbk$, and
$\Cal C_X$ be the category of locally free coherent sheaves on
$X$. Define the subcategories $\Cal C(\Cal X^*), \Cal C(\Cal Y)$
as follows:
$$
\align
\Cal C(\Cal X^*) = \{ \Cal F \in \Cal C_X | H^i(X, \Cal F \otimes
\Cal G^*) = 0 \
\text {for all} \ \Cal G \in \Cal X \}  \hskip 1 in\\
\Cal C(\Cal Y) = \{ \Cal F \in \Cal C_X | H^i(X, \Cal F \otimes
\Cal G) = 0 \
\text {for all} \ \Cal G \in \Cal Y \}  \hskip 1 in
\endalign
$$  We will use the same standard trick, that we used in the case
of ordinary grassmannians, to construct resolutions for all sheaves
$\Cal F \in \Cal C(\Cal Y)$. Denote by $p_1$ and $p_2$ the
projections to the first and second factor:
$$  p_i: X \times X \to X,
$$  and denote by $q: X \to Spec (\Bbbk)$ the structure map. Then
$$  R^i{p_2}_*(p_1^* \Cal F \otimes \Cal O_\Delta) \cong
\left \{
\matrix
\Cal F \ &\text{for} \ i=0 \\ 0 &\text{otherwise}
\endmatrix
\right. 
$$

\proclaim{Lemma 5.5}

For $\Cal F \in \Cal C(\Cal Y) $,  $R^i p_{1*} ( p_2^* \Cal F
\otimes B_. ) $ gives rise to a finite resolution of $\Cal F$:
$$
\gather 0 \to H^0 ( X,\Cal F \otimes ( \overset 2 \to \Lambda S
)^{\otimes 2}) \otimes \overset 2 
\to  \Lambda S^\perp 
\to H^0 ( X,\Cal F \otimes \overset 2 \to \Lambda S ) \otimes \Psi
\\
\to H^0 ( X,\Cal F \otimes S ) \otimes S^\perp 
\to  H^0 ( X,\Cal F ) \otimes O_X \to \Cal F \to 0
\endgather
$$
\endproclaim

\demo{Proof}

Recall that  $B_.$ is the sub-complex resolving $\Cal O_{\Delta} $.
For $\Cal F \in \Cal C(\Cal Y) $, the higher direct images of 
$ ( p_2^* \Cal F \otimes B_. ) $ all vanish, since 
$$ R^i p_{2*} ( p_1^* \Cal F \otimes ( \Sigma^\alpha S \boxtimes
\Psi_{\alpha^*} )) \cong H^i(\Cal F
\otimes \Sigma^\alpha S ) \otimes
\Psi_{\alpha^*} .
$$ Since all the higher direct images vanish, it follows that the
push-forward
$$ p_{2*} ( p_1^* \Cal F \otimes B_.) 
$$ is exact, see also Lemma 1.3.

 q.e.d.

\enddemo

Similar to the case of  ordinary grassmannians, define 
$$
 U_\Cal G : \Cal P \to \Cal C_X , \ U_\Cal G (W) = q^* W \otimes
\Cal G
$$ for all $ \Cal G \in \Cal X$ .  Then the $U_\Cal G$'s are
inducing
 homomorphisms $u_\Cal G$ on the K-theory, i.e.
$$ u_\Cal G :  K_*( \Bbbk ) \to K_*(X) 
$$  Set $u= \underset {\Cal G \in \Cal X } \to \bigoplus \  u_\Cal
G $, then

\proclaim{Theorem 5.6}
$$ u: \underset \Cal G \in \Cal X \to \bigoplus K_*( \Bbbk ) \to
K_*(X) 
\qquad \text {is an isomorphism.}
$$ 
\endproclaim

\demo {Proof} 

  First we show that $u $ is surjective:

Recall that it is enough to define the K-theory on the
subcategory $\Cal C (\Cal Y)$ or $\Cal C (\Cal X^*)$,
Lemma 1.6.
 Consider the finite resolution of $\Cal F \in \Cal C(\Cal Y)$,
then
$$
\gather
 \lbrack {\Cal F } \rbrack  = 
\lbrack \oh^0(\Cal F  ) \otimes \Cal O_X
  \rbrack -
\lbrack \oh^0(\Cal F \otimes S^\perp ) \otimes S
  \rbrack \\ +
\lbrack \oh^0(\Cal F \otimes \Psi ) \otimes \overset 2 \to \Lambda
S
  \rbrack -
\lbrack \oh^0(\Cal F \otimes \overset 2 \to \Lambda S^\perp )
\otimes {(\overset 2 \to \Lambda S)}^ {\otimes 2}
  \rbrack
\endgather
$$  which implies  that $u$ is surjective. 

To show that $u $ is injective we  define a homomorphism
$$ v: K_*(X) \to \underset \Cal G \in \Cal X \to \bigoplus K_*(
\Bbbk ) \ .
$$
 and then show   that $v \circ u$   is injective, consequently
that $u$ is an  isomorphism.

It
suffices to define $v$ on the subcategory
$\Cal C(\Cal X^*) $ rather than defining it on the whole category 
$\Cal C_X $, Lemma 1.7.

Define
$$
\gather
 V_\Cal G : \Cal C(\Cal X^*) \to \Cal P , \ \ V_\Cal G(\Cal F ) =
q_*(\Cal F \otimes
\Cal G^* ) \ \text{for all}\ \Cal G \in \Cal X \  .
\endgather
$$  Note that  $V_\Cal G$ on the  category
$\Cal C_{X} $ is not exact, because $R^i q_* (\Cal F \otimes \Cal
G^* ) $ are in general not zero  for $ i>0$, but the $V_\Cal G$'s
are exact on the subcategory $\Cal C(\Cal X^*)$ . The $V_\Cal G$'s 
induce homomorphisms
$v_\Cal G$  on the K-theory.
$$
 v_\Cal G :  K_*(X) \to  K_*( \Bbbk ).
$$
 
Let us give an ordering to $\Cal X $:
$$
\gather
\Cal X:\ \Cal O_X (-2) \ \leq \Cal O_X (-1) \ \leq S \ \leq \Cal
O_X \ .
\endgather
$$

Let $ v=(...,v_\Cal G,...)$  ,  then consider  $v
\circ u $:

$$ V_\Cal H \circ U_\Cal G (W) =V_\Cal H ( q^* W \otimes \Cal G )
= \lbrack
\oh^0 ( W \otimes \Cal G \otimes ( \Cal H )^* ) \rbrack
$$

 We will show that the matrix of $v \circ u $ with respect  to the
ordering of $\Cal X$ is upper triangular with  ones down the
diagonal, thereby showing that $u$ is an isomorphism. Thus we need
to show that $\Cal X$ is orthogonal, that is: 
\enddemo

\proclaim{Lemma 5.7}

For all $\Cal G , \Cal H \in \Cal X $:

$$
 Ext^i(\Cal H, \Cal G) = 
\cases
\Bbbk & \text {for } \Cal G = \Cal H \ \text{and} \ i=0 \\ 0 &
\text {for } \Cal G < \Cal H  \ \text{and} \ i=0 \\ 0 & \text {for
all} \ \Cal G , \Cal H \in \Cal X \ \text{and}\ i>0 
\endcases
$$

\endproclaim 

\demo{Proof} We need some cohomology computations.

\enddemo

\specialhead {Cohomology Computations}
\endspecialhead

Here we  use the fact that the symplectic grassmannian
$X=SpGr(2,4) $ sits inside the grassmannian $Gr(2,4) $. Therefore
we can apply the cohomology computations for the grassmannian
$Gr(2,4) $ to $X$:

Let
$\ F=\left\{ 0 \subset W_1 \subset W_2 \subset ... \subset W_N = V
\text { with dimension}\ W_i = i  \right \} 
$  be the full flag manifold and $\pi$   the projection  $ \pi: F
\to Gr(2,V)$.  Consider the short exact sequence:
$$
 0 \to \Cal O_{Gr} ( -X) \to \Cal O_{Gr} \to \Cal O_X \to 0 \ .
$$  Recall that 
$\Cal O_{Gr} (-X) \cong \Cal O _{\Bbb P ({ \overset \ 2 \to
\Lambda } V^* )} (-1) \cong \overset \ 2
\to \Lambda S = \pi_* \Cal O_F ( -1,-1,0,...,0).$

Tensor the above short exact sequence by a vector bundle $\Cal V$
on  the grassmannian $Gr(2,4)$ and consider the  resulting long
cohomology sequence:
$$
\aligned
  0 \to & \oh^0 (Gr,\Cal V  \otimes \Cal O_{Gr} (-X)) \to
\oh^0 (Gr,\Cal V) \to \oh^0 (X,\Cal V \otimes \Cal O_X) \\
 \to & \oh^1 (Gr,\Cal V \otimes \Cal O_{Gr} (-X)) \to
\oh^1 (Gr,\Cal V ) \to \oh^1 (X,\Cal V \otimes \Cal O_X)
 \to ...
\endaligned \tag * 
$$   In particular we see from here, that if $\Cal V$ has no
cohomology, then 
$$ 
\oh^{i} (X,\Cal V \otimes \Cal O_X) \cong \oh^{i+1} (Gr,\Cal V(-X)
).
$$ Recall that $\Sigma^\alpha S_{Gr} \cong \pi_* \Cal O_F (
-\alpha_2, -\alpha_1, 0, 0 ) $ and also recall if $\gamma   $ is
not a dominant weight and if
$\gamma + \rho$ with $\rho = (4,3,2,1) $, has a repetition, then 
$\Cal O_F (\alpha) $ has no cohomology.

Although we have done most of the cohomology computations, that are
needed here,  in  chapter 1, we will do  these computations again
in detail:

\proclaim{Lemma 5.8}
$$ H^i(X, \Cal O_X ( k ) ) = 
 \left\{ \matrix
\Bbbk & k=0,&i=0 \\ 0 &  k = -1, -2, & i=0
\\ 0 & k \geq -2,& i>0
\endmatrix
\right.
$$
\endproclaim

\demo{Proof} Recall that
$\Cal O_{Gr}(k) \cong \pi_* ( \Cal O_F (k,k,0,0) )$. The
cohomology of these twisted sheaves is given by:

$$  H^i(Gr, \Cal O_{Gr} ( k ) = 
 \left\{ \matrix
\Bbbk & k=0,&i=0 \\ 0 &  k = -1, -2, -3, & i=0
\\ 0 & k \geq -3,& i>0
\endmatrix
\right.
$$

This follows from:

$(-j,-j,0,0) + \rho = (4-j, 3-j, 2, 0 ) $ has a repetition for
$j=1,2,3 $, thus
$\Cal O_{Gr}(-j) $ has no cohomology for $j=1,2,3 $.  and for
$k>0: (k,k,0,0)  $ is dominant, thus $\Cal O_{Gr}(k) $ has only
$H^0$.

The statement follows now from the long exact cohomology sequence
(*). 

q.e.d.
\enddemo

Let us discuss the terms in Lemma 5.7 involving $S$:

\proclaim{Lemma 5.9}
$$
\gather H^i ( X, S (k) ) = 0 \ \text {for} \ k=0,1,2  \ \text{and}
\ i > 0 , \tag a \\ H^i ( X, S^* (- k) ) = 0 \ \text {for} \
k=0,1,2  \ \text{and}\ i > 0 ,\tag b \\ H^i ( X, S\otimes S^*  ) =
0 \ \text {for} \  i > 0 ,\tag c \\ H^0 ( X, S\otimes S^*  ) =
\Bbbk \tag d \\ H^0 ( X, S  ) = 0  \tag e
\endgather
$$

\endproclaim

\demo{Proof}
$$
\gather  S_{Gr}(-1)  \cong \pi_* \Cal O_F (-1,-2,0,0) \tag {a),(e}
\\
 S_{Gr} \cong  \pi_* \Cal O_F (0,-1,0,0)  \\  S_{Gr}(1) \cong 
\pi_* \Cal O_F (1,0,0,0) \\ S_{Gr}(2) \cong  \pi_* \Cal O_F
(2,1,0,0)
\endgather
$$  After adding $\rho=(4,3,2,1)$ to the first two weights, these
have repetitions and thus $S_{Gr}(-1) $ and $  S_{Gr}$ have no
cohomology.  The last two weights are dominant and thus $S_{Gr}(1)
$ and $S_{Gr}(2) $ have only $H^0$. Using the long exact
cohomology sequence (*), we get that $S$ has no cohomology and 
$S(1)$ and $S(2)$ have no higher cohomology.

$$
\gather   S_{Gr}^*(-3) \cong  \pi_* \Cal O_F (-2,-3,0,0) \tag b \\
 S_{Gr}^*(-2) \cong  \pi_* \Cal O_F (-1,-2,0,0)  \\
  S_{Gr}^*(-1) \cong  \pi_* \Cal O_F (0,-1,0,0) \\ S_{Gr}^* \cong 
\pi_* \Cal O_F (1,0,0,0)
\endgather
$$ The first three weights have repetitions after 
 adding $\rho$ to them, thus they have no cohomology.  The last 
weight is dominant and thus has only $H^0$. The claim  again
follows from the long exact cohomology sequence (*).

$$
\gather   S_{Gr} \otimes S_{Gr}^* (-1) 
\cong S_{Gr} \otimes S_{Gr}  
\cong \overset 2 \to
\Lambda S_{Gr} \oplus  Sym_2 S_{Gr} \tag {c),(d} \\
\cong \pi_* \Cal O_F(-1,-1,0,0)  \oplus \pi_* \Cal O_F(0,-2,0,0).
\endgather
$$  Since $(0,-2,0,0) + \rho = (4,1,2,1 )  $ has a repetition, 
$S_{Gr} \otimes S_{Gr}^* (-1)$ has no cohomology.

$$
 S_{Gr} \otimes S_{Gr}^*   \cong \Cal O_X \oplus  Sym_2 S_{Gr}(1)
\cong \pi_* \Cal O_F  \oplus \pi_* \Cal O_F(1,-1,0,0).
$$  Since $(1,-1,0,0) + \rho = (5,2,2,1 )  $  has a repetition,  
$ S_{Gr} \otimes S_{Gr}^*$ has only $H^0 \cong \Bbbk \cong H^0
(\pi_* \Cal O_F) $.

Using the long exact cohomology sequence (*), it follows that
 $S \otimes S^* $ over $X$ has only $H^0 \cong \Bbbk$.

Hence $\Cal X$ is orthogonal. 

q.e.d.

\enddemo

\subhead {Proof of Lemma 5.7}
\endsubhead Recall that $Ext^i(\Cal F, \Cal G) \cong H^i(X, \Cal
F^* \otimes \Cal G)$. Lemma 5.7 then follows from  Lemma 5.8 and
5.9. q.e.d.

\subhead {Proof of Theorem 5.6}
\endsubhead Consider the matrix of $v \circ u $ with respect  to
the ordering of $\Cal X$:
$$  V_\Cal H \circ U_\Cal G (W) =V_\Cal H ( q^* W \otimes \Cal G )
= \lbrack
\oh^0 ( W \otimes \Cal G \otimes ( \Cal H )^* ) \rbrack
$$ Lemma 5.7 implies that
\item{(a)} $V_\Cal H \circ U_\Cal G$ is well defined for 
$\Cal G, \Cal H \in \Cal X$, since $ \Cal X \subset C(\Cal X^*)$
and
\item{(b)} it implies that this matrix  with respect  to the
ordering of $\Cal X$ is upper triangular with  ones down the
diagonal. 

This proofs that $u$ is injective, thus that $u$ is an
isomorphism. 

q.e.d.

\newpage

\topmatter
\title\chapter{6} Proof of the Main Theorem \endtitle
\endtopmatter

In this chapter we will finish the proof of the Main Theorem 3.2. Let us first
recall some notations and the statement  of the main theorem.

\head {Notations} \endhead

Let $X=SpGr(2,V)$ be the symplectic grassmannian of 2-planes in $V$, and let $V$
be a vector space of dimension $N=2n$ over an algebraically closed field $\Bbbk$
of characteristic 0. We keep the notations of chapter 2 and 3. 

Let $D.$ be the complex of the Tate construction that was defined in chapter 2,
thus 

$$  D_i= \underset { \underset
   N-2 \geq \alpha_1 \geq \alpha_2 \geq 0 \to {| ( \alpha )| =i} 
   } 
\to \bigoplus \Sigma ^\alpha S 
\boxtimes \Psi_{\alpha^*} \ .
$$

and denote by $B.$ the sub-complex:
$$ B.: \ B_{2N-6} \to ... \to B_2 \to S \boxtimes S^\perp \to \Cal O_{X \times
X} ,
$$ where $B_i \subset D_i,  i \leq 2N-6$ is defined as the direct sum  of
$\Sigma ^\alpha S 
\boxtimes \Psi_{\alpha^*} $ over all $\alpha $ of length $i$
 with $\alpha_1 \neq N-2$.

\proclaim{Theorem 6.1}
 
The sub-complex $B.$ is exact.

\endproclaim

\proclaim {Remark}

We  have proven this theorem for $X=SpGr(2,4).$ In this case the sub-complex is
given by 
$$
\overset 2 \to \Lambda S \boxtimes \Psi @> d >> S \boxtimes S^\perp \to \Cal
O_{X \times X}. 
$$ 

Moreover for $X=SpGr(2,4) $ we not only proved the exactness of this complex,
but we also found the kernel of the  map 
$d :\overset 2 \to \Lambda S \boxtimes \Psi @> d >> S \boxtimes S^\perp$, which
is equal to $\Cal O_X(-2) \boxtimes \Cal O_X(-1) $. This enabled us to find the
generating system
$$
 \Cal X= \{ \Cal O_X (-2),\Cal O_X (-1), S, \Cal O_X \}
$$ for the K-Theory of $X$. In the general case we will only show the exactness
of the sub-complex.

\endproclaim

\specialhead {Proof of the Main Theorem:}
\endspecialhead

Recall that at the end of chapter 3, we reduced the proof to showing  that the
quotient complex
$Q. = D. / B. $ is exact up to degree $2N-6$ over points on the diagonal $\Delta
\subset X \times X $. Over a point $W \times W $ on the diagonal, the quotient
complex splits into many  different sequences, all of which come in one way or
another from the sequences of chapter 4.

Let us use the same notations that were used in chapter 2, 3 and 5. Thus the
external  tensor product of  two Young diagrams  of weight $\alpha $
respectively $\beta $ represents the  external tensor product of the
corresponding representations of
$W$ and $(V/W)^*$, that is
$$
\alpha \boxtimes \beta =\Sigma^\alpha W \boxtimes \Sigma^\beta (V/W)^* .
$$

We will study the quotient complex over a point
$W \times W $ on the diagonal $\Delta $.

\head {The maps $ev$ and $\varphi$}
\endhead

Over a point $W \times W $,   $ \Psi_\alpha $  splits into a direct sum 
$$
\Psi_\alpha|_{W \times W} \cong \underset i \geq 0 \to \bigoplus
\Sigma^\alpha W \ \boxtimes \  \Sigma^{\alpha^*-2i} (V/W)^* \ .
$$
 The map $d$ of the  quotient complex, analogously to the
case of $X=SpGr(2,4) $, splits into two kinds of maps:
$ev$ and $\varphi $.

Let $gP$ be the parabolic leaving $W$  invariant. Recall the definition of $
T^i=T^i_{(g,g)} $ and $t= t_{(g,g)} $.
For $g \in G , i > 0 $, 
$$
\gather
  T = (g, v_1 \wedge v_2 ) \boxtimes (g,\eta) \ ,
 T^i = (g, (v_1 \wedge v_2)^{\otimes i } ) 
\boxtimes (g,\eta^i)  \\
 \text {and} \ t= d T \ .
\endgather
$$

$$
\gather
\text{Let} 
 \ w \in \Sigma^{(\alpha_1-i,\alpha_2-i)} W \boxtimes 
\Sigma^{(\alpha^*-2i)} ( V/W)^* \ , \text {then}
\\ d ( w T^i)=  ( dw) T^i + (-1)^{ | w| }i w \cdot t \cdot T^{i-1} \ .
\endgather
$$ Then the maps $ev $ and $\varphi $ are given by:
$$
\alignat 8  & & & \Sigma^{\alpha-1} W \ \boxtimes \   \Sigma^{\beta-1} (V/W)^* & 
\hskip 1 in
\\ & & \overset ev \to \nearrow \\ & d: \Sigma^\alpha W \ \boxtimes \  
\Sigma^\beta (V/W)^* & & \hskip .4 in\oplus
\\ & & \overset \varphi \to \searrow \\ & & & \Sigma^{\alpha-1} W \ \boxtimes
\   
\Sigma^{\beta+1 } (V/W)^*
\endalignat
$$
$$
\gather  ev ( (w_1 \wedge w_2)^{\otimes i} \otimes w )  =  (w_1 \wedge
w_2)^{\otimes i} \otimes d w
\endgather
$$ and 
$$
\varphi ((w_1 \wedge w_2)^{\otimes i} \otimes w) =(-1)^{ | w| }i 
(w_1 \wedge w_2)^{\otimes i} \otimes ( w \cdot t) \
.
$$

Over a point $ W \times W$ on the diagonal $\Delta$, all the evaluation maps
$ev$ are zero, because for $v \in W$ and $f \in (V/W)^*$,
$f(v)=0$, thus $dw=0$. Therefore the quotient complex splits into several
subsequences.

\head {Quotient Complex}
\endhead

Since all the evaluation maps are zero over points on the diagonal, the quotient
complex splits into $N-2$ sequences 

$$ A^0, A^1, \dots A^{N-2} \text { as follows:}
$$

$$ Q|_\Delta \cong A^0 \oplus \dots \oplus A^{N-2} .
$$

$$
\alignat 9 &\text{deg }2N-5 
\\
\\ A^{N-3} : 
\ \ &{\left. \aligned  \undersetbrace  {N-2} \to {\underset \to {\yng{6,5}}}
\boxtimes 
 \yng{2,2,2,2,2,1} \endaligned \right \} {\ssize { N-2 }} }\ & \to & \ 0 &
\hskip 2.0in
\\
  &\text{deg }2N-5 & & \text{deg }2N-6 & & 
\\
\endalignat
$$

$$
\alignat 9
\\ A^{N-4} : \ \  & {\left. \aligned  \undersetbrace  {N-2} \to {\underset \to
{\yng{6,5}}} 
\boxtimes \yng{2,2,2,2,1} \endaligned \right \} {\ssize { N-3 }} }\ & \ \to \ &
\yng{6,4} \boxtimes \yng{2,2,2,2,1,1}\ & \to & \ 0 & \hskip 1in
\\ & \hskip .7in \oplus & \nearrow
\\ & {\left. \aligned  \undersetbrace  {N-1} \to {\underset \to {\yng{7,4}}} 
\boxtimes \yng{2,2,2,1,1,1}  \endaligned \right \} {\ssize { N-2 }} }\  
\\
\endalignat
$$

\newpage

$$
\alignat 9 &A^{N-5}: \\ &\hskip  .4 in \text{deg }2N-5  \hskip  1 in  \text{deg
}2N-6 \hskip  .8 in   \text{deg }2N-7 
\\
\\ &{\left. \aligned  \undersetbrace  {N-2} \to {\underset \to {\yng{6,5}}} 
\boxtimes \yng{2,2,2,1} \endaligned \right \} {\ssize { N-4}} }
 \to \
\yng{6,4}
\boxtimes \yng{2,2,2,1,1}  \ \ 
\to \ \ 
\yng{6,3}
\boxtimes \yng{2,2,2,1,1,1} 
 \to  0
\\ & \hskip .7in \oplus \hskip  .8 in  \nearrow  \hskip .8in \oplus \hskip .6
in  \nearrow 
\\ &{\left. \aligned  \undersetbrace  {N-1} \to {\underset \to {\yng{7,4}}} 
\boxtimes \yng{2,2,1,1,1} \endaligned \right \} {\ssize { N-3 }} }\ 
\to  \yng{7,3} \boxtimes \yng{2,2,1,1,1,1}  
\\ & \hskip .7in \oplus  \hskip  .8 in  \nearrow 
\\ &{\left. \aligned  \undersetbrace  {N} \to {\underset \to {\yng{8,3}}} 
\boxtimes \yng{2,1,1,1,1,1} \endaligned \right \} {\ssize { N-2 }} }\ 
\\
\endalignat
$$

$\dots$

$\dots$

$$
\alignat 8 &\text{deg }2N-5 & & \text{deg }2N-6 & & \text{deg }N-1 
\\
\\ A^1: \ &{\left. \aligned  \undersetbrace  {N-2} \to {\underset \to
{\yng{6,5}}} 
\boxtimes \yng{2,1}\endaligned \right. } &  \to  &{\left. \aligned 
\undersetbrace  {N-2} \to {\underset \to {\yng{6,4}}} 
\boxtimes \yng{2,1,1}\endaligned \right. }
 &\to \dots \to &{\left. \aligned  \undersetbrace  {N-2} \to {\underset \to
{\yng{6,1}}}
 \boxtimes \yng{2,1,1,1,1,1} \endaligned \right. } \to 0
\\ & \hskip .7in \oplus & \nearrow & \hskip .7in \oplus & \nearrow \dots \nearrow
\\ &\yng{7,4} \boxtimes \yng{1,1,1} & \to  &\yng{7,3} \boxtimes \yng{1,1,1,1} &
\ \to \dots \to 0 & &
\endalignat 
$$

\newpage

$$
\alignat 8 &\text{deg }2N-5 &\text{deg }2N-6 \hskip .2 in & & 
\\
\\ A^0: \ & {\left. \aligned  \undersetbrace  {N-2} \to {\underset \to
{\yng{6,5}}} 
\boxtimes \young{1}\endaligned \right. }
 \to  & {\left. \aligned  \undersetbrace  {N-2} \to {\underset \to {\yng{6,4}}} 
\boxtimes \yng{1,1}\endaligned \right. }
\\
\\ &&\text{deg }N-1 \hskip .2 in & \hskip .2 in\text{deg }N-2 &
\\
\\ & \hskip .4 in \to
 \dots \to &\yng{6,1} \boxtimes \yng{1,1,1,1,1}
 \to  &{\left. \aligned
\yng{6} \boxtimes \yng{1,1,1,1,1,1} \endaligned \right \}  } {\ssize {N-2}} \to
0   & \hskip .4 in 
\endalignat
$$
\addline
 
All of these sequences $A^0, ... , A^{N-3} $ are exact up to degree 
$2N-6$. This is obvious for $ A^{N-3}$, because $A^{N-3}$ is zero up to degree 
$2N-6$. For the other sequences  the proof is  similar  to the case of
$X=SpGr(2,4) $. We will use the sequences of Chapter 4,  Corollary 4.5, to show
that
$A^0, ... , A^{N-3} $ are exact.

First recall the construction of these sequences:

\head { Some exact sequences }
\endhead

Similar to the previous chapter, let 
$Y= Hom(T,T')$, 
$Y_r=\{ x \in Y \ |\ \text{rank} (x) \leq r \} $ and
$\overset \wedge \to Y = Y - Y_r $. Set $V=Gr(n-r,T) $,
 $T=W^*$ and $T'=(V/W)^*$.

Then 
$Y=Hom(W^*,(V/W)^* ). $ Again consider the  case when  $V$ is equal to the
projective space 
 $ \Bbb P (W^*) \cong \Bbb P^1  $.

Recall the statement of Corollary 4.5 for $\Cal V = \Cal O_\Bbb P (j)$:

$$
\gather
 \overset m-j \to \Lambda {T'} \to T^* \otimes \overset m-j+1 \to \Lambda {T'}
\otimes 
\Cal O_Y
\to Sym_2(T^*)  \otimes \overset m-j+2 \to \Lambda {T'} \otimes \Cal O_Y \to
...\\ \to Sym_{j-1}(T^*)  \otimes \overset m-1 \to \Lambda{T'} \otimes \Cal O_Y 
\to  Sym_j(T^*)  \otimes \overset m \to \Lambda T' \otimes
\Cal O_Y \\
\text {is exact.}
\endgather 
$$

Over a point $x = s \otimes f \in Y = Hom(T,T') \cong T^* \otimes T', $  the
maps are given by:
$$
\gather d_x : Sym_k (T^*) \otimes \overset m-j+k \to \Lambda {T'}^* \otimes
\Bbbk (\varphi)
\to
 Sym_{k+1 }(T^*) \otimes \overset m-j+{k +1} \to \Lambda {T'}^* \otimes \Bbbk
(\varphi) \\
 d_x ( ( s_1 ... s_k ) \otimes ( f_1 \wedge ... \wedge f_{m-j+k} ) )  = s_1 ...
s_k s \otimes  f \wedge f_1 \wedge ...  \wedge f_{m-j+k}
\endgather
$$

We will use $x = t$. Let us check that $t \in \overset \wedge \to Y: $
 
Let $e_1, ... ,e_4$ be a symplectic basis of $V$, and $y_1, ... ,y_4$ the
corresponding dual basis of $V^*$, such that $W=e_1 \wedge e_2$. 

Recall that $\eta$ is defined as
$\eta =x_1 \otimes x_3 + x_2
\otimes x_4$ and denote by $d$ the map from the Tate construction. Then
$$
\gather t= dT \\ = d ((g, v_1 \wedge v_2) \boxtimes (g, \eta) ) \\ =  d ( (e_1
\wedge e_2) \boxtimes ( y_1 \otimes y_{n+1} + \dots + y_{n+2}
\otimes y_{2n} ) \\  = y_1 (e_1)e_2 \boxtimes y_{n+1} - y_1 (e_2)e_1
\boxtimes y_{n+1} \\ +   y_2(e_1) e_2 \boxtimes y_{n+2} - y_2(e_2) e_1 \boxtimes
y_{n+2} \\ + \dots + \\ y_n(e_1) e_2 \boxtimes y_{2n} - y_n(e_2) e_1 \boxtimes
y_{2n} \\ =e_2 \boxtimes y_{n+1} - e_1 \boxtimes y_{2n} \ .
\endgather
$$ Thus $t $ has rank 2 and therefore $t \in \overset \wedge \to Y$.

\proclaim{Remark} Note that
$$
\varphi ( (e_1 \wedge e_2 ) \otimes w ) = d_x ( w ) \ , 
$$  thus $\varphi$ and $d_x$ just differ by a  twist of 
$\Cal O_{\overset \wedge \to Y}(-1)$. Since $\overset \wedge \to Y $ is affine,
this twist is trivial.
\endproclaim

Let us now finish the proof of the theorem:

\head{The sequences $A^0, A^1 \dots A^{N-3}$ are exact}
\endhead

We will proof this in three steps:

First we show that $A^0$ is exact, then secondly we discuss the exactness of
$A^1$, and finally we consider $A^l$ for any $l$ and proof the exactness  by
induction.

The second step is not really necessary, but it is helpful to recognize a pattern
for the proof.

\head {Step 1}:
\endhead

\proclaim {Lemma 6.2}
$A^0$ is exact.
\endproclaim

\demo {Proof}

Set
$\Cal V = \Cal O_\Bbb P (N-2) \otimes \overset {N-2} \to \Lambda (V/W)^* $ in
4.5.  Then
$$
\gather
  {}^. \boxtimes {}^. @> d_x >> 
 \young{1}
\boxtimes \young{1}
 @> d_x >> 
\young{2}
\boxtimes \yng{1,1} @> d_x >> \dots
\\
\\
 \hskip .4 in 
  @> d_x >>
\yng{5} \boxtimes \yng{1,1,1,1,1}
 @> d_x >> {\left. \aligned
\yng{6} \boxtimes \yng{1,1,1,1,1,1} \endaligned \right \}  } {\ssize {N-2}} \to
0  
\endgather
$$

is exact. 

Recall that the difference between $d_x$ and $\varphi$ is that 
$\varphi$ adds a shift to the map $d_x$ (see above note), thus 

$$
\gather {{}^ {}}^0 {}^{ {}^{ @> \varphi >>} } {\left. \aligned  \undersetbrace  {N-2} \to
{\underset \to {\yng{6,6}}} 
\boxtimes {}^. \endaligned \right. }  {}^{ {}^{ @> \varphi >>} } {\left.
\aligned  \undersetbrace  {N-2} \to {\underset \to {\yng{6,5}}} 
\boxtimes \young{1}\endaligned \right. }  {}^{ {}^{ @> \varphi >>} } {\left.
\aligned  \undersetbrace  {N-2} \to {\underset \to {\yng{6,4}}} 
\boxtimes \yng{1,1}\endaligned \right. }
\\
\\
  @> \varphi >> \dots  @> \varphi >>
\yng{6} \boxtimes \yng{1,1,1,1,1,1} \to 0  \quad \text {is exact}.
\endgather
$$ q.e.d.

\enddemo

\head {Step 2}
\endhead

\proclaim {Lemma 6.3}

$A^1: $ is exact .
\endproclaim

\demo {Proof}

Consider $A^1$:
$$
\alignat 8 A^1: \ &{\left. \aligned  \undersetbrace  {N-2} \to {\underset \to
{\yng{6,5}}} 
\boxtimes \yng{2,1}\endaligned \right. } &   {}^{ {}^{ @> >> } }
 &{\left. \aligned  \undersetbrace  {N-2} \to {\underset \to {\yng{6,4}}} 
\boxtimes \yng{2,1,1}\endaligned \right. }
 &{}^{ {}^{  @> >> } } \dots {}^{ {}^{  @> >> } } &{\left. \aligned 
\undersetbrace  {N-2} \to {\underset \to {\yng{6,1}}}
 \boxtimes \yng{2,1,1,1,1,1} \endaligned \right. } \to 0
\\ & \hskip .7in \oplus & \nearrow & \hskip .7in \oplus & \nearrow \dots \nearrow
\\ &\yng{7,4} \boxtimes \yng{1,1,1} & \to  &\yng{7,3} \boxtimes \yng{1,1,1,1} &
\ \to \dots \to 0 & \quad  \text{is exact.} & &
\endalignat 
$$

The lower row of $A^1 $ is equal to the complex $A^0 $ shifted  by 2, i.e. $A^0
[-2]$, that is in degree $2N-5$  the element of the lower row is equal to $
A^0_{2N-7}$  and so on.
\enddemo

\subhead{ Recall Lemma 3.3}
\endsubhead

Let $\overset \wedge \to D. $ be an exact complex, let $\overset \wedge \to B. $
be a sub-complex and denote by $\overset \wedge \to Q. $ the quotient complex of
$\overset \wedge \to D. $ by   $\overset \wedge \to B. $.

Then  $\overset \wedge \to B. $ is exact up to degree $m$ if and only if 
$\overset \wedge \to Q. $ is exact up to degree $m+1 $.

Since $A^0[-2]$ is exact up to degree $2N-3$, this reduces the proof to showing
that the upper row of $A^1$ is exact.

Use $\Cal V = \Cal O_\Bbb P ((N-2)-1) \otimes \overset N-2 \to \Lambda (V/W)^* $.
Then by 4.5 the following complex is exact:

$$
\align {\left. \aligned  \undersetbrace  {N-2} \to {\underset \to {\yng{6,6}}} 
\boxtimes \young{1}\endaligned \right. }
 @> \varphi >>  {\left. \aligned  \undersetbrace  {N-2} \to {\underset \to
{\yng{6,5}}} 
\boxtimes \yng{1,1}\endaligned \right. } @> \varphi >> \dots
 @> \varphi >> {\left. \aligned  \undersetbrace  {N-2} \to {\underset \to
{\yng{6,1}}} 
\boxtimes
\yng{1,1,1,1,1,1} \endaligned \right \} {\ssize N-2 } } @> \varphi >> 0
\endalign
$$ Tensor this by $(V/W)^* $ on the right:
$$
\alignat 8 (*) \   &\yng{6,6} \boxtimes \young{2} & @> \varphi >>  &\yng{6,5}
\boxtimes \yng{2,1} & @> \varphi >> ... @> \varphi >> &\yng{6,1} \boxtimes
\yng{2,1,1,1,1,1} & @> \varphi >> 0
\\ & \hskip .5 in \oplus & \overset \varphi \to \nearrow & \hskip .5 in \oplus
\hskip .3 in \overset \varphi \to \nearrow
\\ &\yng{6,6} \boxtimes \yng{1,1} & @> \varphi >>  &\yng{6,5} \boxtimes
\yng{1,1,1}  @> \varphi >> & \dots \ \ \
\endalignat
$$

Then the top row (*) is what we would like to be exact.

The bottom row is the sequence of 
$\Cal V=\Cal O_\Bbb P ((N-2) -2 ) \otimes \overset N-2 \to \Lambda ( V/W)^* $,
which is exact up to degree $2N-5 $ by
 Corollary 4.5.

Using Lemma 3.3 again this implies the exactness of $(*)$, up to degree $2N-6 $,
as desired. q.e.d.

\addline

For the completion of the proof of the Main Theorem, we will  proof the
exactness in general for these sequences $A^l$. The theorem then follows.

\newpage

\head General Step:
\endhead

\proclaim {Lemma 6.4} The sequence $A^l $ is exact for all $ 0 \leq l \leq N-3$.
\endproclaim

\demo {Proof}

Consider $A^l$:

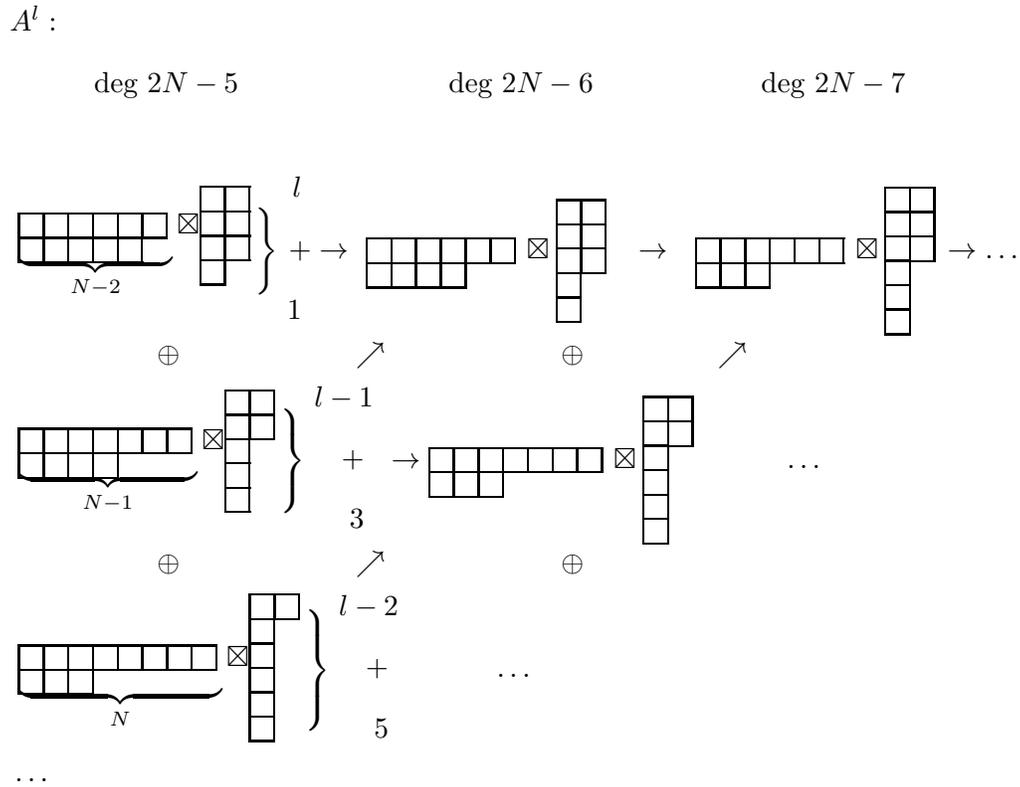
$$
\alignat 9 &A^{l}: \\ &\hskip  .4 in \text{deg }2N-5  \hskip  1 in  \text{deg
}2N-6 \hskip  .8 in   \text{deg }2N-7 
\\
\\ &{\left. \aligned  \undersetbrace  {N-2} \to {\underset \to {\yng{6,5}}} 
\boxtimes \yng{2,2,2,1} \endaligned \right \}  {\ssize 
\aligned l \ \\
 + \\ 1 \
\endaligned } }
 \to \
\yng{6,4}
\boxtimes \yng{2,2,2,1,1}  \ \ 
\to \ \ 
\yng{6,3}
\boxtimes \yng{2,2,2,1,1,1} 
 \to  \dots 
\\ & \hskip .7in \oplus \hskip  .8 in  \nearrow  \hskip .8in \oplus \hskip .6
in  \nearrow 
\\ &{\left. \aligned  \undersetbrace  {N-1} \to {\underset \to {\yng{7,4}}} 
\boxtimes \yng{2,2,1,1,1} \endaligned \right \}  {\ssize 
\aligned l-1 \\
 + \ \\ 3 \
\endaligned }
 }\ 
\to  \yng{7,3} \boxtimes \yng{2,2,1,1,1,1}  \hskip .4 in \dots 
\\ & \hskip .7in \oplus  \hskip  .8 in  \nearrow \hskip .8in \oplus 
\\ &{\left. \aligned  \undersetbrace  {N} \to {\underset \to {\yng{8,3}}} 
\boxtimes \yng{2,1,1,1,1,1} \endaligned \right \}  {\ssize 
\aligned l - 2\\
 + \ \\ 5 \
\endaligned } }\  \hskip .4in \dots
\\ & \dots
\endalignat
$$

Then the lower rows are  the sequence
$A^{l-1}[-2]$, thus using Lemma 3.3 we reduce the proof to showing that the top
row  is exact.

Let us denote the top row by $T^l$.

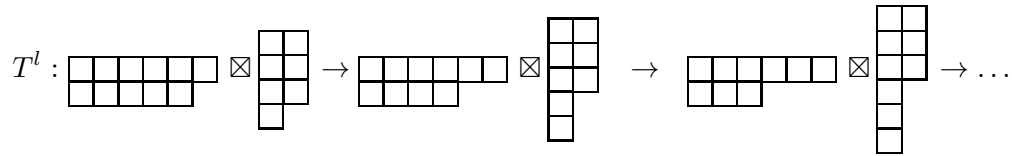
$$
\alignat 9 T^{l}: \ & \yng{6,5} 
\boxtimes \yng{2,2,2,1} 
\to 
\yng{6,4}
\boxtimes \yng{2,2,2,1,1}  \ \ 
\to \ \ 
\yng{6,3}
\boxtimes \yng{2,2,2,1,1,1} 
 \to  \dots 
\endalignat
$$

Use $\Cal V = \Cal O_\Bbb P ((N-2)-l) \otimes 
\overset N-2 \to \Lambda (V/W)^*$, then by 4.5 the following is exact:

$$
\gather {\left. \aligned  \undersetbrace  {N-2} \to {\underset \to {\yng{6,6}}} 
\boxtimes \yng{1,1,1} \endaligned \right \}  {\ssize  l } }
 \to {\left. \aligned  \undersetbrace  {N-2} \to {\underset \to {\yng{6,5}}} 
\boxtimes \yng{1,1,1,1} \endaligned \right \}  {\ssize 
\aligned l \ \\
 + \\ 1 \
\endaligned } }
  \to  \dots \to  {\left. \aligned  \undersetbrace  {(N-2-l) + l} \to {\underset
\to {\yng{6,3}}} 
\boxtimes \yng{1,1,1,1,1,1} \endaligned \right \}  {\ssize  N-2 } } \\
 \to 0 .
\endgather
$$

Tensor this by $\overset l \to \Lambda (V/W)^* $ on the right and call the
corresponding sequence $(**)$, which again is exact.

$$
\alignat 9 &{\left. \aligned  \undersetbrace  {N-2} \to {\underset \to
{\yng{6,6}}} 
\boxtimes \yng{2,2,2} \endaligned \right \}  {\ssize  l  } }
 \to \
\yng{6,5}
\boxtimes \yng{2,2,2,1}  \ \ 
\to \ \ 
\yng{6,4}
\boxtimes \yng{2,2,2,1,1} 
 \to  \dots 
\\ & \hskip .7in \oplus \hskip  .8 in  \nearrow  \hskip .8in \oplus \hskip .6
in  \nearrow 
\\ &{\left. \aligned  \undersetbrace  {N-1} \to {\underset \to {\yng{7,5}}} 
\boxtimes \yng{2,2,1,1} \endaligned \right \}  {\ssize 
\aligned l-1 \\
 + \ \\ 2 \
\endaligned }
 }\ 
\to  \yng{7,4} \boxtimes \yng{2,2,1,1,1}  \hskip .4 in \dots 
\\ & \hskip .7in \oplus  \hskip  .8 in  \nearrow \hskip .8in \oplus 
\\ &{\left. \aligned  \undersetbrace  {N} \to {\underset \to {\yng{8,4}}} 
\boxtimes \yng{2,1,1,1,1} \endaligned \right \}  {\ssize 
\aligned l - 2\\
 + \ \\ 4 \
\endaligned } }\  \hskip .4in \dots
\\ & \dots
\endalignat
$$

Note that the top row is equal to $T^l$. Apply again Lemma 3.3: To show that
$T^l$ and thus $A^l$ is exact up to degree
$2N-6$, it is enough to show that the quotient of $(**)$ by $T^l$ ( the bottom
rows of (**)) is
 exact up to degree $2N-5 $.

Consider the sequence for
$\Cal V = \Cal O_\Bbb P ((N-2)-(l+1)) \otimes 
\overset N-2 \to \Lambda (V/W)^*$:

\newpage

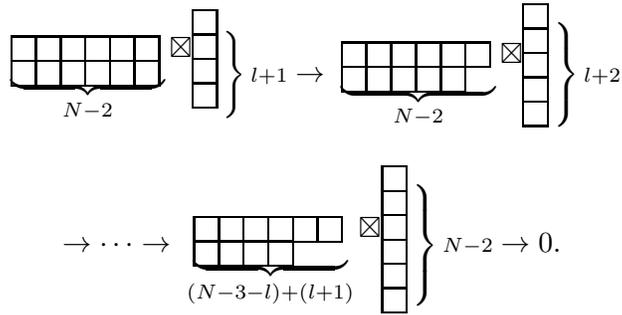
$$
\gather {\left. \aligned  \undersetbrace  {N-2} \to {\underset \to {\yng{6,6}}} 
\boxtimes \yng{1,1,1,1} \endaligned \right \}  {\ssize  l+1 } }
 \to {\left. \aligned  \undersetbrace  {N-2} \to {\underset \to {\yng{6,5}}} 
\boxtimes \yng{1,1,1,1,1} \endaligned \right \}  {\ssize  l+2 } }
\\
\\
 \to  \dots \to  {\left. \aligned  \undersetbrace  {(N-3-l) + (l+1)} \to
{\underset \to {\yng{6,4}}} 
\boxtimes \yng{1,1,1,1,1,1} \endaligned \right \}  {\ssize  N-2 } } 
 \to 0 .
\endgather
$$

Twist this by $\overset l-1 \to \Lambda (V/W)^* $ on the right:

$$
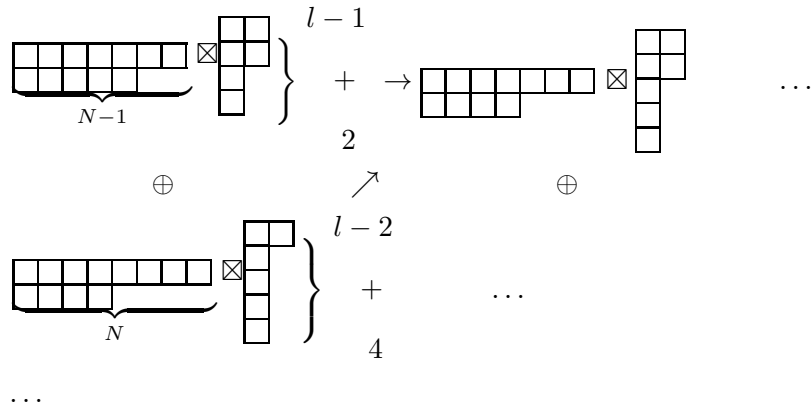
\alignat 8 &{\left. \aligned  \undersetbrace  {N-1} \to {\underset \to
{\yng{7,5}}} 
\boxtimes \yng{2,2,1,1} \endaligned \right \}  {\ssize 
\aligned l-1 \\
 + \ \\ 2 \
\endaligned }
 }\ 
\to  \yng{7,4} \boxtimes \yng{2,2,1,1,1}  \hskip .4 in \dots 
\\ & \hskip .7in \oplus  \hskip  .8 in  \nearrow \hskip .8in \oplus 
\\ &{\left. \aligned  \undersetbrace  {N} \to {\underset \to {\yng{8,4}}} 
\boxtimes \yng{2,1,1,1,1} \endaligned \right \}  {\ssize 
\aligned l - 2\\
 + \ \\ 4 \
\endaligned } }\  \hskip .4in \dots
\\ & \dots
\endalignat
$$

These are just the bottom rows of (**), thus the quotient of $(**)$ by $T^l$.
This completes the proof of the Lemma and thereby the proof of the Main Theorem.
q.e.d.

\enddemo

\specialhead X=SpGr(2,6):
\endspecialhead

Last let us look at one example, $X=SpGr(2,6)$, the symplectic grassmannian of
2-planes in 6-space.

Consider the Tate construction in this case:

$$
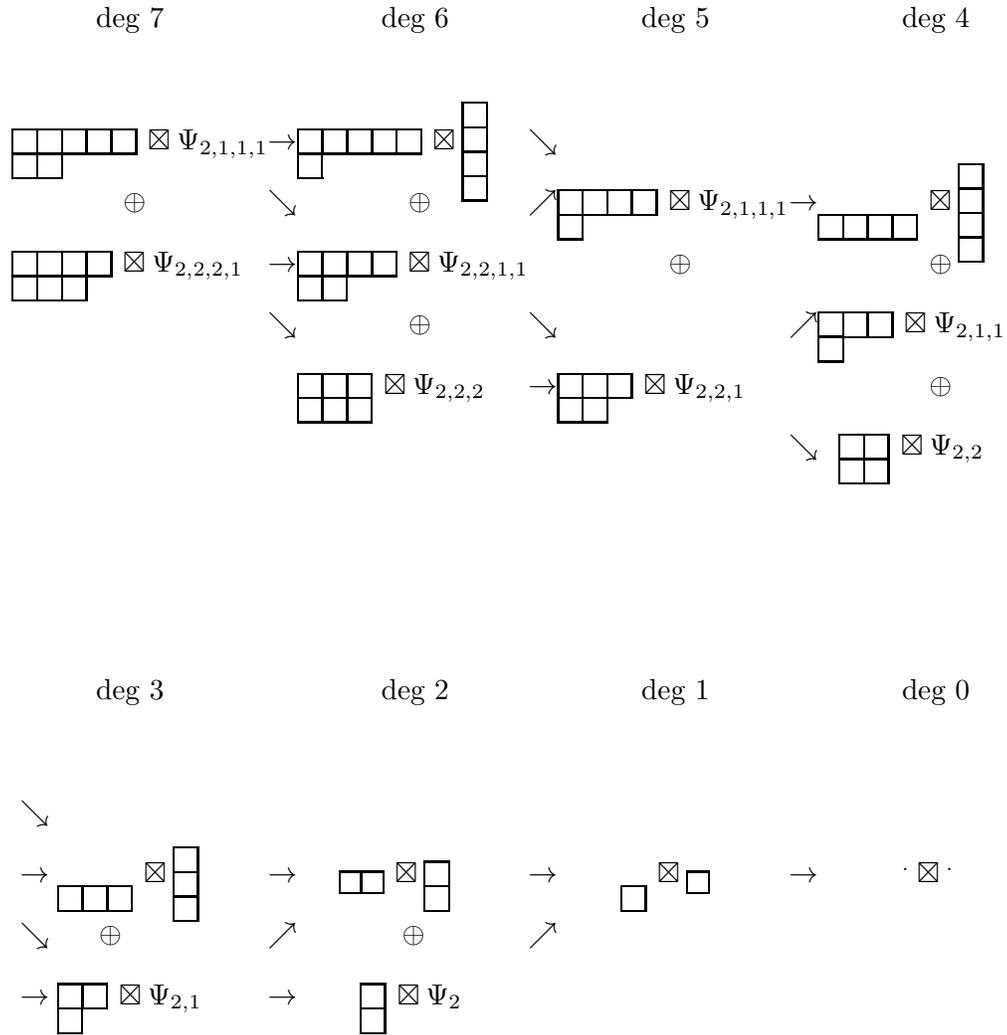
\alignat 8 &\hskip .4in \text{deg 7} &&\hskip .4in \text{deg 6} &&\hskip .4in
\text{deg 5}  && \hskip .4in \text{deg 4}
\\
\\ &\yng{5,2} \boxtimes \Psi_{2,1,1,1} &\to  &\yng{5,1} \boxtimes \yng{1,1,1,1}
&\searrow 
\\ & \hskip .53in \oplus & \searrow & \hskip .53in \oplus & \nearrow &\yng{4,1}
\boxtimes \Psi_{2,1,1,1} &\to &\yng{4}  \boxtimes \yng{1,1,1,1} &  
\\ &\yng{4,3} \boxtimes \Psi_{2,2,2,1} &\to &\yng{4,2} \boxtimes \Psi_{2,2,1,1}
& & \hskip .53in \oplus & & \hskip .53in \oplus & 
\\ &  & \searrow & \hskip .53in \oplus &   \searrow & & 
\nearrow & \yng{3,1} \boxtimes \Psi_{2,1,1} & 
\\ &&& \yng{3,3} \boxtimes \Psi_{2,2,2} & \to &
\yng{3,2} \boxtimes \Psi_{2,2,1} &  & \hskip .53in \oplus &
\\ &&&&&& \searrow & \hskip .1in \yng{2,2} \boxtimes \Psi_{2,2}
\\
\\
\\
\\ &\hskip .4in \text{deg 3} &&\hskip .4in \text{deg 2} &&\hskip .4in \text{deg
1}  && \hskip .4in \text{deg 0}
\\
\\ & \searrow
\\ & \to \yng {3} \boxtimes \yng{1,1,1} & \to  & \hskip .2in \young{2} \boxtimes
\yng{1,1} & \to & \hskip .3 in \yng {1} \boxtimes \young{1} & \to & \hskip .4
in{}^. \boxtimes {}^.
\\ &\searrow \hskip .2 in \oplus & \nearrow & \hskip .5 in \oplus &\nearrow
\\ & \to \yng{2,1} \boxtimes \Psi_{2,1} & \to & \hskip .3in \yng{1,1} \boxtimes
\Psi_2 
\endalignat
$$

\addline
\addline

The sub-complex starts in degree 6, with the term 
$$
\gather
\yng{3,3} \boxtimes \Psi_{2,2,2} \cong \Cal O_X(-3) \boxtimes \Psi_{2,2,2} \ .
\endgather
$$ Thus Theorem 6.1 implies the exactness of the  sub-complex $B.$ :
$$ B. : B_6 \to B_5 \to ... \to B_1 \to \Cal O_{X \times X} \ ,
$$ that is:

\newpage

\proclaim{Corollary 6.5}
$$
\alignat 8 & &&& \nearrow & (Sym_2 S) (-1) S \boxtimes \Psi_{2,1,1} &\to & Sym_3
S \boxtimes \overset 3 \to \Lambda S^\perp & 
\\ & \Cal O_X (-3) \boxtimes \Psi_{2,2,2} & \to & S(-2) \boxtimes \Psi_{2,2,1} &
& \hskip .4in \oplus & & \hskip .4in \oplus  
\\ &&&&\searrow & \hskip .1in \Cal O_X (-2)  \boxtimes \Psi_{2,2} & \to &  S(-1)
\boxtimes \Psi_{2,1} 
\\
\\ &&\to & Sym_2 S \boxtimes \overset 2 \to \Lambda S^\perp   &\ \searrow 
\\ &&\nearrow & \hskip .3in\oplus &&\hskip .2in S \boxtimes S^\perp &\to &
\Cal O_X \boxtimes \Cal O_X 
\\
\\ &&\to &\ \ \overset 2 \to \Lambda S \boxtimes \Psi_{2} &\ \nearrow  
\endalignat
$$ is exact.

\endproclaim

\newpage
\topmatter
\title Appendix: Tensor Product of two Young Diagrams  \endtitle
\endtopmatter

\specialhead  The Young symmetrizer
\endspecialhead

Let $E$ be a vector space of dimension $M$ over an algebraically closed field
$\Bbbk$ of characteristic $0$ and denote by
$\frak S_d$ the symmetric group. 

The symmetric group $\frak S_d$ acts on $E^{\otimes d}$ by permuting the
factors, that is for $\sigma \in \frak S_d$,
$$
(e_1  \otimes e_2 \otimes \dots \otimes e_d ) \cdot \sigma
=e_{\sigma(1)}  \otimes e_{\sigma(2)} \otimes \dots \otimes e_{\sigma(d)} \ .
$$
This action commutes with the left action of $Gl (E)$ .

Let $\lambda=(\lambda_1, \lambda_2, \dots, \lambda_M)$ be  an ordered
partition of length $d$ and $\lambda_M \geq 0$, that is
$\lambda_1 \geq	 \lambda_2 \geq \dots \geq \lambda_M $ and
$\lambda_1 +	 \lambda_2 + \dots + \lambda_M = d $. To a partition
$\lambda$ associate the Young diagram with $\lambda_i$ boxes in the 
$i$-th row, the rows of boxes are lined up on the left.
Denote by $\lambda^*$ the partition defined by interchanging rows and columns
in the Young diagram $\lambda$.

Define two subgroups of the symmetric group:
$$
\gather
P=P_\lambda=\{ p \in \frak S_d | p \ \text {preserves each row} \} \\
Q=Q_\lambda=\{ q \in \frak S_d | q \ \text {preserves each column} \}
\endgather
$$
Set
$$
\gather
a_\lambda=\underset p \in P \to \Sigma e_p \ \text{and} \
b_\lambda=\underset q \in Q \to \Sigma sgn(q) \cdot e_q  .
\endgather
$$

\proclaim{Definition}

\roster
\item
Define the Young symmetrizer of $\lambda,
 c_\lambda$ as
$$
 c_\lambda
=a_\lambda \cdot b_\lambda \in \Bbb C
\frak S_d \ \text{and}
$$ 
\item
denote by $\Sigma^\lambda E$ the image of
$$
c_\lambda: E^{\otimes d} \to E^{\otimes d} \ .
$$
$ E \mapsto \Sigma^\lambda E $ is called the Schur functor of $\lambda$.
\endroster
\endproclaim

\subhead {Examples}
\endsubhead
\roster
\item
For $\lambda=(2,0, \dots 0)$ and its associated Young
diagram 
$
\young{2}
$,
the Young symmetrizer is equal to 
$c_\lambda= 1 + e_{(1 \ 2)} $. Then the image of 
$c_\lambda$ is the subspace of $E \otimes E$ spanned by all vectors:
$$
v_1 \otimes v_2 + v_2 \otimes v_1.  
$$
Thus $\Sigma^{(2,0 \dots 0 )} \cong Sym_2 E$.

\item
For $\lambda=(1,1,0 \dots 0)$, the associated Young diagram is equal to
$$
\gather \yng{1,1} \ \text{and} \
c_\lambda= 1 - e_{(1 \ 2 )} \ . \text{ The image of }  \ c_\lambda
\ \text{ is spanned by all vectors} \\
v_1 \otimes v_2 - v_2 \otimes v_1. 
\endgather
$$
Thus  $\Sigma^{(1,1,0 \dots 0 )} \cong \overset 2 \to \Lambda E$.

\item
For $\lambda=(2,1,0 \dots 0)$, the associated Young diagram is equal to
$$
\gather \yng{2,1} \ \text{and} \ 
c_\lambda= 1 + e_{(1 \ 2 )} - e_{(1 \ 3 )} - e_{(1 \ 3 \ 2 )}\ . 
\endgather
$$ The
image of $ \ c_\lambda $ is the subspace of $E^{\otimes 3}$
 spanned by all vectors
$$
 v_1 \otimes v_2 \otimes v_3
+  v_2 \otimes v_1 \otimes v_3  
-  v_3 \otimes v_2 \otimes v_1
-  v_3 \otimes v_1 \otimes v_2
$$ 

\endroster

\proclaim{Lemma}
\roster
\item
$\Sigma^{(m,0 \dots 0)} E \cong Sym_m E$,
\item
$\Sigma^{(1,\dots , 1,0 \dots 0)} E \cong \overset m \to \Lambda E$.
\endroster
\endproclaim

\demo{Proof}
This follows immediately from the definition of $c_\lambda$, since
in (1)
$c_\lambda=a_\lambda$,
and
in (2)
$c_\lambda=b_\lambda$.

q.e.d.

\enddemo

Note that the definition of the Young-symmetrizer $c_\lambda$ 
and the Schur functor $\Sigma^\lambda E$ also make  sense
for an ordered partition
 $\lambda=(\lambda_1, \lambda_2 \dots, \lambda_M, \lambda_{M+1} \dots )$,
but if $\lambda_{M+1} \neq 0$, then  $\Sigma^\lambda E=0$ .

\specialhead
Littlewood-Richardson Rule
\endspecialhead

Let $\lambda$ and $\mu$ be ordered partitions. Then
$$
\Sigma^\lambda E \otimes \Sigma^\mu E \cong \bigoplus N_{\lambda \ \mu \ \nu }
\Sigma^\nu E \ ,
$$
where $N_{\lambda \ \mu \ \nu } $ can be determined by the
Littlewood-Richardson rule.

\proclaim { Littlewood-Richardson Rule}

$N_{\lambda \ \mu \ \nu }= $ number of ways the Young diagram can be
extended to the Young diagram $\nu$ by a strict $\mu$-expansion. 

A 
$\mu$-expansion of a Young diagram is obtained by first adding $\mu_1$ boxes
to the rows of the Young diagram $\lambda$, but with no two boxes in the
same column and in such a way that one obtains a Young diagram for each box
added. Then put the integer $1$ in each of these $\mu_1$ boxes.
Then adding similarly $\mu_2$ boxes, continuing until finally 
$\mu_M$ boxes are added with the integer $M$. 

The expansion is called strict
if, when the integers in the boxes are listed from right to left, starting
with the top row and working down, and one looks at the first $t$ entries in
this list ( for any $t$ between $1$ and $ \mu_1 + \dots \mu_k$ ), each
integer $p$ between $1$ and $k-1$ occurs as least as many times as the next
integer
$p+1$. 
\endproclaim

\demo{Proof}
See [FH], page 456, and [M], \S 1.9.
\enddemo

\proclaim{Cauchy-Formula}
Let $E$ and $F$ be two vector spaces over $\Bbbk$. Then
$$
\overset p \to \Lambda ( E \otimes F ) \cong
\underset |\lambda|=p \to \bigoplus \Sigma^\lambda E \otimes
\Sigma^{\lambda^*} F 
$$
as $Gl(E) \times Gl(F)$ representations.
\endproclaim

\demo{Proof}
The proof is essentially an application of the Littlewood-Richardson rule, see
[KI], Lemma 0.5. 
\enddemo

\proclaim{Lemma 1}

Let $E$ be a vector space of dimension $M$ over $\Bbbk$, then
$$
\Sigma^\lambda E \otimes \Sigma^{(1, \dots ,1)} E \cong 
\Sigma^{(\lambda_1 + 1, \dots , \lambda_M +1)} \ .
$$
Particularly if dim $E = 2$, then
$$
\Sigma^\lambda E \otimes \Sigma^{(1 ,1)} E \cong 
\Sigma^{(\lambda_1 + 1,  \lambda_2 +1)} \ .
$$
\endproclaim

\demo{Proof}
The Young diagram of $(1, \dots ,1) $ is of the form
$$
\gather
\yng{1,1}
\\
\vdots
\\
\yng{1,1}
\endgather
$$
Thus 
the tensor product of $\lambda$ by $(1, \dots ,1) $ is given by
adding $M$ boxes to each  $\lambda$ following the combinatorial rule of
the Little-Richardson formula. This can only be done in one way, which is to
add one box to each row, thus
$$
\lambda \otimes (1, \dots ,1) \cong 
(\lambda_1 + 1, \dots , \lambda_M +1) \ .
$$
q.e.d.
\enddemo

Set $\Sigma^\mu=0$, if $\mu$ is not ordered, that is
$\mu_i < \mu_{i+1}$ for some $i$.

\proclaim{Lemma 2}
Let $2 \geq \lambda_1 \geq \dots \geq 0$ and set
$\alpha=\lambda^*$, then
$$
\gather
\Sigma^{\alpha^*} \otimes  E  \\ \cong 
\Sigma^{(\alpha_1+1, \alpha_2 , 0 \dots 0 )^*} E   \oplus
\Sigma^{(\alpha_1, \alpha_2+1, 0 \dots 0  )^*} E \oplus
\Sigma^{(\alpha_1, \alpha_2 ,1, 0 \dots 0  )^*} E  \tag 1
\endgather
$$
and
$$
\gather
\Sigma^{\alpha^*} \otimes Sym_2 E   \cong \tag 2 \\
\Sigma^{(\alpha_1+1, \alpha_2 +1, 0 \dots 0 )^*} E  
 \oplus
\Sigma^{(\alpha_1+1, \alpha_2, 1, 0 \dots 0  )^*} E \oplus \\
\Sigma^{(\alpha_1, \alpha_2 ,1,1, 0 \dots 0  )^*} E \oplus
\Sigma^{(\alpha_1, \alpha_2+1,1, 0 \dots 0  )^*}  E 
\endgather
$$
with each summand occurring with multiplicity one.

\endproclaim 

\demo{Proof}

\item{(1)}
Consider the tensor product of $\lambda$ by $\young{1}$. 

This is given by
adding one box to $\lambda$. The box can be added either to the
first, second or third column of $\lambda$. Since
$\alpha_i= $ number of boxes in the $i$-th column of $\lambda$, it follows
that adding a box in the $j$-th column just adds one to $\alpha_j$.
therefore the summands that occur in the tensor product are
$$
(\alpha_1+1, \alpha_2 , 0 \dots 0 )^* \oplus
(\alpha_1, \alpha_2+1, 0 \dots 0  )^* \oplus
(\alpha_1, \alpha_2 ,1, 0 \dots 0  )^* , 
$$
all occurring precisely once.

\item{(2)}
Consider the tensor product of $\lambda$ by $\young{2}$.

This is given by adding
the two boxes of $\young{2}$ to $\lambda$
according to the combinatorial rules, that is at most one box to each column.
This adds either (a) no,(b) one or (c) two boxes to the first two columns of
the Young diagram
$\lambda$.
\item{(a)}
If we add no box to the first two columns, this can only be done by adding
$2$ boxes to the first row of $\lambda$. The thus obtained weight is equal to
$(\alpha_1, \alpha_2,1,1, 0 \dots 0  )^*$.
\item{(b)}
Adding one box to the third column and adding a box to the first column or
second column of $\lambda$ results in the weights
$(\alpha_1+1, \alpha_2, 1, 0 \dots 0  )^* $ and $(\alpha_1, \alpha_2+1, 1, 0
\dots 0  )^* $, which both just occur once.
\item{(c)}
If both boxes are added to the first two columns, then there is only one way
of doing so and the resulting weight is equal to
$(\alpha_1+1, \alpha_2+1,  0 \dots 0  )^*  $.

q.e.d.

\enddemo

\proclaim{Remark} 
\roster
\item
Denote by $det(E)$ the vector bundle $\overset M \to \Lambda
E \cong 
\Sigma^{(1, \dots ,1)} E$. Let $\lambda$ be any ordered partition
$(\lambda_1, \dots , \lambda_M)$. Suppose that the $\lambda_i$'s are not all
positive. Set
$$
\Sigma^\lambda E = \Sigma^{(\lambda_1 - \lambda_M, \dots ,\lambda_M -
\lambda_M) } \otimes det (E) ^{\otimes  \lambda_M } \ ,
$$ then this is well defined by   Lemma 1 and because
$ \lambda_1 - \lambda_M \geq \dots \geq \lambda_M - 
\lambda_M \geq 0$.

\item
The dual of the representation $\Sigma^\lambda E $ is given by:
$$
(\Sigma^\lambda E)^* \cong \Sigma^\lambda E^* \cong 
\Sigma^{(-\lambda_M, -\lambda_{M-1}, \dots , -\lambda_1)} \ ,
$$
see [KI],(0.1).

\item
The maximum number of rows or of columns in each Young diagram
$\nu$ appearing in the tensor product $\lambda \otimes \mu$ does not exceed
the sum of these numbers for the separate factors, that is:

 If $\Sigma^\nu
\subset
\Sigma^\lambda
\otimes
\Sigma^\mu$, then
$\nu_i \leq \lambda_1 + \mu_1$ for all $i$ and 
$(\nu^*)_i \leq (\lambda^*)_1 + (\mu^*)_1 $.

This follows from the Littlewood-Richardson rule, since we add at most
$\mu_1$ boxes to $\lambda_1$ and since we add at most $(\mu^*)_1 $
boxes to the first column of $\lambda$.

\endroster

\endproclaim

\newpage

\topmatter
\rightheadtext{References}
\endtopmatter

\enddocument